\documentclass[12pt]{article}
\usepackage[left=3cm,right=3cm,top=2cm,bottom=2cm]{geometry}
\usepackage[french]{babel}

\usepackage{amsmath, amsfonts, amssymb,latexsym, bm}
\usepackage{wasysym,graphicx,stmaryrd,tikz, tipa, bclogo}
\usepackage{pgfplots}

\usetikzlibrary{arrows.meta}
\tikzset{>={Latex[width=2mm,length=2mm]}}

\addto\captionsfrench{}

\usepackage{amsmath, amsfonts, amssymb,latexsym}

\input epsf.sty

\newtheorem{Lemma}{Lemma}[section]

\newtheorem{Corollary}[Lemma]{Corollary}
\newtheorem{Definition}[Lemma]{Definition}

\newcommand{\BEQ}{\begin{equation}}     
\newcommand{\BEA}{\begin{eqnarray}}
\newcommand{\BD}{\begin{displaymath}}
\newcommand{\EEQ}{\end{equation}}       
\newcommand{\EEA}{\end{eqnarray}}
\newcommand{\ED}{\end{displaymath}}

\newcommand{\del}{\delta}
\newcommand{\Del}{\Delta}
\newcommand{\eps}{\varepsilon}          

\newcommand{\Tr}{{\mathrm{Tr}}}

\newcommand{\R}{\mathbb{R}}

\newcommand{\Z}{\mathbb{Z}}
\newcommand{\N}{\mathbb{N}}

\newcommand{\Id}{{\mathrm{Id}}}

\def\T{{\mathbb{T}}}

\newcommand{\eop}{\hfill $\Box$}        
\newcommand{\Medskip}{\medskip\noindent}
\newcommand{\Bigskip}{\bigskip\noindent}

\renewcommand{\Re}{{\rm Re\ }}          
\newcommand{\half}{{1\over 2}}          
\newcommand{\diag}{{\mathrm{diag}}}

\catcode`\@=11
\def\numberbysection{\@addtoreset{equation}{section}
        \def\theequation{\thesection.\arabic{equation}}}
\numberbysection

\begin{document}
\renewcommand{\contentsname}{Table of contents}
\renewcommand{\bibname}{References}

\title{\bf General multi-scale estimates for Lyapunov data  of Perron-Frobenius matrices.
The case of 
diluted autocatalytic chemical reaction networks}
\vskip -2cm
\date{}
\vskip -2cm
\maketitle

\vspace{2mm}
\begin{center}
{\bf  J\'er\'emie Unterberger}
\end{center}

\vspace{2mm}
\begin{quote}

\renewcommand{\baselinestretch}{1.0}
\footnotesize
{Autocatalytic chemical reaction networks are dynamical systems  whose linearization around zero
$dX/dt=AX$ is represented by a Perron-Frobenius matrix $A$
with positive Lyapunov exponent ; this exponent gives the growth rate of the species concentration vector $X$  in the diluted regime, i.e. in a  vicinity of zero.
We introduce here a new, general recursive procedure
providing precise quantitative information about Lyapunov data, namely, the Lyapunov eigenvalue,
and left and right eigenvectors. Our estimates are
based on a multi-scale  algorithm inspired from Wilson’s renormalization
group method in quantum field theory, and Markov chain arguments introduced in Nghe \&
Unterberger \cite{NgheUnt1}. They are  compatible with the very scarce knowledge of
kinetic rates (coefficients of $A$) generally available in chemistry, and take on the form of simple rational functions
of the latter.

 }

\end{quote}
\vspace{4mm}
\noindent



\tableofcontents


\section{Introduction} 


Perron-Frobenius matrices are square matrices $A= (A_{i,j})_{1\le i,j\le n}$ such that 
off-diagonal coefficients are non-negative, $A_{i,j}\ge 0$ for all $i\not=j$. Such matrices arise
in particular in a probabilistic context. Namely, fix $\Sigma\simeq \{1,\ldots,n\}$; 
continuous-time
Markov chain generators -- and  discrete-time Markov matrices alike -- on state space $\Sigma$ are Perron-Frobenius matrices. Replacing nonzero off-diagonal coefficients of a  Perron-Frobenius
matrix $A$ by $1$, we get the adjacency matrix of a directed graph $G=(\Sigma,E)$; by definition, $(\sigma,\sigma')$ is  an edge of $G$ if $A_{\sigma',\sigma}>0$. Then $A$ is   irreducible if and only
if $G$ is strongly connected, i.e. $\forall \sigma,\sigma'\in\Sigma$, there exists a directed path 
of $G$ connecting $\sigma$ to $\sigma'$. It is well-known that irreducible Perron-Frobenius matrices $A$ enjoy the following property:

\begin{itemize}
\item[(i)] there exists a unique eigenvalue $\lambda^*=\lambda^*(A)$ (called: {\em Lyapunov eigenvalue}) such that $\Re\lambda^*= \max 
\{\Re \lambda, \lambda\in {\mathrm{Spec}}(A)\}$, where Spec$(A)$ is the spectrum of $A$;
furthermore, $\lambda^*$ is real;
\item[(ii)] the algebraic multiplicity of $\lambda^*$ is one, and the corresponding eigenvector
may be chosen with strictly positive components.
\end{itemize}

Property (ii) implies unicity (up to normalization) of the eigenvector $v^*$ associated to 
the Lyapunov eigenvalue, called {\em Lyapunov eigenvector}; furthermore, $v^*_{\sigma}>0$ for 
all $\sigma\in\Sigma$.  Since the adjoint matrix $^t A$ 
is also Perron-Frobenius irreducible, there is also a unique (up to normalization)
adjoint (left) eigenvector $v^{\dagger,*}$.  Properties (i) and (ii) together
imply that Spec$(A)=\{\lambda_1,\ldots,\lambda_n\}$, with $\lambda_1=\lambda^*$ and 
$\Re \lambda_i<\lambda^*$, $i=2,\ldots,n$. Also, choosing the normalization of $v^*,v^{\dagger,*}$
in such a way that $\langle v^{\dagger,*},v^*\rangle=1$, we get for every vector $u>0$, i.e. 
$u\in \R^{\Sigma}$ such that $u_{\sigma}>0$ for all $\sigma$:
\BEQ  e^{tA} u \sim_{t\to\infty} \langle v^{\dagger,*},u\rangle \ e^{\lambda^*t} v^*  
\label{eq:lin-asympt}
\EEQ
where $\langle\cdot,\cdot\rangle$ is the canonical scalar product on $\R^{\Sigma}$.     
    
\Medskip We are interested here in estimating the {\em Lyapunov data}, i.e. Lyapunov 
eigenvalue $\lambda^*$, eigenvector $v^*$ and adjoint eigenvector $v^{\dagger,*}$ of a large 
subset $\mathbb A$ of Perron-Frobenius matrices relevant in a particular chemical context which we 
describe below. Let us first describe  $\mathbb A$. An element of $\mathbb A$ is 
of the form $A\equiv \sum_{\rho} A^{\rho}$ (sum over $\rho$ ranging in a finite, abstract set $R$), 
where $A^{\rho}$ can be of the following form:

\begin{itemize}
\item[(i)] either there exist $\sigma\in\Sigma$ and $\beta_{\sigma}>0$ such that  $A^{\rho}_{\sigma',\sigma} = -\beta_{\sigma} \del_{\sigma',\sigma}$, $\sigma'\in\Sigma$;    
\item[(ii)] or there exist $\sigma\not=\sigma'\in\Sigma$ and $k^{\rho}>0$ such that $A^{\rho}_{\sigma',\sigma} = k^{\rho}$ and $A^{\rho}_{\sigma,\sigma}=-k^{\rho}$;
\item[(iii)] or  there exist $\sigma\not=\sigma',\sigma''\in\Sigma$ and $k^{\rho}>0$ such that $A^{\rho}_{\sigma',\sigma}, A^{\rho}_{\sigma'',\sigma} = k^{\rho}$ and $A^{\rho}_{\sigma,\sigma}=-k^{\rho}$
\end{itemize}
(unspecified coefficients are zero). 
The rationale comes from the kinetic modeling of open chemical reaction networks \cite{Fein}, 
\cite{Kamp}, \cite{RaoEsp}.  (Footnotes in the Introduction give minimal information about the relevant chemical context, but may be skipped altogether.)
Namely, assume that $\Sigma$ indexes a finite set of {\em chemical species} related by 
{\em reactions}
$\rho$ of type (i) $\sigma \overset{\beta_{\sigma}}{\to} \emptyset$, 
(ii) $\sigma\overset{k^{\rho}}{\to} \sigma'$ or (iii) $\sigma\overset{k^{\rho}}{\to} \sigma'+\sigma''$
($\sigma\overset{k^{\rho}}{\to} 2\sigma'$ when $\sigma''=\sigma'$), with kinetic rates $k^{\rho}$. 
Species $\sigma$ (left of the arrow) is called {\em reactant}, species $\sigma',\sigma''$ 
(right of the arrow) are called {\em products}\footnote{The reaction network may be represented as 
a hypergraph \cite{Sri}, where edges connect a reactant to its product(s). However, only the graph structure $G$ will be relevant in this work.} (extension of (iii) to reactions with $\ge 3$ 
products is possible but chemically irrelevant). Let $X=(X_{\sigma})_{\sigma\in \Sigma}\ge 0$ be the vector of {\em concentrations} of the chemical
species. Then the time evolution of $X$ is given by 
\BEQ dX/dt = AX := \sum_{\rho\in R} F^{\rho} \label{eq:kin},
\EEQ
where $F^{\rho} = A^{\rho} X$ is {\em linear}, i.e. (depending on the type of $\rho$) 

\begin{itemize}
\item[(i)] $F^{\rho}_{\sigma} = -\beta_{\sigma} X_{\sigma}$;   
\item[(ii)] $F^{\rho}_{\sigma} = -k^{\rho} X_{\sigma}, \qquad F^{\rho}_{\sigma'} = +k^{\rho} X_{\sigma}$;
\item[(iii)] $F^{\rho}_{\sigma} = -k^{\rho} X_{\sigma},\qquad F^{\rho}_{\sigma'} = F^{\rho}_{\sigma''} =  +k^{\rho} X_{\sigma}$.
\end{itemize}
Type (i), (ii) dynamics are equivalent to the time evolution of the measure $\mu(t)$ of a continuous-time
Markov chain on $\Sigma$ with  transition rates (ii) $k^{\rho}$ for the {\bf 1-1 reaction} $\rho:\sigma\to\sigma'$,
and killing rates (i) $\beta_{\sigma}$ for the {\bf degradation reaction} $\rho: \sigma\to 
\emptyset$, where $\emptyset$ denotes the cemetery state. Thus, in absence of {\bf 1-2 reactions}, i.e. reactions of type (iii), the kinetic equations 
(\ref{eq:kin}) model the {\em average time-evolution of a deficient Markov chain with adjoint generator 
$A$}. {\em 'Deficient'} refers to the fact that the time-evolution given by $A$ does not preserve in general
the space of probability measures, contrary to standard ({\em 'non-deficient'}) Markov generators.          

\Medskip Type (iii) equations can be simulated by a Gillespie-type algorithm for branching particle systems similar to that
 of Markov chains, the space state being now $\N^{\Sigma}$; namely, each pair (particle
 of type $\sigma$, reaction $\rho$ with reactant $\sigma$) is equipped with a clock with law Exp$(\beta_{\sigma})$
 or Exp$(k^{\rho})$; if the first reaction to ring is of type (iii), then one particle of type $\sigma$ is replaced by one particle of type $\sigma'$ and one of type $\sigma''$. In the kinetic limit
we are interested in (where the number of particles goes to infinity), the law of large numbers
implies the deterministic equation  (\ref{eq:kin}). Although the time-evolution in $\R^{\Sigma}$ is deterministic in the kinetic limit, it helps to think of it as the time-evolution of the measure $\mu(t)$ of a continuous-time
Markov chain on $\Sigma$ with  transition rates (i), (ii) as above, to which, for every
reaction   $\sigma\overset{k^{\rho}}{\to} \sigma'+\sigma''$ of type (iii), the following
 transitions (later on called: {\em split reactions}) must be added,

\BEQ (ii) \   \sigma\overset{k^{\rho}}{\to} \sigma', \ \sigma\overset{k^{\rho}}{\to} \sigma''; 
\qquad (i)\  \sigma\overset{-k^{\rho}}{\to} \emptyset  \label{eq:1.3}
\EEQ
with now an equivalent {\em negative} killing rate $-k^{\rho}$. Edges $\sigma\to\sigma',\sigma\to\sigma''$ are part of the adjacency graph $G$. 
In absence of degradation reactions, and provided $A$ is irreducible, it was proved in \cite{NgheUnt1} that the Lyapunov 
exponent of $A$ is $>0$ beyond the standard, textbook Markov case, i.e. in presence of at least one
branching reaction $\sigma\to\sigma'+\sigma''$. For chemistry, this means that any
initial positive perturbation (i.e. $X(0)>0$)  is exponentially amplified through (\ref{eq:lin-asympt}) \footnote{In reality, the presence of reverse transitions $\bar{\rho}:\ \sigma'+\sigma''\to\sigma$ with
quadratic (so-called mass-action) rates $F^{\bar{\rho}}$, negligible in a small vicinity of zero but not any
more away from zero, makes the whole picture much more complicated. However, the linear behavior described here is experimentally observable in the {\em diluted regime} \cite{NgheUnt1}, i.e. starting
from small enough $X(0)$.}. This behavior, called  
{\bf autocatalytic}, has been observed experimentally, the foremost example being that of the
celebrated formose reaction network discovered by Butlerov
\cite{But} in 1861, and is currently under close examination by physicists, chemists and biologists interested
in the {\em origin of life} problem, i.e. in primitive forms of evolution observed in 
prebiotic (non-living but organic) chemistry, for its analogy with replication \cite{SmiMor}, 
\cite{SmithSzath}. The standard
definition of autocatalysis, see \cite{BloLacNgh}, is stoechiometric, i.e. based on the 
reaction set $R$, but independent of the choice of the kinetic rates (it simply states that there exists some linear combination of the reactions such that the total balance is strictly positive for all species, a criterion which disregards in particular degradation rates). However, it has been recognized \cite{PasPro} that kinetics 
(i.e. the magnitude of the coefficients of $A$) play an essential role in general.  Here, we keep to our dynamical definition in terms of positivity of the Lyapunov exponent of  
$A$.   We exhibit in \cite{NgheUnt1} a necessary and sufficient topological condition for stoechiometric autocatalysis to hold, which is however not sufficient for dynamical autocatalysis to hold, except in the case when $A$ irreducible (equivalently, $G$ is strongly connected) and all 
degradation rates vanish.    Related works \cite{DesLacUnt1}, \cite{DesLacUnt2} study on chemically motivated examples the {\em autocatalysis threshold}, namely, the maximum value of
degradation rates for which autocatalysis holds. The current study arose from a search for general
quantitative formulas for $\lambda^*$, in particular, for autocatalysis thresholds.


\subsection{Aim of the article} \label{subsection:aim}


\Medskip A typical situation encountered in prebiotic chemistry, when dealing with {\em complex
mixtures}, i.e. with large reaction networks ($n=|\Sigma|$ is typically of order $10^2-10^4$ or
even larger), is that the set of relevant reactions $R$  -- in particular, the adjacency matrix --
is by and large known (chemical expertise has made it possible to 'catalog' possible reactions
according to a large but finite set of observed rules), but very little is known about kinetic rates $k^{\rho}$, $\rho\in R$. One of the main difficulties (for a mathematician, at least) does not come from the large number of rates, but from the fact that rate magnitudes vary in a wild range, typically $k^{\rho} \sim 10^{-12}-10^{12}$ (the unit is $s^{-1}$), with a turnover 
time $(k^{\rho})^{-1}$ ranging from the picosecond to a year or more \footnote{There is no real upper limit, some reactions having a theoretical time-scale larger than the age of the universe.} (and most observations
relating to times compatible with laboratory experiments, from the minute to the week). Thus
the ODE giving the time-dynamics of the system is {\em stiff}, making it difficult to solve
numerically; further, uncertainty on the rates can be amplified through the equation flow, yielding unreliable predictions \cite{Grassi}.

\Bigskip We turn this difficulty into an asset here by assuming  {\em separation of scales}.
Namely, we assume that {\em reaction rates do not 'accumulate' around any
given kinetic scale.} Following Eyring's equation\footnote{The celebrated Eyring equation from
transition state theory suggests the Ansatz $k^{\rho}\sim {\cal A}\, e^{-\Del 
G^{\ddagger}_{\rho} /RT} \sim {\cal A}\,  10^{-\Del G^{\ddagger}_{\rho}/RT\ln(10)}$,
 where ${\cal A}$ (Arrhenius constant) is a constant prefactor, $\Del G^{\ddagger}_{\rho}$
is the activation energy of the reaction (in kcal/mol), measured up to chemical precision, and $RT\ln(10) \simeq 1.36$ kcal/mol at room temperature $T=298 K$.} and standard chemical precision
conventionally set at $\sim 1$ kcal/mol, we assume 
kinetic constants to be determined in the best case up to a multiplicative factor $b\sim 10$.  
 This suggests that we define  the {\em scale} $n$ of a (reaction with)
kinetic rate $k$  to be the integer part 
(or some rounding) of the base $b$ logarithm
of $k$, 
\BEQ n = \lfloor \log_{b} (k) \rfloor \label{eq:scale}
\EEQ
Thus, knowing $n=n(k)$ yields an indetermination of $k$ comparable
to that obtained when chemical precision is reached.
The choice of sign in (\ref{eq:scale}) implies that 
the larger $n$, the faster the reaction. The parameter $b>1$ will be called   {\bf scale parameter}.  Separation of 
scales  means then that {\em for each integer $n$, there are 'not too many'
reactions with kinetic rate $k$ such that $\lfloor \log_{b} (k) \rfloor =n$.}  A somewhat
more formal
statement would be that our estimates deteriorate, slowly but in general
irremediably, with the number of reactions with scale $n$ in a 
given {\em cluster} (see \S \ref{subsection:multi-scale-method}). 
There are some further restrictions, namely {\em resonance cases}
when scales associated with growth rates of 
subnetworks happen to coincide with their external scales, namely, with
the rate scales of outgoing reactions; close to resonances 
(which are hypersurfaces of codimension one, i.e. defined by one equation, in the space of rate constants), growth rates and vectors shift
rapidly, and therefore, depend on the fine-tuning of the rate
parameters. Another resonance case is when there appear several {\em autocatalytic cores}, i.e.
autocatalytic clusters whose Lyapunov eigenvalues both have maximal order of magnitude.    Generally speaking, we shall use the word {\em generically}
when it is assumed that the network is far from resonance regimes.

\Bigskip Hypothesizing separation of scales, and leaving out the
delicate discussion of resonances, our results may be put bluntly
in this way.  {\em Assume the reaction network   is known, and
so are all reaction scales $n=n(k)$. Then we describe an explicit
algorithm\footnote{The algorithm has been implemented in Python, and its predictions validated
on half a dozen non trivial examples.}   giving  the Lyapunov data, at or close
to chemical precision, i.e. the error on $\log_b (\lambda^*) , \log_b (v^*_{\sigma}),
\log_b (v^{\dagger,*}_{\sigma}), \ \sigma\in\Sigma$ is of order $O(1)$.} The final outcome
of the algorithm is streamlined in the key equations (\ref{eq:pisigma1})--(\ref{eq:vdaggersigma}), which we call {\bf hierarchical formulas}. 
 
\noindent Note that this makes no sense if $\lambda^*<0$; our algorithm actually detects automatically if
degradation rates are below threshold, and provided it be the case, estimates $\lambda^*$, which is
de facto $>0$. However, it says nothing about the sign of $\lambda^*$ in a shadow zone where degradation does not dominate. In particular, it fails to detect that $\lambda^*=0$ in the strictly Markov case (all reactions of type (ii)), but does provide estimates for $\log_b (v^*_{\sigma}), \sigma\in\Sigma$, where $v^*$ is the stationary measure, e.g. when $A$ is irreducible, which is valuable when $A$ is not reversible.  
 
 \Bigskip Furthermore, our algorithm can be run
with an even more minimal input, namely, the set of reactions, plus {\em an ordering of the 
scales of all reactions}. The latter is in principle more or less 
folklore knowledge for the chemist, who knows which reactions are 
'fast' (higher scale) and which are 'slow' (lower scale). This may be enough for simple networks; a more detailed grouping of reactions
into several scales is wanted for more complex ones. Now,  resonance loci 
split the space of rate scales into  a finite (but possibly large)
number of {\em kinetic domains}, and a separate expression is found  for Lyapunov exponents/eigenvectors in each domain.  Therefore, one must 
simultaneously write the resonance equations, which define the 
boundaries of the domains. Fortunately, resonance equations are 
 linear in the  scale parameters.  Let $\sigma$, resp. $e$ range over the set of species $\Sigma$, resp. edges $E$. The general outcome is therefore  a {\em decision tree} $\T=(V(\T),E(\T))$ in the form of successive
linear inequalities with integer coefficients in the reaction scale parameters $(n_e)_e$, $\sum_{e} c^{\T,i}_e n_e\gtrless 0$ $(i\in V(\T))$, defining the domains, and obtained inductively by going
up the branches of $\T$ from the root;
and then, {\em  estimates for  Lyapunov data}, in the form of logarithmic
expressions $\log_{b}(\lambda^*)
\sim 
\sum_{e}  c_e n_e$, $\log_{b} (v^*_{\sigma}) \sim \sum_e C_{\sigma,e}
n_e$, $\log_{b} (v^{\dagger,*}_{\sigma}) \sim \sum_e C^{\dagger}_{\sigma,e}
n_e$.  Coefficients $c_e, C_{\sigma,e}, C^{\dagger}_{\sigma,e}$ are domain-dependent integer numbers. Let us mention that, for some 
edges $e$, called {\em secondary edges}, all these coefficients vanish; the associated (secondary) reactions have therefore little 
influence on the growth of the system, and may be discarded.    Remaining,  main reactions are grouped into
{\em hierarchical clusters} (once again, see \S \ref{subsection:multi-scale-method}), within which all reactions are equally important. Overall,
by  simply erasing secondary reactions, and clustering main ones, 
one  wins a much simpler
vision of the network, and a detailed insight into its organization.

\vskip 2cm In the rest of the section, we loosely describe some
of the main mathematical tools, and say some words on the multi-scale
method.


\subsection{Resolvent formula}  \label{subsection:resolvent}


One of the main tools used here (also in \cite{NgheUnt1}) is the {\bf resolvent formula}. Namely,
let $\alpha>\lambda^*$; then all eigenvalues of $A(\alpha):=A-\alpha \Id$ have negative
real part, and the entries of the resolvent $(\alpha\Id-A)^{-1}$ may be computed as 
\BEQ ((\alpha \Id - A)^{-1})_{\sigma',\sigma} =  \Big(\int_0^{+\infty} e^{t(A-\alpha \Id)}\, dt\Big)_{\sigma',\sigma}= \sum_{\gamma:\sigma\to\sigma'} 
\Big(\prod_{k} w(\alpha)_{\sigma_k\to\sigma_{k+1}}\Big) \, \times\, \frac{1}{|A_{\sigma',
\sigma'}|+\alpha}.   \label{eq:resolvent-formula}
\EEQ
where $\gamma$ ranges over the set of all paths $\sigma\equiv \sigma_1 \to \sigma_2\to
\cdots\to \sigma_{\ell+1}\equiv \sigma'$ of arbitrary length $\ell\ge 0$ from $\sigma$ to $\sigma'$, and 
$w(\alpha)_{\sigma\to\sigma'} = \frac{k_{\sigma\to\sigma'}}{|A_{\sigma,\sigma}|+\alpha}$ are
{\bf  transition weights}, equal to the transition probabilities of the skeleton (discrete-time) 
Markov chain associated to $A$ when equivalent killing rates and $\alpha$ vanish, so that 
$|A_{\sigma,\sigma}|= -A_{\sigma,\sigma} = \sum_{\sigma'\not=\sigma} k_{\sigma,\sigma'}$. The
path sum is always defined in $\R_+\cup\{+\infty\}$, but converges if and only if $\alpha>\lambda^*$. The product $w(\alpha)_{\gamma}:= \prod_{k} w(\alpha)_{\sigma_k\to\sigma_{k+1}}$ is interpreted as a {\bf path weight}.   We get in particular the self-consistent equation for $\lambda^*$, 
\BEQ \Phi(\lambda^*)_{\sigma}:=\sum_{\gamma:\sigma\to\sigma} \prod_k w(\lambda^*)_{\sigma_k\to\sigma_{k+1}} = 1, \label{eq:excursion}
\EEQ where
$\gamma$ now ranges over {\em excursions}, i.e. cycles from $\sigma$ to $\sigma$ which do not
contain $\sigma$ as an intermediate vertex. The unrestricted sum $\gamma:\sigma\to\sigma'$ in
(\ref{eq:resolvent-formula}) involves the geometric series $\sum_{n=0}^{+\infty} (\Phi(\alpha)_{\sigma})^n$, which is equal to $\frac{1}{1-\Phi(\alpha)_{\sigma}}$ for $\alpha>\lambda^*$,
but diverges at $\alpha=\lambda^*$. 

\Medskip All formulas in this article are based on transition weights, and often take on the
form of a product of edge weights along a path. Hill \cite{Hill} proved a long time ago that the stationary
measure of an irreducible discrete-time Markov chain could be written as a sum ranging over spanning trees of similar weight products. We know no extension of Hill's formula to Markov chains with nonzero killing rates, but our approximate formulas are in the same spirit. They extract a hierarchical, tree-like structure (see below) encoding the leading contribution in 
(\ref{eq:resolvent-formula}), obtained by a partial resummation over paths.


\subsection{Merging patterns and hierarchical structures} \label{subsection:merging}


The appropriate description of our algorithm is through {\em hierarchical structures}. These are particular 
instances of hierarchical graphs, which can be defined through {\em merging patterns} and 
{\em encapsulated graphs}. A general merging pattern is  a finite sequence of the form $\Sigma\equiv \Sigma^{(0)}\to \Sigma^{(1)}\cdots \to \Sigma^{(j_{max})}$, where elements of $\Sigma^{(j)}$ are either elements of $\Sigma^{(j-1)}$ or 
subsets of $\Sigma^{(j-1)}$; e.g. (see Fig. \ref{fig:merging}) $\Sigma=\{0,\cdots,8\}$, $\Sigma^{(1)} = \{\{0,1,2\},3,\cdots,8\}$, 
$\Sigma^{(2)} = \{\{0,1,2\}, 3, \{4,5\}, 6,7,8\}$, $\Sigma^{(3)}= \{\{\{0,1,2\},\{3\},\{4,5\}\},6,7,8\}$,\\ $\Sigma^{(4)} = \{\{\{0,1,2\},\{3\},\{4,5\}\},6,\{7,8\}\}$, $\Sigma^{(5)} =
\{\{\{0,1,2\},\{3\},\{4,5\}\},\{6,\{7,8\}\}\}$ and finally, $\Sigma^{(6)} =\{\{\{\{0,1,2\},\{3\},\{4,5\}\},\{6,\{7,8\}\}\}\}$  . The recursive nesting, matryoshka-like, structure, may be drawn as a coalescence tree, see Fig. \ref{fig:merging}.1, with compound vertices given new names, $G_1=\{0,1,2\}, G_2=\{4,5\}, G_3= 
\{G_1,3,G_2\}, G_4 = \{7,8\}, G_5= \{6,G_4\}, G_6=\{G_3,G_5\}$.

\begin{center}
\begin{tikzpicture}[scale=0.6] \label{fig:merging}
\draw(0,0) node {$0$}; \draw(0,-1) node {$1$}; \draw(0,-2) node {$2$}; \draw(0,-3) 
node {$3$}; \draw(0,-4) node {$4$}; 
\draw(0,-5) node {$5$};   \draw(0,-6) 
node {$6$}; \draw(0,-7) node {$7$}; 
\draw(0,-8) node {$8$}; 
\draw(0.2,0)--(2.1,-0.9); \draw(0.2,-1)--(2.1,-1); \draw(0.2,-2)--(2.1,-1.1);
\draw(2.5,-1) node {$G_1$};

\draw(0.2,-4)--(4.6,-4.45); \draw(0.2,-5)--(4.6,-4.55); \draw(5,-4.5) node {$G_2$};

\draw(5.5,-4.5)--(7.5,-3); 
\draw(0.2,-3)--(7.3,-3);
\draw(2.7,-1)--(7.5,-3); 
\draw(8,-3) node {$G_3$};

\draw(0.2,-7)--(10,-7.45); \draw(0.2,-8)--(10,-7.55); \draw(10.5,-7.5) node {$G_4$};

\draw(0,1) node {\small $j\!=\!0$}; \draw(2.5,1) node {\small $j\!=\!1$}; 
\draw(5,1) node {\small $j\!=\!2$}; \draw(8,1) node {\small $j\!=\!3$}; 
\draw(10.5,1) node {\small $j\!=\!4$}; \draw(12.5,1) node {\small $j\!=\!5$}; 
\draw(15,1) node {\small $j\!=\!6$}; 

\draw(0.2,-6)--(12,-6.45); \draw(11,-7.5)--(12,-6.55); \draw(12.5,-6.5) node {$G_5$};

\draw(8.5,-3)--(14.5,-4.95); \draw(13,-6.5)--(14.5,-5.05); \draw(15,-5) node {$G_6$};
\end{tikzpicture}

FIGURE \ref{fig:merging}.1 -- Example of coalescence tree
\end{center}

\bigskip

\begin{center}
\begin{tikzpicture}[scale=1.5] 
\draw(0,0) node {$0$}; \draw[->](0.2,0)--(0.5,0); \draw(0.5,0)--(0.8,0); \draw(1,0) node {$1$}; 
\draw(0.5,-0.7) node {$2$}; \draw[->](1,-0.2)--(0.8,-0.4); \draw(0.8,-0.4)--(0.6,-0.6);
\draw[->](0.4,-0.6)--(0.2,-0.4); \draw(0.2,-0.4)--(0,-0.2);
\draw(-0.2,-1) rectangle(1.2,0.2); \draw(0.5,0.4) node {$G_1$};

\draw(2,-2) node {$3$}; 
\draw[->](0.5,-0.8)--(1.2,-1.3); \draw(1.2,-1.3)--(1.9,-1.8);
\draw[-<](0.5,-0.9)--(1.2,-1.4); \draw(1.2,-1.4)--(1.9,-1.9);
\draw[->](2.1,-1.7)--(2.55,-0.9); \draw(2.55,-0.9)--(3,-0.1);
 \draw[-<](2.1,-1.9)--(2.55,-1.1); \draw(2.55,-1.1)--(3,-0.3);

\draw(3,0) node {$4$}; \draw(4,0) node {$5$}; 
\draw[->](3.2,0.05)--(3.5,0.05); \draw(3.5,0.05)--(3.8,0.05);
\draw[-<](3.2,-0.05)--(3.5,-0.05); \draw(3.5,-0.05)--(3.8,-0.05);
\draw(2.8,-0.5) rectangle(4.2,0.2); \draw(3.5,0.4) node {$G_2$};

\draw[->](4.05,-0.2)--(4.05,-1.9); \draw(4.05,-1.9)--(4.05,-3.8);
\draw[-<](3.95,-0.2)--(3.95,-1.9); \draw(3.95,-1.9)--(3.95,-3.8);

\draw(-0.4,-2.3) rectangle(4.4,0.6); \draw(2,0.8) node {$G_3$};

\draw(4,-4) node {$6$};
\begin{scope}[shift={(2,-3)}]
\draw(3,0) node {$7$}; \draw(4,0) node {$8$}; 
\draw[->](3.2,0.05)--(3.5,0.05); \draw(3.5,0.05)--(3.8,0.05);
\draw[-<](3.2,-0.05)--(3.5,-0.05); \draw(3.5,-0.05)--(3.8,-0.05);
\draw(2.8,-0.5) rectangle(4.2,0.2); \draw(3.5,0.35) node {$G_4$};
\end{scope}
\draw[->](4.1,-3.9)--(4.5,-3.5); \draw(4.5,-3.5)--(4.9,-3.1);
\draw[-<](4.2,-4)--(5.1,-3.6); \draw(5.1,-3.6)--(6,-3.2);
\draw(3.8,-4.2) rectangle(6.4,-2.5);  \draw(5.1,-2.3) node {$G_5$};

\draw[->](0.5,-0.9)--(0.5,-2); \draw(0.5,-2)--(0.5,-4)--(3.9,-4);

\draw(-0.6,-4.4) rectangle(6.6,1); \draw(3,1.2) node {$G_6$};
\end{tikzpicture}

\vskip 0.5cm

FIGURE \ref{fig:merging}.2 -- Possible hierarchical graph associated to the  coalescence tree of Fig. 
\ref{fig:merging}
\end{center}


   We now let elements of $\Sigma^{(j)}$, $j=0,\ldots,j_{max}$ be vertices of a graph $G^{(j)} = (\Sigma^{(j)},E(G^{(j)}))$, which is constructed from $G^{(0)}\equiv G$ and $\Sigma$ by 
   successive rewiring steps as follows: if 
$\sigma\to\sigma'$ is an edge of $G^{(j)}$ $(\sigma,\sigma'\in \Sigma^{(j)})$, then $G^{(j+1)}(\sigma)\to G^{(j+1)}(\sigma')$ is an 
edge of $G^{(j+1)}$, where $G^{(j+1)}(\sigma)$ is the vertex of $G^{(j+1)}$ equal to or containing
$\sigma$. Self-edges $\sigma\to\sigma$ are removed.  If e.g. (see Fig. \ref{fig:merging}.2) $E(G^{(0)}) = \{0\to 1\to 2\to 0, 2\rightleftarrows 3\rightleftarrows 4,  4\rightleftarrows 5,
5\rightleftarrows 6, 6\to 7, 7\rightleftarrows 8, 8\to 6, 2\to 6\}$, then $E(G^{(1)}) =
 \{G_1\rightleftarrows 3\rightleftarrows 4,  4\rightleftarrows 5,
5\rightleftarrows 6, 6\to 7, 7\rightleftarrows 8, 8\to 6, G_1\to 6\}$,
$ E(G^{(2)}) = \{G_1 \rightleftarrows 3\rightleftarrows  G_2, G_2\rightleftarrows 6, 6\to 7, 7\rightleftarrows 8, 8\to 6, G_2\to 6\}$, $E(G^{(3)}) = \{G_3\rightleftarrows 6,  6\to 7, 7\rightleftarrows 8, 8\to 6\}$, $E(G^{(4)}) = \{G_3\rightleftarrows 6, 6\rightleftarrows G_4\}$,
$E(G^{(5)}) = \{G_3\rightleftarrows G_5\}$,
$E(G^{(6)}) = \emptyset$; $G^{(0)}(0)=0, G^{(1)}(0)=G^{(2)}(0)=G_1$, $G^{(3)}(0) = G^{(4)}(0) = 
G^{(5)}(0)= G_3$, $G^{(6)}(0)= G_6$.  Encapsulated graphs are the internal (self-) edges of 
compound vertices, here $E(G_{1}) = \{0\to 1\to 2\to 0\}, E(G_2) = \{4\rightleftarrows 5\}, 
E(G_3) = \{G_1\rightleftarrows 3 \rightleftarrows G_2\}$, $E(G_4) = \{7\rightleftarrows 8\}$, 
$E(G_5) = \{6\rightleftarrows G_4\}, E(G_6) = \{G_3 \rightleftarrows G_5\}$. 

\Medskip  The  encapsulated graph  $\{\sigma_1,\ldots,\sigma_n\}\in \Sigma^{(j)}$, $n\ge 1$ is formed only if  $\sigma_1,\ldots,\sigma_n\in \Sigma^{(j-1)}$ form a  {\bf strongly connected} 
 structure, i.e.  there exists a path
connecting any pair $(\sigma_i,\sigma_{i'})$. Strongly connected {\em components} (SCCs for short) are {\em maximal} path-connected vertex 
subsets (they are equivalent to communication classes in the Markov chain terminology); here, 
encapsulated graphs $G_i, i=1,2,\ldots$ are not SCCs in the usual sense, however, they are  maximal
 dominant SCCs of subgraphs of $G$ defined inductively by lowering the cut-off scale (see below).  Sequences $\sigma \to G^{(1)}(\sigma)\to G^{(2)}(\sigma)\to\cdots \to G^{(q)}(\sigma) = G^{(q+1)}(\sigma) = \cdots \equiv G^{(\infty)}(\sigma)$, $\sigma\in \Sigma^{(0)}$ become stationary after a finite number of steps. Compound vertices
 $G^{(\infty)}(\sigma)$ are called {\em terminal vertices}. Eliminating repeated steps yields
 the {\bf descending sequence} $\{\sigma\}=G^{(i_0)}(\sigma)\to G^{(i_1)}(\sigma)\to G^{(i_2)}(\sigma) \to \cdots\to G^{(i_p)}(\sigma) =  G^{(\infty)}(\sigma)$, $0=i_0<i_1<i_2<\ldots<i_p$, for instance,
 $(0 = G^{(0)}(0))\to (G_1 = G^{(1)}(0)=G^{(2)}(0))\to (G_3 = G^{(3)}(0)=G^{(4)}(0)=G^{(5)}(0)) \to (G_6 = G^{(6)}(0))$, with $i_1=1,i_2=3,i_3=6$. To this
 sequence is associated  an ensemble of paths, here 
\BEQ \{(0\to 1\to 2)\to G_5, ((0\to 1\to 2)\to 3\to G_2)\to G_5\}
\EEQ   
 called {\bf descending compound vertex  path}, by concatenating successively $\sigma$, an excursion  from $\sigma$  (i.e. simple path started from $\sigma$)
 in $G^{i_1}(\sigma)$, an excursion from $G^{(i_1}(\sigma)$ in $G^{(i_2)}(\sigma)$, and so on. 
A descending compound vertex path may be further detailed as a number of associated {\bf descending
paths} by replacing each compound vertex $G_i$ by a simple path connecting vertices in the set 
$(G^{(i)})^{-1}(G_i) \subset \Sigma$ of $G$-vertices contained in $G_i$.      For instance, 
on Fig. \ref{fig:merging}.2, the   sequence $0\to G_1 \to G_3\to G_6$ descending
from $0$ may be seen as the ensemble of descending paths $\{\gamma_i\cdot(\gamma_i^{end}\to 6), \gamma_i\cdot(\gamma_i^{end}\to 6) \cdot (6\to 7), 
\gamma_i\cdot(\gamma_i^{end}\to 6)\cdot (6\to 7\to 8)\}, i=1,13\}$, where $\gamma_1=(0\to 1\to 2)$, resp.  $\gamma_{13} =  (0\to 1\to 2)\cdot(2\to 3)\cdot(3\to 4) \cdot (4\to 5)$, and $\gamma_i^{end} = \begin{cases}
2 \\ 5 \end{cases}$ is the end-point of $\gamma_i$. Note how $\gamma_{13}$ has been
 decomposed: $(0\to 1\to 2)$ is a path {\em circulating within} $G_1$; $(2\to 3)$ is an edge {\em leaving} (outgoing from) $G_1$;
 $(3\to 4)$ is an edge {\em reentering} $G_2$; and $(4\to 5)$ is a path {\em circulating within} $G_2$. The edge
 $(\gamma_1^{end}\to 6) = (2\to 6)$ {\em leaves} $G_1$, then $G_3$, then reenters $G_5$ at $6$; the
 edge $(\gamma_{13}^{end}\to 6) = (5\to 6)$ {\em leaves} $G_3$, then reenters $G_5$ at $6$. Thus
 we have {\bf circulating paths}, and  {\bf leaving}, resp. {\bf reentering} edges. 
 {\em Leaving edges} cross a number of {\em barriers} before possibly reentering some compound
 vertex, for instance, $(\gamma_1^{end}\to 6)$, resp. $\gamma_{13}^{end}\to 6$, crosses 
 2, resp. 1 barrier. 

\Medskip  By construction, if $\Sigma^{(j_{max})} = \{G_{j_{max}}\}$ is made up of a single
vertex, and $\sigma\in \Sigma$, paths descending from $\sigma$ connect $\sigma$ to all
other vertices $\sigma'$ in $\Sigma$. {\em For fixed $\sigma'\in \Sigma$, it is actually possible
to connect $\sigma$ to $\sigma'$ within $G^{(i_{\sigma,\sigma'})}(\sigma)$}, where 
$i_{\sigma,\sigma'} = \min\{i \ |\ G^{(i)}(\sigma)=G^{(i)}(\sigma')\}$, for instance (see 
Fig. \ref{fig:merging}.2), $0$ and $5$ may be connected within $G_3$ since 
$G^{(2)}(0)=G_1\not=G_2=G^{(2)}(5)$ and $G^{(3)}(0) = G^{(3)}(5) = G_3$.  
 The box around $G_i$ (see Fig. \ref{fig:merging}.2) in the hierarchical representation separates
compound vertices $\sigma_{int}$ inside $V(G_i)$ from compound vertices $\sigma_{ext}$ outside $V(G_i)$; it also separates  vertices in the set $(G^{(i)})^{-1}(G_i) \subset \Sigma$
from other vertices in $\Sigma$. The important point is, our multi-scale construction
ensures that $G_i$'s are {\em basins}, i.e.  the resolvent formula \S \ref{subsection:resolvent} give the most weight to sums over paths that explore systematically and a large 
number of times  each $G_i$ they come across.  For instance, the descending path 
$0\to 1\to \cdots \to 5$ within $G_3$, coming from the  descending compound vertex path $(0\to 1\to 2)\to 3
\to G_2$ associated to the stopped descending sequence $0\to G_1\to G_3$,   does not have a large weight (see below (\ref{eq:resolvent-formula})).
On the other hand, let ${\cal C}_1=(2\to 0\to 1\to 2)$, resp. ${\cal C}_2 = (4\to 5\to 4)$,
circle around $G_1$, resp. $G_2$. Then the sum of the weights of the paths $(0\to 1\to 2)\cdot  ({\cal C}_1^{n_1})\cdot(2\to 3)\cdot \  \prod_{i_3=1}^{n_3-1} \Big\{ (3\to \sigma_{i_3})\cdot 
\Big(({\cal C}_{\sigma_{i_3}/2}^{n_{\sigma_{i_3}/2,i_3}})\cdot(\sigma_{i_3}\to 3)\Big) \Big\} \cdot 
(3\to 4) \cdot \Big( {\cal C}_2^{n_{2,n_3}}\cdot(4\to 5)\Big) $ with $\sigma_{i_3}=2$ or $4$, and  $n_3, (n_{1,i_3})_{i_3=1,\ldots,n_3},
(n_{2,i_3})_{i_3=1,\ldots,n_3}$ 
ranging from 1 to large enough values,  which circle around $G_1,G_2$ and $G_3$ a 
large number of times, {\em may} have a large weight. Such paths are called {\em leading paths}
in \S \ref{subsubsection:full-non-auto}. Conversely, {\em weights of
leaving edges}, such as e.g.  
$(2\to 3)$ {\em leaving} $G_1$, or $(2\to 6)$ {\em leaving}  $G_1$  (and then, $G_3$), 
later on (see \S \ref{subsubsection:exit-proba}) leveraged to produce exit probabilities,
{\em are  small}. 

\Bigskip Compared to a hierarchical graph, a hierarchical structure has two main additional
structures:
\begin{itemize}
\item[(i)] all encapsulated graphs are derived from a graph $G^{(0)}$ on $\Sigma$;
\item[(ii)] merging patterns are ordered (from left to right on Fig. \ref{fig:merging}) 
by the upper index $j=0,\ldots,j_{max}$. To the $j$-th {\bf merging step} leading to $\Sigma^{(j)}$ from
$\Sigma^{(j-1)}$ is associated a  kinetic scale $n(j)$, called {\bf step $j$ cut-off scale} in the description of the renormalization procedure (see \S \ref{section:description}), with $(n(j))_{j\ge 1}$ strictly
increasing. 
\end{itemize}


\subsection{Multi-scale method} \label{subsection:multi-scale-method}


The main line of thought in this work, coming from a long experience in the subject \cite{Unt-ren1,MagUnt1,MagUnt2},  is directly borrowed from
 rigorous  renormalization group methods in 
quantum field theory, see e.g. \cite{Salm}, or \cite{FVT}, chap. 1 for a quick, 
self-contained introduction. The idea of the multi-scale method is to deal with transitions not all at once, but
starting from the highest scales (fastest transitions), i.e. from rates $k$ such that $
\lfloor \log_b (k) \rfloor = n^{(0)}_1 := \max_{\rho}  \lfloor \log_b(k^{\rho})\rfloor$, and incorporating successively rates of scales $n^{(0)}_2>n^{(0)}_3>\cdots$
equal to the transition scales of the reaction network in decreasing order. The simplest
non-trivial example is the so-called Michaelis-Menten model of catalysis\footnote{The reactions involve an enzyme $E$ and are usually written $ES\rightleftarrows S+E, ES\rightleftarrows P+E$. The
enzyme merely facilitates the transition from the substrate $S$ to the product $P$. In the regime
(here considered) where $E$ is abundant, $X_E$ may be considered as constant, and enters the kinetic rate $k_+$ as a mere multiplicative 
constant.} 
$(\Sigma,E) = (\{S,ES,P\}, \{ES\overset{k_-}{\to} S, S\overset{k_+}{\to} ES,  
ES\overset{k_2}{\to} P, P\overset{\bar{k}_2}{\to} ES\})$, see e.g. \cite{Chor}, in the regime where e.g.  $n^{(0)}_1= 
\lfloor \log_b (k_-) \rfloor, n^{(0)}_2 = \lfloor \log_b (k_+)\rfloor, n^{(0)}_3 = \lfloor \log_b (k_2)\rfloor,  n^{(0)}_4 = \lfloor \log_b (\bar{k}_2)\rfloor $, with $n^{(0)}_1>n^{(0)}_2>n^{(0)}_3>n^{(0)}_4$. Since $(e^{tA})_{\sigma,\sigma'}=
tk_{\sigma'\to\sigma} + o_{t\to 0}(t)$ $(\sigma'\not=\sigma)$, it is natural to think
of {\em time-scales} as the opposite of the kinetic scales, with transition times $t_{\sigma\to\sigma'} \sim 
k^{-1}_{\sigma\to\sigma'} \sim b^{-n_{\sigma\to\sigma'}}$, where $n_{\sigma\to\sigma'} = 
\lfloor \log_b (k_{\sigma\to\sigma'})\rfloor$. Thus the multi-scale method considers short-time evolution first, gradually shifting to long-time evolution. Considering only reactions
with scale $n^{(0)}_1$, we  get the reducible network $ES\overset{k_-}{\to} S$. 
Incorporating  scale $n^{(0)}_2$, on the other hand, yields a cycle $ES \underset{k_+}{\overset{k_-}{\rightleftarrows}} S$ at scale $n(1)=n_2^{(0)}$. Dismissing the in- and out-going edges $ES\rightleftarrows P$ 
which operate on a longer time scale, we may conjecture that the subnetwork $(ES,S)$ is at near equilibrium with relative concentration  $\frac{X_{ES}}{X_S} = \frac{k_+}{k_-}$ equal
to the ratio of components of the Markov chain stationary measure. Indeed (as follows from elementary computation), the presence of 
the outgoing edge multiplies the above ratio at stationarity by $1+O(\frac{k_2}{k_-}) = 1+ O(b^{-(n_1^{(0)}-n_3^{(0)})})$.  The correcting term $\frac{k_2}{k_-}$, an instance
of {\bf renormalization factor}, is equal to
the ratio $k_{ES\to P}/k_{ES\to S}$ of the outgoing rate $ES\to P$ and the internal rate
$ES\to S$.
  
\Medskip The above fact may be proved as a particular case of a perturbative fixed-point argument 
(see Appendix, Section \ref{section:app-fixed-point}) which is mathematically rigorous {\bf when $b$} (playing the role of the inverse of 
a coupling constant between scales) {\bf is large enough}. The argument implies that the Lyapunov
data of the network are a perturbation of those of an {\bf effective network}  $G^{(1)}= (\Sigma^{(1)},E^{(1)}) = (\{G_1=\{S,ES\},P\}, G_1\rightleftarrows P)$ with {\bf renormalized
kinetic rates} $k_{G_1\to P} \sim \frac{k_+}{k_-} k_2, \  k_{P\to G_1} \sim 
k_{P\to ES} = \bar{k}_2$. Thus $\frac{X_P}{X_{G_1}}\sim \frac{k_{G_1\to P}}{k_{P\to G_1}}$ 
at stationarity, with $X_{G_1}\sim X_S + X_{ES} \sim X_S$, completing the determination
of the Lyapunov eigenvector, here stationary measure.  Following general practice, we call $G^{(0)}$, resp. 
$n^{(0)}_i, i\ge 1$ the {\bf bare}
network, resp. kinetic rates. 

\Medskip In more complex networks, the above procedure must be repeated recursively, letting
the {\bf cut-off scale} (lowest scale under consideration) decrease. It is important to understand that scales are step-dependent; for instance, the two-scales of $G^{(1)}$ in the previous
example are $n_{G_1\to P} = \lfloor \log_b(k_{G_1\to P})\rfloor$ and $n^{(0)}_4$; thus 
$n^{(1)}_1 = n_{G_1\to P}, n^{(1)}_2= n^{(0)}_4$, resp. $n^{(1)}_1 = n^{(0)}_4, n^{(1)}_2= n_{G_1\to P}$ depending on the relative ordering of these scales.    Step $j$ consists
in {\em collapsing} {\bf maximal dominant SCCs} (also called: {\bf clusters}, by analogy
 with cluster decompositions, see \cite{Salm}) with lowest scale $n(j)$. These structures are extracted from the 
subgraph $G^{(j)}_{dominant}$ of dominant edges of $G^{(j)}$, defined as edges $(\sigma,\sigma')
\in E(G^{(j)})$ such that $n(k_{\sigma\to\sigma'})
= \max_{\sigma''} n(k_{\sigma\to\sigma''})$.  Maximal dominant SCCs are maximal classes (in the
terminology of Markov chains) of $G^{(j)}_{dominant}$.  This procedure selects in priority edges $(\sigma,\sigma')$ such that $w_{\sigma\to\sigma'}$ is of order 1, since (see   \S \ref{subsection:resolvent}, and Section \ref{section:notations}
for details) $w_{\sigma\to\sigma'} \approx \frac{k_{\sigma\to\sigma'}}{\sum_{\sigma''} k_{\sigma\to\sigma''}} \approx \frac{k_{\sigma\to\sigma'}}{\max_{\sigma''} (k_{\sigma\to\sigma''})}$ in general; thus, dominant edges give the maximum contributions to sums over paths.  

\Medskip A maximal dominant SCC which is autocatalytic changes the procedure as described above. 
Namely, add to the previous network the doubling reaction $S\overset{\nu_+}{\to} ES + ES$ with 
$\nu_+\ll k_{\pm}$. Dismissing as above the outgoing reaction $ES\to P$ so that the subgraph
$G_1$ with vertex set $\{S,ES\}$ is closed, it follows from  elementary
computations, or from our general procedure (see Example 1 in \S \ref{subsection:examples}) that the network is autocatalytic, with Lyapunov exponent $\lambda_G\sim \nu_+ $. However, this 
rate enters in competition with the rate $k_2$ of the outgoing reaction; taking this into account, one can prove that the SCC is autocatalytic, resp. is not, if  $\nu_+ \gg k_2$, resp. $\ll k_2$. This leaves undetermined a region (called: {\bf resonance region}) where 
$\nu_+  \sim k_2$.  It is convenient to symbolize the Lyapunov exponent $\lambda_G$ of $G$ by
a self-edge $G\to G$ with same scale, since it has the same net effect as a doubling reaction 
$G\overset{\lambda_G}{\to} G+G$. 

\Medskip  In the autocatalytic case, $w(\lambda_G)_{ES\to P}\ll 1$ is small, i.e. of order $O(b^{-1})$ at most. Generalizing, all 
outgoing weights $w(\lambda_G)_{G\to \sigma'}$ with $G$ autocatalytic are {\em small}, because
rates $G\to\sigma'$ are subdominant compared 
to the scale of the self-edge $G\to G$. This forbids the formation of further hierarchical structures
involving $G$. 

\Medskip On the contrary, it may happen that some degradation reaction (described as $\sigma\to\emptyset$, where $\emptyset$ is some formal vertex) becomes dominant. Then the trivial
maximal dominant SCC $\{\emptyset\}$
is recognized as an absorbing state, and the algorithm returns simply $v^*_{\sigma} \sim \del_{\sigma,\emptyset}$ up
to normalization, implying in particular $\lambda^*<0$ (but $\lambda^*$ is not estimated). The case when  
degradation rates remain non-dominant throughout, but no autocatalytic non-trivial maximal dominant SCC is formed,  defines the 'shadow zone' first mentioned
in \S \ref{subsection:aim}.  

\Medskip In the end, all dominant SCCs have been merged, and the behavior of the r.-h.s. in 
the resolvent formula (\ref{eq:resolvent-formula}) may be captured by summing over the a subset 
of the finite set of {\em simple} (i.e. non-self-intersecting) paths defining a directed acyclic graph (DAG). This reduces
the task of estimating the Lyapunov data to a finite combinatorial problem.

\Medskip The recursive nature of the algorithm also produces Lyapunov data estimates for
intermediate, cut-off graphs $G_{cut}(j)$ associated to each step $j$, with lowest scale $n(j)$. These provide an adequate summary of the network at time scale  $\sim b^{n(j)}$, and may be 
expected to play a role in transient dynamics. 
Actual computations (not provided here) imply that this is indeed the case, though the different time-scale descriptions get mixed up.

\Medskip The whole procedure may be understood as a $b\to \infty$ expansion. However, numerical
 simulations (not included here) for networks which are low-dimensional, but already have
 an involved, multi-step hierarchical structure, suggest that our estimates are good enough
 even for low values of $b$ (including $b\simeq 10$, or even $b\simeq 2$).


\subsection{Overview of the article}


Fundamental graph-theoretic concepts and notations are introduced in Section 2, including 
{\em  cut-off graphs}. Section 3  shows how renormalization works on a simple one-step example, and
gives heuristics for the renormalization algorithm.  The general, multi-step renormalization procedure, leading to the hierarchical formulas (\ref{eq:pisigma1})--(\ref{eq:vdaggersigma}),  is described precisely in Section 4; two detailed examples are provided at the end of the section, supported by numerical simulations. The general skeleton of the proof is found in Section 5; some technical lemmas which may be skipped in a first reading, including the fixed-point arguments underlying the proof, have been postponed to the appendices. An index
of main definitions and notations is provided at the end of the article.


\section{Main notations}  \label{section:notations}


\noindent {\bf Graph of split reactions.}  The starting point
is the graph of split reactions $G=(\Sigma,E)$ (later on denoted: $G^{(0)}$, and called: {\bf bare graph}).   By construction, 
$\sigma\to \sigma'$  $(\sigma\not=\sigma')$ is an edge in $E$ if (i) either there exists
a $1 \to 1$ reaction $\sigma\to \sigma'$;  or (ii) $\sigma
\to \sigma'$ is a {\em 'split reaction'} coming from a $1\to 2$
reaction $\sigma\to \sigma' + \sigma''$.  There are no
multiple edges. The rate  $k_{\sigma\to\sigma'}$ attached to $(\sigma,\sigma')\in E$ is equal to the sum of the rates $k^{\rho}$ of reactions  (i) or
split reactions (ii) from $\sigma$ to $\sigma'$, counting
multiplicities, and gives the off-diagonal coefficients of the adjoint generator $A$, 
$A_{\sigma',\sigma} = k_{\sigma\to\sigma'}$ $(\sigma\not=\sigma')$.  Degradation reactions $\sigma\overset{\beta_{\sigma}}{\to} \emptyset$ are
considered as $1\to 1$ reactions with $\sigma'=\emptyset$; only the kinetic rates receive a 
specific notation.  The degree $p(\rho)$ of a reaction $\rho$ is either (i) 1 or (ii) 2.
Diagonal coefficients of $A$ are $<0$. By construction, (i)  a reaction $\rho : \sigma 
\overset{k^{\rho}}{\to} \sigma'$,  resp. (ii) $\rho:\sigma \overset{k^{\rho}}{\to} 
\sigma'+\sigma''$,  
 contributes $-k^{\rho}$ to $A_{\sigma,\sigma}$ and  $0$, resp. $+k^{\rho}$ to
$\sum_{\sigma^*\in\Sigma} A_{\sigma^*,\sigma}$. 

\Medskip {\bf Autocatalysis.} We say that $G$ is autocatalytic if $\lambda^*_G:=\lambda^*(A)>0$. 

\Medskip {\bf Total outgoing rate.} Let 
\begin{equation} k_{\sigma}:= \sum_{\sigma'\not=\sigma} k_{\sigma\to\sigma'}  \label{eq:k} \end{equation}
be the sum of the rates of reactions outgoing from $\sigma$. 

\Medskip {\bf Deficiency rate.} Let
\begin{equation} \kappa_{\sigma} := \sum_{\rho\ |\ p(\rho)=2,\ \rho: X_{\sigma}\to \cdots}
k^{\rho}  \label{eq:deficiency-rate}
\end{equation}
be the sum of the kinetic rates of all  $1\to 2$ reactions $\rho$ with reactant $\sigma$. 
If $\kappa_{\sigma}>0$, we add to $G$ the {\em self-edge} $\sigma\overset{\kappa_{\sigma}}{\to} \sigma$ (then edges $\sigma\to\sigma'\not=\sigma$ are called {\em non-trivial}).  By
construction,
\BEQ \sum_{\sigma'} A_{\sigma',\sigma} = \kappa_{\sigma}.
\label{eq:deficiency-rate2}
\EEQ
 If there are no $1\to 2$ reactions,
 then $\kappa_{\sigma}=0$, and $\sum_{\sigma'} A_{\sigma',\sigma}=0$ expresses Markov chain probability
 conservation. In general, the generator $A$ is {\bf defective}, meaning that probability is not conserved; we call $A$  a {\em defective adjoint Markov generator}.  This brings us to the next point.

\Medskip {\bf Associated Markov chain.} We let $\tilde{A}=(\tilde{A}_{\sigma',\sigma})_{\sigma',\sigma\in \Sigma}$
be the matrix with  off-diagonal coefficients $\tilde{A}_{\sigma',\sigma} =A_{\sigma',\sigma}$ $(\sigma\not=\sigma')$ and negative diagonal coefficients
 
\begin{equation} \tilde{A}_{\sigma,\sigma}=-\sum_{\sigma'\not=\sigma} A_{\sigma',\sigma} = A_{\sigma,\sigma}-\kappa_{\sigma} = -k_{\sigma}.
\end{equation}
 By construction \cite{Anderson}, $\tilde{A}$ is a "true" (i.e. non-defective) adjoint Markov generator, so
that its Lyapunov exponent (i.e. eigenvalue with largest real part)
is $0$. 

\Medskip {\bf Overall degradation rate.} Let $\alpha\ge 0$ be a constant. Then 
$A(\alpha):=  A - \alpha \Id$ has same off-diagonal coefficients as $A$, and diagonal coefficients
$|A(\alpha)_{\sigma,\sigma}|= |A_{\sigma,\sigma}| +\alpha$.  

\Medskip {\bf Edge scales.} {\em Fix a scale parameter $b>1$ once and for all.}  If $(\sigma,\sigma')\in E$,
we let 
\begin{equation} n_{\sigma\to\sigma'} := \lfloor \log_{b} (k_{\sigma\to
\sigma'}) \rfloor \in \Z  \label{eq:edge-scales}
\end{equation}
Because of the truncation error, all our results are valid up to one or a few scale units. We
write $k\prec k'$, resp. $\preceq k'$ when $\log_b (k/k')<-1$ (resp. $<0$) up to a few units, and use similarly curved
symbols $\succ,  \succeq$; symmetrizing, $k,k'$ have same order if
\begin{equation} (k\sim k') \Leftrightarrow  (k\preceq k'\preceq k)
\end{equation}
In particular, we often use the 'sum-max' substitution rule for positive quantities (kinetic rates or kinetic rate ratios)
\begin{equation}  (k_1\succeq k_2\succeq\cdots\succeq k_m) \Rightarrow 
(k_1+\cdots+k_m \sim k_1) 
\end{equation}
 Errors do not accumulate as a result of successive approximations provided $b$ is large enough. How large depends intricately on the size (number of species) of the 
maximal dominant SCCs, but simulations on small- to average-sized networks tend to suggest
that truncation errors do not invalidate our estimates as soon as $b> 3$ or $4$. Thus, 
$k\prec 1$ means in principle: $k$ small.

\Medskip {\bf Dominant edges.} An edge
$\sigma\to\sigma'$ is dominant  if 
\begin{equation} n_{\sigma\to\sigma'}  = \max_{\sigma''\in \Sigma}
n_{\sigma\to\sigma''}  \label{eq:dominant edges}
\end{equation}
Generalizing, if $\alpha>0$, $n_{\alpha} = \lfloor \log_b \alpha \rfloor$, then 
$\sigma\to\sigma'$ is $\alpha$-dominant if $n_{\sigma\to\sigma'}\succeq n_{\alpha}$ and 
$\sigma\to\sigma'$ is dominant (letting $\alpha\to 0$, the two notions coincide). 
 As a general rule, $\alpha$-dominant edges appear as {\bf bold} lines on all our graphs.
An important particular case is when  $\sigma$ is {\bf autocatalytic}, i.e.  $\kappa_{\sigma}\succ k_{\sigma}$, for then 
$n_{\sigma\to\sigma}$ is dominant and the equivalent doubling reaction $\sigma\overset{\kappa_{\sigma}}{\to} \sigma +\sigma$ ensures that $X_{\sigma}(t)$ increases
at least as fast as $\sim e^{\kappa_{\sigma}t}$ (possibly faster if there exists a path
from $\sigma'$ to $\sigma$, with $\kappa_{\sigma'}\succ \kappa_{\sigma}$).

\Medskip {\bf Vertex scales.} A {\em vertex} $\sigma$ is an element
of $\Sigma$. The scale of $\sigma$ is 
\begin{equation} n_{\sigma}:= \max_{\sigma'\in \Sigma}\,  n_{\sigma\to\sigma'} 
\label{eq: vertex scales}
\end{equation}
i.e. the maximum scale of all edges $\sigma\to \sigma'$
 of $G$ outgoing from $\sigma$. 
 
\Medskip Note that an edge $\sigma\to\sigma'$ is dominant if and only if 
$n_{\sigma\to\sigma'} = n_{\sigma}$.  
 By construction,  $\log_{b}  (k_{\sigma}) = \log_b (|\tilde{A}_{\sigma,\sigma}|) 
\sim \log_b (|A_{\sigma,\sigma}|) 
 \sim n_{\sigma}$ if $\sigma$ is {\em not} autocatalytic.

\Medskip {\bf Deficiency weight.}  Let $\alpha\ge 0$ be an overall degradation rate. The $\alpha$-deficiency
weight is the ratio
\begin{equation} \eps_{\sigma} : = \frac{\kappa_{\sigma}}{|A(\alpha)_{\sigma,\sigma}|} = 
\frac{\kappa_{\sigma}}{|A_{\sigma,\sigma}| + \alpha}.  \label{eq:deficiency-weight}
\end{equation}
 For example, a couple consisting of a  $1\to 1$ reaction $\sigma\overset{k_{off}}{\to} \sigma'$ and a $1\to 2$ reaction
$\sigma\overset{\nu_+}{\to} \sigma' + \sigma'$ produces a deficiency weight $\eps_{\sigma} = \frac{\nu_+}{k_{off}+\nu_+ + \alpha}$. If $\alpha\succ k_{off},\nu_+$, then $\eps_{\sigma}\sim
\frac{\nu_+}{\alpha} \prec 1$; this  regime is called {\em degraded}. In the contrary case, 
  $\eps_{\sigma}\sim \min(\frac{\nu_+}{k_{off}},1)$. 
  
  \Medskip Note that, by construction (see first paragraph of the present section), it always holds: $0\le \eps_{\sigma}\le 1$.

\Medskip {\bf Scale ordering.} We let $(n_i)_{i=1,2,\ldots,n_{max}}$ be the
set of edge and vertex scales, and  order them by decreasing order,
$n_1> n_2>\cdots$.  Representing them {\bf from top to bottom} yields a {\bf multi-scale graph} (see Supp. Info. for examples). By construction, there is 
a reaction $\sigma\to\cdots$ of scale $n_{\sigma}$ {\em above} any
reaction with reactant $\sigma$. In particular, there is at least 
one vertex of scale $n_1$.  

\Medskip {\bf Transition weights.}  Let, for $\alpha\ge 0$, 
\begin{equation} w(\alpha)_{\sigma\to\sigma'}= \frac{k_{\sigma\to\sigma'}}{|A_{\sigma,\sigma}| + \alpha} \ \ (\sigma\not=\sigma')
\end{equation}
In particular, $w_{\sigma\to\sigma'} := w(0)_{\sigma\to\sigma'}$. Similarly, we let 
$\tilde{w}_{\sigma\to\sigma'} := \frac{k_{\sigma\to\sigma'}}{|\tilde{A}_{\sigma,\sigma}|} = 
\frac{k_{\sigma\to\sigma'}}{k_{\sigma}}$.

\Medskip {\bf Path weights. }  The {\em weight} $w(\alpha)_{\gamma}$ of a path $\gamma: \sigma_1\to\cdots\to \sigma_{\ell}$  is the
product of the transition weights $w$ along its edges,
\begin{equation} w(\alpha)_{\gamma}= \prod_{i=1}^{\ell-1} w(\alpha)_{\sigma_i\to\sigma_{i+1}}.
\end{equation}

\Medskip {\bf Stationary measure, stationary weights.}  Let ${\bf 1}$ be the constant vector 
with components ${\bf 1}_{\sigma}=1$. The zero-diagonal matrix $\tilde{\cal W}= (\tilde{w}_{\sigma\to 
\sigma'})_{\sigma,\sigma'}$ is a Markov transition matrix since $\tilde{\cal W}{\bf 1} = {\bf 1}$ by construction.  A stationary measure
of the discrete-time Markov chain with transition weights
$\tilde{w}$ is a positive distribution  $\tilde{\pi}$ such that  $\tilde{\pi} \tilde{\cal W} = 
\tilde{\pi}$.  
 Let  $\tilde{\mu}_{\sigma} := k_{\sigma}^{-1} \tilde{\pi}_{\sigma}$; 
plugging the definition of  $\tilde{w}_{\sigma\to\sigma'} $  into the above eigenvector identity 
yields $\sum_{\sigma\not=\sigma'} \tilde{\mu}_{\sigma} k_{\sigma\to\sigma'} = |\tilde{A}_{\sigma',\sigma'}|\, \tilde{\mu}_{\sigma'}$, from 
which $\tilde{A}\tilde{\mu}=0$. Thus $\tilde{\pi}$ is a stationary
measure of the discrete-time Markov chain if and only if $\tilde{\mu}$   is a stationary 
composition for $\tilde{A}$. 
We work most of the time with discrete-time transition weights, and
our results will provide estimates for the associated discrete-time
weights, called {\bf Lyapunov weights}, 
\begin{equation} \pi^*_{\sigma} := (|A_{\sigma,\sigma}|+\lambda^*) v^*_{\sigma}. \label{eq:pi*}
\end{equation}   
Then $Av^* = \lambda^* v^*$ if and only if $\sum_{\sigma\not=\sigma'}
\pi^*_{\sigma} w^*_{\sigma\to\sigma'} = \pi^*_{\sigma'}$, where 
$w^*_{\sigma\to\sigma'} \equiv w(\lambda^*)_{\sigma\to\sigma'}$.

\Medskip 
{\bf Dominant paths.} From (\ref{eq:dominant edges}), (\ref{eq: vertex scales}), the edge $\sigma\to\sigma'$  is 
dominant if  $n_{\sigma\to\sigma'}=n_{\sigma}$. 
If $\sigma\not=\sigma'$, the condition is  equivalent to 
 $w_{\sigma\to\sigma'}\sim 1$.  Iterating, we say 
that $\sigma$ is connected to $\sigma'$ by a dominant path
if there exists a path $\gamma:\sigma=\sigma_1\to \cdots \to 
\sigma_{\ell}=\sigma'$ whose edges $\sigma_i\to\sigma_{i+1}, \ 
i=1,\ldots,\ell-1$ are all dominant. This defines a notion of 
{\em dominant path} (weak or strong) {\em connectivity}.

\Medskip {\bf Dominant graph.} A   dominant graph  is a  graph 
connected by dominant paths, i.e. a graph $G=(V(G),E(G))$ such that all oriented
pairs of vertices $(\sigma,\sigma')$
are connected by dominant paths. The {\bf dominant subgraph} of a graph $G$ is obtained 
by selecting only dominant edges of $G$; even if $G$ is connected, it may be disconnected.
 
 \Medskip{\bf Dominant SCCs.} If $\sigma,\sigma' \in V(G)$, 
we let $\sigma\sim\sigma'$ (including the case $\sigma=\sigma'$) if there exists a path 
$\sigma=\sigma_1\to \cdots\to \sigma_{\ell}=\sigma'$ of dominant 
edges of $G$ connecting $\sigma$ to $\sigma'$, and
similarly, a path of  dominant 
edges of $G$ connecting $\sigma'$ to $\sigma$. The
relation $\sim$ is an equivalence relation on vertices of 
$G$. We call {\bf dominant SCCs} (dominant strongly
connected components) its equivalence classes. {\bf Maximal dominant SCCs} of $G$ are 
dominant SCCs which are 'downstream', i.e. which are not the source of any dominant edge (in the terminology of 
Markov chains \cite{Norris}, they are maximal classes of the dominant subgraph of $G$).  {\bf Non-trivial dominant SCCs} are dominant SCCs
 containing at least two vertices; equivalently, containing a dominant cycle.

\Medskip The main objects of consideration for renormalization are {\bf non-trivial maximal
dominant SCCs}, that is, maximal dominant SCCs which are non-trivial. Note that the presence
of a dominant cycle  implies the existence of a  non-trivial dominant SCC, but not necessarily
maximal.

\Medskip {\bf Infra-red cut-off.} Fix a scale $n$ (called: {\em infra-red cut-off scale}, or 
simply, {\em cut-off scale}). Let $\Sigma^{\ge n}:=\{\sigma\in \Sigma \ |\ n_{\sigma}\ge n\}$ and 
$E^{\ge n}:=\{(\sigma,\sigma')\in E\ |\ n_{\sigma\to\sigma'}\ge n\}$. Then $G_{\searrow n}:=(\Sigma^{\ge n}, E^{\ge n})$ is called {\em scale $n$ (infra-red) cut-off graph}.

\Medskip {\bf Cut-off graphs.} Split $\Sigma$ into a disjoint union
$\Sigma=\Sigma^{int}\uplus \Sigma^{ext}$ (internal species, versus external species).
The cut-off graph ${\cal G}^{int} = (\Sigma^{int}, {\cal E}^{int})$ has internal edge set
${\cal E}^{int} := \{(\sigma,\sigma')\in E\ |\ \sigma,\sigma'\in \Sigma^{int}\}
\subset E$, and degradation rates 
\BEQ \beta^{int}_{\sigma} := \beta_{\sigma} + \sum_{\sigma'\in \Sigma^{ext}} k_{\sigma\to\sigma'}.
\label{eq:beta_int}
\EEQ
  Roughly speaking, external
 species have been 'frozen': influxes have been discarded, and  outfluxes from $\Sigma^{int}$ to $\Sigma^{ext}$ 
 have been added to degradation rates.
 
 \Medskip  Instead of increasing degradation rates, one may add to the internal edges in $G^{int}$ the
 set ${\cal E}^{out}$ of outgoing edges $k_{\sigma\to\sigma'},\sigma\in \Sigma^{int},\sigma'\in 
 \Sigma^{ext}$. This yields an {\em extended cut-off  graph} ${\cal G}$, which is the same as the cut-off graph ${\cal G}^{int}$, except for the fact that the target of outgoing edges is specified. The distinction
 between $\cal G$ and ${\cal G}^{int}$ is sometimes immaterial. Both graphs are associated
 to the Markov chain on $G$ started inside $\Sigma^{int}$ and  stopped upon hitting the 
 "boundary" $\Sigma^{ext}$. 
Working on ${\cal G}$ instead of ${\cal G}^{int}$ makes it possible to know at which point of
$\Sigma^{ext}$ the chain stops.  
 
\Medskip The cut-off graphs produced by our algorithm have by construction the following essential
properties:  (i) every pair $\sigma,\sigma'$ of vertices in $\Sigma^{int}$ is connected by a 
path of dominant edges of ${\cal G}^{int}$; (ii) no outgoing edge  is dominant,
which means that every outgoing edge $(\sigma,\sigma'_{ext}), \sigma\in \Sigma^{int}, \sigma'_{ext}\in 
 \Sigma^{ext}$, is dominated by some internal edge $(\sigma,\sigma'), \sigma,\sigma'\in \Sigma^{int}$, 
 i.e. $k_{\sigma\to\sigma'_{ext}}\prec k_{\sigma\to\sigma'}$.    Such structures will be uncovered inductively by lowering a cut-off scale.


\section{Renormalization: a heuristic introduction}  \label{section:resummation}


This short section contains no element of proof, but motivates the renormalization algorithm detailed 
in Section \ref{section:description} on the ground of an elementary  example
(\S \ref{subsection:cycle-elementary}) and heuristic computations based on the resolvent formula
(\S \ref{subsection:heuristics}).


\subsection{Cycles: an elementary example}  \label{subsection:cycle-elementary}


The simplest possible example of dominant SCC is that of a cycle ${\cal G}^{int}=\{1\to 2\to \cdots \to p\to 1\}$
with dominant edges. We restrict to the case $p=2$ for the discussion (generalization to 
$p>2$ is straightforward). Write $k_1 = k_{1\to 2}, k_2 = k_{2\to 1}$. To make the cycle 
possibly autocatalytic, we include a doubling reaction $1\overset{\nu_+}{\to} 2+2$, with 
$\nu_+ \prec k_1$ equal to the deficiency rate $\kappa_1$. We assume that the cycle
is connected to the rest of the network  through outgoing edges $1\to \cdots$ and $2\to \cdots$, 
with total rates $\xi_1$, $\xi_2$, yielding a graph $\cal G$. Since we want ${\cal G}^{int}$ to be a dominant SCC, we assume that 
$\xi_1\prec k_1, \xi_2\prec k_2$; then ${\cal G}^{int},{\cal G}$ are cut-off graphs (see 
Section \ref{section:notations}, just above) satisfying properties (i) and (ii). By definition, $-A_{1,1}= k_1+\nu_++\xi_1\sim k_1$, 
$-A_{2,2}= k_2+\xi_2\sim k_2$; thus, the matrix $A$ restricted to the cycle, 
$A^{int} = \left(\begin{array}{cc} -k_1-\nu_+-\xi_1 & k_2 \\  k_1+2\nu_+ & -k_2-\xi_2 
\end{array}\right)$ is almost probability-preserving, in the sense that 'extended' deficiency
weights  $\frac{A_{1,1}+A_{2,1}}{|A_{1,1}|} = \frac{\nu_+-\xi_1}{|A_{1,1}|} \sim \frac{\nu_+-\xi_1}{k_1}, \frac{A_{2,2}+A_{1,2}}{|A_{2,2}|} = \frac{-\xi_2}{|A_{2,2}|}\sim 
\frac{-\xi_2}{k_2}$ (including the total outgoing rates $\xi_1,\xi_2$, compare with (\ref{eq:deficiency-rate2}) and (\ref{eq:beta_int}))   are $\prec 1$.   

\Medskip Renormalization requires the additional knowledge of the external structure of 
${\cal G}^{int}$; see discussion on cut-off graphs in Section \ref{section:notations}. We choose for
specificity
$\Sigma^{ext} = \{1_{ext},2_{ext}\}$, and assume that outgoing edges of ${\cal G}$ are 
$E^{out} = \{1\to 1_{ext},2\to 2_{ext}\}$, with rates $\xi_1,\xi_2$. Let $k_{min}=\min(k_1,k_2)$ 
be the lowest rate of ${\cal G}^{int}$, and $n_{min}$ its scale, called {\em lowest scale}.

\Medskip Visualizing sometimes helps understanding.   If (for instance) $\xi_1,\xi_2\prec k_{min}$, then ${\cal G}^{int} = {\cal G}_{\searrow n_{min}}$ is an infra-red cut-off graph; see multi-scale representation below, 
where scales decrease from top to bottom, and vertices are 'split' over several scales. The
edge-splitting convention is as in (\ref{eq:1.3}); the negative killing rate $1\overset{-\nu_+}{\to}
\emptyset$, equivalent to a doubling reaction $1\overset{\nu_+}{\to} 1+1$,  is  represented as a self-edge $1\to 1$ with rate $\nu_+$; we assumed here that $k_1\succ k_2\succ \xi_2\succ \xi_1$
have scales $n_1>n_2>n_3>n_4$; the zigzag line separates 'upper' scales $\geq n_{min} = n_2$ from the 'lower' scales $<n_{min}$.

\begin{center}
\begin{tikzpicture}

\draw(0,2) node {$2$};
\draw(3,2) node {$1$};

\draw[dotted](0,1)--(0,-0.65); \draw[dotted](0,-1.35)--(0,-3);
\draw[dotted](3,0.65)--(3,-3);

\draw[<-, ultra thick](0,1)--(3,1);  \draw(1.5,1.4) node {$k_1$};

\draw(4,1) node {$n_1$};

\draw[-<](3,0) arc(-160:160:0.4 and 0.25) ; \draw(4.25,0) node {$\nu_+$};
 
\draw[<-, ultra thick](3,-1)--(0,-1); \draw(1.5,0.4-1) node {$k_2$};
\draw(4,-1) node {$n_2$};

\draw[decorate, decoration = zigzag](-1,-1.5)--(4,-1.5);

\draw[->](0,-2.25)--(0.8,-2.25); \draw(0.4,-2.25+0.4)
node {$\xi_2$};  \draw(1,-2.25) node {\textbullet};  \draw(1,-1.85) node {$2_{ext}$};

\draw(4,-2) node {$n_3$};

\draw[->](3,-3.25)--(2.2,-3.25);
\draw(2.6,-3.25+0.4) node {$\xi_1$}; \draw(2,-3.25) node {\textbullet}; \draw(2,-2.85) node {$1_{ext}$};
\draw(4,-3) node {$n_4$};

\end{tikzpicture}
\end{center}

We emphasize that the computations below only depend on the assumptions made in the first and 
second paragraph, and not on the particular scale ordering chosen above; in particular, they
also hold true when ${\cal G}^{int}$ is not an infra-red cut-off graph. 

\Medskip The general purpose is to estimate the Lyapunov eigenvector $\lambda_{\cal G}$ of ${\cal G}^{int}$, and derive heuristically a renormalized graph equivalent to $\cal G$ to leading order, in which all vertices in $V({\cal G}^{int})=\{1,2\}$ have been merged to a compound vertex denoted $G$.

\Medskip {\bf Excursion weight.} In this setting, the total weight of  ${\cal G}^{int}$-excursions 
$\Phi(\alpha)_{\cal G} := \Phi(\alpha)_1 = \Phi(\alpha)_2$,  see (\ref{eq:excursion}), is close to 1 when $\alpha=0$,
so that the associated geometric series  either diverges or has a large sum.  When $\alpha=0$ 
and $\nu_+=0$,  
\BEQ \sum_{n=0}^{+\infty} (\Phi(0)_{\cal G})^n =  \frac{1}{1-w(0)_{1\to 2} w(0)_{2\to 1}}
= \frac{1}{1- \frac{k_1}{k_1+\xi_1} 
\frac{k_2}{k_2+\xi_2}}
\sim \frac{1}{1- (1-\frac{\xi_1}{k_1}) 
(1-\frac{\xi_2}{k_2})}  \sim   \frac{1}{\max(\frac{\xi_1}{k_1},
\frac{\xi_2}{k_2})} \label{eq:geo-series-0}
\EEQ
Extending to $\alpha>0$ and deficiency weight $\eps:=\frac{\kappa_1}{|A_{1,1}|} = \frac{\nu_+}{|A_{1,1}|} \sim \frac{\nu_+}{k_1}>0$, we get  for $\alpha$ large enough so that the series
converges
\BEA  \sum_{n=0}^{+\infty} (\Phi(\eps,\alpha)_{\cal G})^n &:=& \frac{1}{1-w(\alpha)_{1\to 2} w(\alpha)_{2\to 1}} =  \frac{1}{1- \frac{k_1+2\nu_+}{k_1+\xi_1 + \nu_+ +\alpha} 
\frac{k_2}{k_2+\xi_2+\alpha}}  \nonumber\\
&\sim&  \frac{1}{1- (1-\frac{\xi_1+\alpha}{k_1}) 
(1-\frac{\xi_2+\alpha}{k_2}) - \eps} \nonumber\\
&  \sim &   \frac{1}{\max(\frac{\xi_1+\alpha}{k_1},
\frac{\xi_2+\alpha}{k_2}) - \eps} \sim  
\frac{1}{\max(\frac{\xi_1}{k_1},
\frac{\xi_2}{k_2}) + \frac{\alpha}{k_{min}} - \eps}
\label{eq:geo-series-1}
\EEA

\Medskip 
\begin{Definition}[weight renormalization factor of ${\cal G}$)]  \label{def:cycle-Z}
 Let 
\BEQ Z(\eps,\alpha)_{{\cal G}} := \frac{1}{ \sum_{n=0}^{+\infty} (\Phi(\eps,\alpha)_{{\cal G}})^n}\EEQ
\end{Definition}

From the above,
\BEQ Z(\eps,\alpha)_{{\cal G}} \sim Z(0)_{\cal G} - \eps + \frac{\alpha}{k_{min}}.  \EEQ
where
\BEQ Z(0)_{{\cal G}} := Z(0,0)_{{\cal G}} \sim   \max(\frac{\xi_1}{k_1},
\frac{\xi_2}{k_2})
\EEQ
may be interpreted as something like the probability
of leaving ${\cal G}^{int}$ each time one goes around the cycle (see below), and represents a mass loss ratio. {\em By assumption, both $Z(0)_{\cal G}\prec 1$ and $\eps\prec 1$ are small factors}; this is  a key fact. Then,   
if $\alpha\prec k_{min}$, $Z(\eps,\alpha)_{\cal G}\prec 1$ is also small.

\Medskip We let the {\em weight} of $G$
be
\BEQ Z^{-1}_G := (\max(Z(0)_{\cal G}, \eps))^{-1}. \EEQ
Finally, the {\em Lyapunov exponent} $\lambda_{\cal G}$ is by definition the value of $\alpha$ such
that $Z(\eps,\alpha)_{\cal G}\equiv 1$ since this is equivalent to the condition  (\ref{eq:excursion})
on the total excursion weight. Thus 
\BEQ \lambda_{\cal G} \sim k_{min}(\eps-Z(0)_{\cal G}) \EEQ
when the two terms do not compensate. The main focus is on the quantitative estimate obtained in the autocatalytic case: when
$\eps\succ Z(0)_{\cal G}$, $\lambda_{\cal G}\sim k_{min}\eps$ is $>0$.

\Medskip We get three easy regimes: 

\qquad\qquad \textbullet\ {\em (free regime)} Assume $\frac{\alpha}{k_{min}}\ ,\eps \prec Z(0)_{\cal G}$. Then $Z(\eps,\alpha)_{\cal G} \sim Z(0)_{\cal G}$ is $>0$, and so is  $Z^{-1}_G \sim (Z(0)_{\cal G})^{-1}$.

\qquad\qquad \textbullet\ {\em (autocatalytic regime)} Assume $\eps \succ Z(0)_{\cal G}, \frac{\alpha}{k_{min}}$ (in particular, $\cal G$ is autocatalytic: $\lambda_{\cal G}>0$). Then 
$Z(\eps,\alpha)_{\cal G} \sim -\eps$ is $<0$, but $Z^{-1}_G \sim \eps^{-1}$ is $>0$.

\qquad\qquad \textbullet\ {\em (degraded regime)} Assume $\frac{\alpha}{k_{min}} \succ 
Z(0)_{\cal G},\eps$. Then $Z(\eps,\alpha)_{\cal G}\sim \frac{\alpha}{k_{min}} \succ Z_G$.   

\Medskip Then there are a priori two delicate  boundary cases,

\qquad\qquad \textbullet\   $ Z(0)_{\cal G}\sim \eps\succ \frac{\alpha}{k_{min}}$. In this 
case, the mass loss ratio (due to outgoing edges) is approximately compensated by  the mass increase ratio represented by the  deficiency weight (due to internal autocatalysis). This is a {\bf resonance case}, as discussed in \S \ref{subsection:aim}. Excluding this case where our algorithm is inconclusive, we also have the case when

\qquad\qquad \textbullet\ 
  $\eps\sim \frac{\alpha}{k_{min}}
\succ Z(0)_{\cal G}$. Then the mass increase ratio, and the term 
$\frac{\alpha}{k_{min}}$ due to the overall degradation rate $\alpha$, potentially compensate.  
The weight renormalization factor $Z(\eps,\alpha)_{\cal G}$ is never used as such in the estimates; rather (see discussion of {\em equivalent transition rates} below), it is the
ratio of weights $w(\alpha)_{G\to\sigma_{ext}}/w_{\sigma\to\sigma_{ext}}$, see (\ref{eq:walphaGsigma'Z0alpha}),   and  (\ref{eq:walphaGsigma'Z}), that plays a role. So there 
is actually no issue here.

\Bigskip Thinking of the concentration vector $X(t)$ as the time-evolution of the Markov 
chain as discussed in \S \ref{subsection:aim}, and using the Gillespie algorithm to generate the
continuous-time chain from the discrete one, suggests the following intuitive elements.

\Medskip {\em Characteristic time of ${\cal G}$.} We let 
\BEQ \tau_{{\cal G}} \sim 1/k_{min} \EEQ
Namely, each cycle of the Markov chain takes an average time $\sim \tau_{{\cal G}}$, whence the name. The scale of the {\em characteristic rate}, $1/\tau_{{\cal G}}\sim k_{min}$, is equal to the
lowest scale $n_{min}$.

\Medskip {\em External scale of ${\cal G}$ .}  We let
\BEQ n^{ext}_{{\cal G}} = \lfloor \log_b k^{ext}_{{\cal G}}\rfloor, \qquad 
\ k^{ext}_{{\cal G}}\sim \frac{1}{\tau_{{\cal G}}} \, \times\,  Z(0)_{{\cal G}} \label{eq:n_G}
\EEQ
Namely, the probability to leave the cycle from $\sigma$ (in the Gillespie realization) is 
 $\sim\frac{\xi_{\sigma}}{k_{\sigma}}$; since the discrete-time chain is alternatively in either
 state $\sigma=1,2$, one  leaves ${\cal G}^{int}$ after going around the cycle  $\sim 1/Z(0)_{{\cal G}}$ times in 
 average. Since each cycle takes a time $\sim\tau_{{\cal G}}$, $1/k^{ext}_{\cal G}$ represents
 the average time is takes to leave the cycle. Note the following (general, as turns out) fact: {\em the external scale of ${\cal G}^{int}$ is
strictly lower than its lowest scale.}

\Bigskip Next, we would like to define renormalized transition rates/weights for the 
network obtained after collapsing the cycle, i.e. merging all the internal vertices into
a compound vertex denoted $G$.

\Medskip {\em Equivalent transition rates.} Let
\BEQ k_{G\to \sigma_{ext}}\sim \frac{1}{\tau_{{\cal G}}}  \, \times\, 
 \frac{\xi_{\sigma}}{k_{\sigma}}  \label{eq:ETR}
 \EEQ
Namely, the average time needed to make one cycle is $\sim \tau_{{\cal G}}$. Then, if the chain 
state is $\sigma$, the chain exits 
to $\sigma_{ext}$ with probability $\frac{\xi_{\sigma}}{k_{\sigma}}$.   Note that 
\BEQ k^{ext}_{\cal G} \sim \sum_{\sigma=1,2} k_{G\to\sigma_{ext}} \EEQ
Thus $k^{ext}_{\cal G}$ may be thought as the total outgoing rate $k_G$ of the compound
vertex $G$,
in coherence with the general definition (\ref{eq:k}).

\Medskip {\em Equivalent transition weights} $w(\alpha)_{G\to \sigma_{ext}}$ are obtained 
from (\ref{eq:ETR}) by dividing by $k_G+\alpha\sim k^{ext}_{\cal G}+\alpha$. 
In the free regime,  $k^{ext}_{\cal G}+\alpha\sim k^{ext}_{\cal G}$, yielding 
\BEQ w(\alpha)_{G\to \sigma_{ext}} \sim w(0)_{G\to \sigma_{ext}} \sim Z(0)_{{\cal G}}^{-1} 
\, \times\, 
 \frac{\xi_{\sigma}}{k_{\sigma}}.  \qquad {\mathrm{(free \ regime)}} 
 \EEQ
In the autocatalytic/degraded regimes, the
 prefactor $Z(0)_{{\cal G}}^{-1} $ is replaced by $\sim (\alpha \tau_{{\cal G}})^{-1}$.

\Medskip {\em Equivalent deficiency rate.} By a reasoning similar to (\ref{eq:ETR}), we get
\BEQ \kappa_G \sim \frac{\eps}{\tau_{{\cal G}}}.
\EEQ
Then the {\em renormalized deficiency  weight}    is
\BEQ \eps_G \sim \frac{\kappa_G}{k_G} \sim \frac{\kappa_G}{k^{ext}_{\cal G}} \sim \eps 
\ \times\ \frac{k_{min}}{k^{ext}_{\cal G}} \sim \frac{\eps}
{Z(0)_{{\cal G}}}, \label{eq:epsG'}
\EEQ
equal to $\eps$,  renormalized by the large
factor $1/Z(0)_{{\cal G}}$ due to the replacement  of the lowest scale
by the external
scale in the denominator.

\vskip 2cm\noindent   We report in the Table below the most important expressions obtained so far,
depending on the regime.  

\Bigskip 
\begin{center}
\begin{tabular}{c|c|c|c} \label{table:regimes}
Regime & free &  autocatalytic & degraded \\ \hline
Conditions & $Z(0)_{{\cal G}}\succ \alpha \tau_{\cal G},\eps$ & 
$\alpha \tau_{{\cal G}} \sim \eps \succ Z(0)_{{\cal G}}$ & 
$\alpha\tau_{{\cal G}} \succ Z(0)_{{\cal G}}, \eps $ \\
   \hline
Value of $Z(\eps,\alpha)_{{\cal G}}$ & $Z(0)_{{\cal G}}$ & $\alpha \tau_{{\cal G}}$ & $\alpha\tau_{{\cal G}}$ \\
Value of $w(\alpha)_{G\to \sigma_{ext}}$ & $Z(0)_{{\cal G}}^{-1} 
 \frac{\xi_{\sigma}}{k_{\sigma}}
 $ & $\frac{1}{\alpha\tau_{{\cal G}}}  \frac{\xi_{\sigma}}{k_{\sigma}}$
& $\frac{1}{\alpha\tau_{{\cal G} }} \frac{\xi_{\sigma}}{k_{\sigma}}$ \\ \hline
Value of $k_{G\to \sigma_{ext}}$ & $\frac{1}{\tau_{{\cal G}}} 
 \frac{\xi_{\sigma}}{k_{\sigma}}$ & 
$\frac{1}{\tau_{{\cal G}}}  \frac{\xi_{\sigma}}{k_{\sigma}}
$
  & 
$\frac{1}{\tau_{{\cal G}}} 
 \frac{\xi_{\sigma}}{k_{\sigma}}   $ 
\end{tabular}

\bigskip {\textsc{TABLE \ref{table:regimes} -- The three different regimes for a 
cycle.}}
\end{center}

\Bigskip Here $\sigma=1,2$, but the extension to cycles $1\to 2\to\cdots\to p\to 1$ is straightforward.

\Medskip We get the following renormalized graph,

\bigskip
\begin{center} 
\begin{tikzpicture}
\begin{scope}[shift={(10,0)}]

\draw[dotted](1.5,-1)--(1.5,-2.25);

\draw(1.5,-0.5) node {$G$};
\draw[->, ultra thick](1.5,-1)--(2.5,-1);  \draw(2.5,-0.5) node {$1_{ext}$};
\draw(5,-1) node {$k_{G\to 1_{ext}} \sim k_{min} \frac{\xi_1}{k_1}  $};

\draw[-<](1.25,-1) arc(20:340:0.4 and 0.25) ; \draw(0,-1) node {$\kappa_G$};

\draw[dashed](0,-1.5)--(5,-1.5);

\draw(1,-1.85) node {$2_{ext}$};
\draw[->](1.5,-2.25)--(1,-2.25);
\draw(4.5,-2.25) node {$k_{G\to 2_{ext}} \sim k_{min} \frac{\xi_2}{k_2} \sim \xi_2$};

\end{scope}
\end{tikzpicture}
\end{center}

\Bigskip if $\xi_1/k_1\succ \xi_2/k_2$, resp.

\bigskip
\begin{center} 
\begin{tikzpicture}
\begin{scope}[shift={(10,0)}]

\draw[dotted](1.5,-1)--(1.5,-2.25);

\draw(1.5,-0.5) node {$G$};
\draw[->, ultra thick](1.5,-1)--(2.5,-1);  \draw(2.5,-0.5) node {$2_{ext}$};

\draw[-<](1.25,-1) arc(20:340:0.4 and 0.25) ; \draw(0,-1) node {$\kappa_G$};

\draw[dashed](0,-1.5)--(5,-1.5);

\draw(1,-1.85) node {$1_{ext}$};
\draw[->](1.5,-2.25)--(1,-2.25);

\draw(4.8,-1) node {$k_{G\to 2_{ext}} \sim \xi_2$};

\draw(4.5,-2.25) node {$k_{G\to 1_{ext}} \sim k_{min} \frac{\xi_1}{k_1} $};

\end{scope}
\end{tikzpicture}
\end{center}

\Bigskip  if  $\xi_1/k_1\prec \xi_2/k_2$.  It is assumed implicitly that we are in the free case,
for then, the new edge $G\to 1_{ext}$, resp. $G\to 2_{ext}$ is now dominant, and consequently, drawn
in boldface. In the autocatalytic, resp. degraded
case, the scale of this new edge is dominated either by the scale of $\kappa_G$ or by $\alpha$.


\subsection{Heuristics  for the merging procedure}  \label{subsection:heuristics}


\medskip We refer the reader back to Section \ref{section:notations} for notations. Let ${\cal G}^{int} =
(V^{int}, E^{int})$ be a cut-off graph connected by dominant edges, and ${\cal G}$ its completion 
by outgoing edges. The generator $A^{int} = A|_{{\cal G}^{int}}$ is not probability-preserving, but can be made so by changing the diagonal coefficients. We call $\tilde{A}^{int}$ the modified generator, which is a Markov chain generator, $\tilde{\cal W}^{int}$ the Markov transition matrix, and  $\tilde{\pi}^{int}$  the associated stationary probability measure. 

\Medskip We generalize here the quantities introduced in \S \ref{subsection:cycle-elementary}, based
on an expansion of the resolvent formula, and provide a bridge to the renormalization algorithm
of \S \ref{subsection:description}. This subsection is not rigorous; the proofs developed in Section
\ref{section:proofs} rely instead on partial resummations of the path series.

\Bigskip {\bf Characteristic rate/scale of ${\cal G}$.} As well-known from the
theory of finite Markov chains, $\tilde{\pi}^{int}_{\sigma}$ represents the
average fraction of instants spent in state $\sigma$ in a long
discrete-time trajectory (of the discrete-time Markov chain with generator $\tilde{\cal W}^{int}$) started from any point. Renormalizing
by the average waiting times $(1/k_{\sigma})_{\sigma}$, we get
the equivalent quantities for the time-continuous Markov chain;
namely, its stationary measure  is $\mu_{\sigma}=\frac{1}{\tau_{\cal G}}
\frac{\tilde{\pi}_{\sigma}}{k_{\sigma}}$, where the normalization 
constant 
\begin{equation} \tau_{{\cal G}}=\sum_{\sigma} \frac{\tilde{\pi}_{\sigma}}{k_{\sigma}} 
\equiv \langle \frac{1}{k_{\sigma}} \rangle_{\tilde{\pi}} \end{equation}
  an averaged waiting time,  is a
 characteristic time for ${\cal G}^{int}$. Its inverse $ 1/\tau_{\cal G}$ is called {\em characteristic rate}, see (\ref{eq:char-rate}). 

\medskip\noindent  {\bf Deficiency weight of $\cal G$.} Let $\bar{\eps}_{\cal G}:=\sum_{\sigma\in V^{int}}
\tilde{\pi}^{int}_{\sigma} \eps_{\sigma} = \langle \eps\rangle_{\tilde{\pi}}$ be the average deficiency weight of ${\cal G}$.
This may be interpreted as an average growth rate per time unit
$\tau_{{\cal G}}$.

\medskip\noindent {\bf Resolvent formula} (see \S \ref{subsection:resolvent}). Let $\alpha>\lambda^*(A)$, then 
{\small \begin{equation}
((\alpha \Id - A)^{-1})_{\sigma',\sigma} = \sum_{\gamma:\sigma\to\sigma'} 
\Big(\prod_{k} w(\alpha)_{\sigma_k\to\sigma_{k+1}}\Big) \, \times\, \frac{1}{|A_{\sigma',
\sigma'}|+\alpha}.  \label{eq:resolvent}
\end{equation} }
is a sum over all paths $\gamma$ connecting $\sigma$ to $\sigma'$ (the condition on $\alpha$ ensures that  the series is convergent). This
classical formula \cite{RevYor}, \cite{NgheUnt1} is obtained by rewriting the matrix coefficient as $\int_0^{+\infty} dt\, (e^{tA(\alpha)})_{\sigma',\sigma}$ and using Trotter's formula to separate
transition-defining off-diagonal coefficients from the diagonal, which is exponentiated and
integrated out.  

\medskip\noindent  We note that (see (\ref{eq:5.14}) for details)
\begin{eqnarray} && w(\alpha)_{\sigma\to \sigma'}  = \frac{k_{\sigma\to\sigma'}}{|A_{\sigma,\sigma}|+\alpha}  = \frac{k_{\sigma}}{|A_{\sigma,\sigma}|} \, \times\,  \frac{|A_{\sigma,\sigma}|}{|A_{\sigma,\sigma}|+\alpha} \, \times\, 
\frac{|\tilde{A}^{int}_{\sigma,\sigma}|}{k_{\sigma}}  \nonumber\\
&& \qquad  \times\, 
\frac{k_{\sigma\to\sigma'}}{|\tilde{A}^{int}_{\sigma,\sigma}|}
\sim \Big(1+ \eps_{\sigma} - \frac{k^{ext}_{\sigma} + \alpha}{k_{\sigma}}\Big)\, 
\tilde{\cal W}^{int}_{\sigma,\sigma'}.  \label{eq:wsigmaalpha-pert}
\end{eqnarray}
when $\eps_{\sigma}$ (because $\sigma$ is not autocatalytic), $\frac{k_{\sigma}^{ext}}{k_{\sigma}}$ (because $\Sigma^{int}$ is a dominant SCC) and $\frac{\alpha}{k_{\sigma}}$ 
are $\prec 1$ and can be treated perturbatively. 
Replacing $w(\alpha)_{\sigma\to\sigma'}$ by $\tilde{\cal W}^{int}_{\sigma,\sigma'}$ in 
(\ref{eq:resolvent}) would lead to a divergence since we get (by resummation) a geometric series with parameter $ \sum_{\gamma^*:\sigma\to\sigma} \prod_{k} \tilde{w}^{int}_{\sigma_k\to\sigma_{k+1}}=1$ 
equal to a sum over {\em excursions} $\gamma^*:\sigma=\sigma_1\to\cdots\to \sigma_{\ell+1}=\sigma$, i.e. paths such that $\sigma_i\not=\sigma$, $2\le i\le \ell$; the above expression is
1 because it can be interpreted as the return probability to $\sigma$ for a finite, irreducible Markov chain.
Summing over excursions in $\Sigma^{int}$ and multiplying factors as in (\ref{eq:wsigmaalpha-pert})  yields 1, minus (to leading order) a $O(1)$ weight factor $Z(\eps,\alpha)_{\cal G}$ that we
now discuss. If   $Z(\eps,\alpha)_{\cal G}>0$, then the geometric series $\sum_{\ell\ge 0} 
(1-Z(\eps,\alpha)_{\cal G})^{\ell}$ converges to $\frac{1}{Z(\eps,\alpha)_{\cal G}}$. If 
it is $<0$, then the series diverges, which means that $\alpha\prec \lambda^*(A^{int})$. This
suggests to estimate $\lambda^*(A^{int})$  as the threshold value $\alpha_{thr} := \inf\{\alpha>0 \ |\ Z(\eps,\alpha)_{\cal G}>0\}$. Rigorous arguments (see Section \ref{section:proofs} and Appendix) rely instead on the fixed point theorem and a classical gap estimate for Markov chains attributed to W. Doeblin.  

\medskip\noindent  {\bf Weight renormalization of a subgraph ${\cal G}$.} 
The weight renormalization is the average with respect to
stationary measure $\tilde{\pi}^{int}$  of the outgoing transition weights from 
vertices $\sigma$ in $\Sigma^{int}$, namely (see (\ref{eq:Z(0)}) or Lemma \ref{lem:Lya(1)})
\begin{equation} Z(0)_{{\cal G}} = \sum_{\sigma\in V^{int}}
\tilde{\pi}_{\sigma} \sum_{\sigma'\in V^{ext}}   \frac{k_{\sigma}^{ext}}{k_{\sigma}}, \qquad
k_{\sigma}^{ext} = 
 \sum_{(\sigma,\sigma')\in E^{out}} k_{\sigma\to\sigma'} 
\label{eq:6.65}
\end{equation} 
Expanding the sum over excursions yields rather the  $(\eps,\alpha)$-corrected weight 
\begin{eqnarray}
 && Z(\eps,\alpha)_{\cal G} = \sum_{\sigma\in V^{int}}
\tilde{\pi}_{\sigma} \sum_{\sigma'\in V^{ext}}  \Big( \frac{k_{\sigma}^{ext} + \alpha}{k_{\sigma}} - \eps_{\sigma}\Big)  \nonumber\\
&& = Z(0)_{\cal G} + \tau_{{\cal G}}\alpha - \bar{\eps}_{\cal G}
\end{eqnarray}
where 
\BEQ \bar{\eps}_{\cal G}= \sum_{\sigma \in V^{int}} \tilde{\pi}^{int}_{\sigma} \eps_{\sigma} \EEQ
is an average deficiency weight, or an average growth rate per time unit $\sim \tau_{\cal G}$. 
Hence $\cal G$ is autocatalytic when $\bar{\eps}_{\cal G} \succ Z(0)_{\cal G}$, and then 
$\alpha_{thr} \sim \lambda_{{\cal G}}\sim   \bar{\eps}_{{\cal G}}/\tau_{{\cal G}_{int}}$ (see (\ref{eq:lambdaGp})).

\medskip\noindent  {\bf External transition rate of $\cal G$.} Define (see (\ref{eq:ext-rate}))
\begin{equation} k^{ext}_{{\cal G}} \sim 
\frac{1}{\tau_{G}} \, \times\, Z(0)_{{\cal G}}.
\label{eq:kextG}
\end{equation}
 The
rate $k^{ext}_{{\cal G}}$, in the product form (characteristic rate) $\times $ 
(average transition weight),  represents the {\em external
transition rate} from ${\cal G}^{int}$. 

\medskip\noindent  {\bf Transition weights/rates from $G$.} The transition rate from $G$,
\begin{equation} k_{G\to \sigma'} \sim  \frac{1}{\tau_{\cal G}} \sum_{\sigma \in V^{int} }
\tilde{\pi}_{\sigma} \frac{k_{\sigma\to\sigma'}}{k_{\sigma}} 
\end{equation}
see (\ref{eq:ren-ingoing-rate}), represents an average transition probability from $G$ 
to $\sigma'$ per time unit $\tau_{\cal G}$, multiplied by the characteristic rate $1/\tau_{\cal G}$. Note that the total outgoing rate from $G$ is
\BEQ k_G \sim \sum_{\sigma'\in V^{ext}} k_{G\to \sigma'} \sim k^{ext}_{\cal G} 
\EEQ 
Then, dividing by the external transition rate $k^{ext}_{\cal G}$, we get

\begin{equation} w(0)_{G\to \sigma'} \sim \frac{k_{G\to\sigma'}}{k^{ext}_{\cal G}} \sim   \frac{1}{Z(0)_{{\cal G}}} \sum_{\sigma \in V^{int} }
\tilde{\pi}_{\sigma} \frac{k_{\sigma\to\sigma'}}{k_{\sigma}} .
\end{equation}
Letting $\sigma^*$ maximize the set $\{\frac{k_{\sigma\to\sigma'}}{k_{\sigma}}, \sigma \in V^{int}\}$,
 $w(0)_{G\to \sigma'} \sim (Z(0)_{\cal G})^{-1} \, \times\, w_{\sigma^*\to\sigma'}$
is renormalized by the prefactor $(Z(0)_{\cal G})^{-1}$. For all useful $\alpha$, see remark below
(\ref{eq:w(alpha)-free-autocata-regimes}), this is generalized as 
\BEQ w(\alpha)_{G\to \sigma'} \sim 
(Z(0,\alpha)_{\cal G})^{-1} \, \times\, w_{\sigma^*\to\sigma'}.  \label{eq:walphaGsigma'Z0alpha}
\EEQ  

\medskip\noindent  {\bf Renormalized deficiency weight of $G$} (see (\ref{eq:ren-def-weight})).  Let $\eps_{G} = \frac{\bar{\eps}_{{\cal G}}}{Z(0)_{{\cal G}}}$. 
Comparing with $\bar{\eps}_{{\cal G}}$, this new quantity is  an average growth rate per time unit
$1/k^{ext}_{{\cal G}} \sim \tau_{{\cal G}}/ Z(0)_{{\cal G}}$.

\medskip\noindent {\bf Lyapunov weights and Green functions.} We sketch a proof of formulas (\ref{eq:pisigma1},
 \ref{eq:pisigma2}).  Let $P:={\cal W}(\lambda^*)$ be
the discrete-time transition matrix on $V(G)$ with coefficients $P_{\sigma,\sigma'}= w(\lambda^*)_{\sigma\to\sigma'}$. For $\sigma^*\in V(G)$ fixed, 
\BEQ G^{\sigma^*}_N(\sigma'):=
\sum_{\ell =0}^{N-1} (P^{\ell})_{\sigma^*,\sigma'} \label{eq:Green0}
\EEQ
 is a truncated {\bf Green kernel}. Both $\lim_{N\to\infty}
\frac{1}{N}G^{\sigma^*}_N$ and $\pi^*$ (see (\ref{eq:pi*})) are left eigenvectors of $P$ with eigenvalue $1$, hence they are equal up to normalization (see \S \ref{subsubsection:full-auto} for details).  Expanding $\frac{1}{N}G^{\sigma^*}_N$
gives a C\'esaro summation of the sequence $(P^{\ell})_{\sigma^*,\sigma}$, whose general term is equal to the  sum
of the weights of paths $\gamma:\sigma^*\to \sigma$ of length $\ell$.  We choose  a vertex $\sigma^*\in \Sigma^{(0)}$ (the original vertex set)
belonging to one of the cores (i.e. maximal dominant SCCs with highest deficiency scale, see 
\S \ref{section:description} for details), and of highest scale. Then the total weight
$w(\lambda^*)_{\gamma'}$
of paths $\gamma':\sigma^*\to \sigma^*$ of lengths inside a finite interval centered around 
$\ell$ is $\sim 1$ for $\ell$ large enough. Then a typical path contributing to the C\'esaro summation will be
a repetition of combinations of paths $\gamma'$, followed by a path $\gamma$ from $\sigma^*$ to $\sigma$
which does not return to $\sigma^*$. The procedure may be iterated while trickling down from $\sigma^*$ 
by finite, non self-intersecting paths $\sigma^*\to \sigma_1\to\cdots\to \sigma_q$ obtained
by concatenating paths in  $G^{(i)}(\sigma^*)\setminus G^{(i-1)}(\sigma^*)$, 
(excluding repetitions); see \S \ref{subsection:merging}. Dominant path fragments circulating within $G^{(i)}(\sigma^*)\setminus G^{(i-1)}(\sigma^*)$ have weight $\sim 1$. On the other hand, steps $\sigma_q\to\sigma_{q+1}$ 
going down a coalescence tree, i.e. from $G^{(i)}(\sigma^*)\setminus G^{(i-1)}(\sigma^*)$ to 
$G^{(i+1)}(\sigma^*)\setminus G^{(i)}(\sigma^*)$, face  weight renormalization by a 
renormalization prefactor $(Z(0,\alpha)_{\cal G})^{-1}$ as in (\ref{eq:walphaGsigma'Z0alpha}) . Conserving the original relative
normalization is ensured by the product structure (\ref{eq:pisigma1}) along the sequence of 
downsteps. For $\sigma$ not in the terminal vertex $G^{(\infty)}(\sigma^*)$, the relative prefactor is 
simply $1$, and the product on the second line of (\ref{eq:pisigma2})  is obtained by
computing the products of factors of the type  $(Z(0,\alpha)_{\cal G})^{-1}$ obtained by concatenating excursions as in Definition \ref{def:cycle-Z}.
On the other hand,  paths $\gamma$ with maximal weight in the first line of (\ref{eq:pisigma2}) may be found
following the DAG construction procedure of Appendix \ref{subsection:DAG}.


\section{Description of the renormalization algorithm}  \label{section:description}


This section  contains no proofs, however, the logic should be accessible to the reader
through the heuristic subsection \S \ref{subsection:heuristics}. All proofs are to be found 
in Section \ref{section:proofs}. We provide two explicit examples at the end.


\subsection{Description}  \label{subsection:description}


\Medskip The starting point is the split graph $G^{(0)}$ with reaction rates $k^{(0)}_{\sigma\to\sigma'} = k_{\sigma\to\sigma'}, \beta^{(0)}_{\sigma} = \beta_{\sigma} =  k_{\sigma\to\emptyset}$ obtained after linearizing at 0. Call $n^{(0)}_1> n^{(0)}_2>\cdots$ the reaction scales of $G^{(0)}$. We start from the cut-off scale $n^{(0)}_1$, i.e.
eliminate  reactions with scale $< n^{(0)}_1$, so that all reactions in the infra-red cut-off graph 
$G: = G^{(0)}_{\searrow n^{(0)}_1}$   (see Section \ref{section:notations})
have same scale $n^{(0)}_1$; in particular, they are all dominant. If $G$ contains no non-trivial maximal
dominant SCC, then we replace the cut-off scale by $n^{(0)}_2$. 
The cut-off graph $G = G^{(0)}_{\searrow n^{(0)}_2}$ now contains reactions of both scales 
$n^{(0)}_1,n^{(0)}_2$. We consider {\em only} those reactions which are dominant, yielding a 
dominant subgraph, still denoted $G$. If $G$ contains no non-trivial maximal
dominant SCC, then we replace the cut-off scale by $n^{(0)}_3$, and so forth.

\Medskip\noindent {\bf One renormalization step $(i=1)$.} The process usually stops at some scale $n = n(i)$ (here $i=1$), called {\em step $i$ cut-off scale}, at which  a  non-trivial maximal
dominant SCC appears. By construction, $n(i)\equiv n^{(i-1)}_{j_i}$ is one of the scales of 
the previous step graph. We first let $G^{(i-1)}_{cut} := G^{(i-1)}_{\searrow n^{(i-1)}_{j_i-1}}$ be the graph cut off at the scale just above $n$; this cut-off graph has no maximal dominant SCC.      Then (adding scale $n$) we let $(G_p)_{p=1,2,\ldots}$ be the non-trivial maximal dominant SCCs of $G_{\searrow n}$.  By construction, they involve only dominant  edges, and their minimal scale is $n$; also, none of the vertices of $G_p$ is autocatalytic, because there is no dominant edge $\sigma\to\sigma'$, $\sigma'\not=\sigma$ when $\sigma$ is autocatalytic.
 Fix $p$. Adding to the internal edges 
in $G_p$    the set  ${\cal E}^{out}_p$ of  outgoing edges $k_{\sigma\to\sigma'}$, $\sigma\in G_p, \sigma\not\in G_p$ yields a
new graph, ${\cal G}_p$, which is the same as the cut-off graph $G^{int}$, except for the fact that the target of outgoing edges is specified.   By construction, {\em none of the outgoing edges  is dominant}.  
{\em We shall now collapse  ${\cal G}_p$, $p=1,2,\ldots$ inside $G^{(0)}$}. Graphically, one (1) merges all vertices $\sigma$ in 
$G_p$ into a compound vertex denoted $G_{i=1,p}$, so that $G^{(i)}(\sigma) = G_{i,p}$ in the
notation of \S \ref{subsection:merging}; for all other vertices $\sigma$ (including 
those in non-maximal dominant SCCs), 
the map $G^{(i)}(\sigma)=\sigma$ is trivial; (2) 
redirects every outgoing edge $\sigma\to\sigma'$ in ${\cal E}^{out}_p$ $(\sigma\in V(G_p), 
\sigma'\not\in V(G_p))$ into an edge $G_{i,p}\to G^{(i)}(\sigma')$ from $G_{i,p}$; (3) and further, 
redirects every ingoing edge $\sigma'\to\sigma$  $(\sigma\in V(G_p)$, 
$\sigma'$ 'trivial', i.e. $G^{(i)}(\sigma')=\sigma'$) into an edge $\sigma'\to G_{i,p}$.  We need now 
specify the kinetic rates of the new coarse-grained graph $G^{(i)}$, called {\bf step $i$ effective graph}.  Define first

\begin{equation}  {\mathit{(characteristic\ rate)}} \qquad  1/\tau_{{\cal G}_p} = k_{p,min} \label{eq:char-rate} 
\end{equation}
 where $k_{p,min}$ is any internal
rate with scale $n(i)$, equal to the lowest scale of all internal edges in $G_p$;

\begin{equation} {\mathit{(bare\ deficiency\ weight)}}  \qquad \bar{\eps}_{{\cal G}_p} \sim \max\{\eps_{\sigma},\ \sigma\in V(G_p)\};  \label{eq:bare-def-weight}
\end{equation}
\begin{equation} {\mathit{(renormalization\ factor)}}  \qquad 
 Z(0)_{{\cal G}_p} \sim \max_{(\sigma,\sigma') \in {\cal E}_p^{out}}
\frac{k_{\sigma\to\sigma'}}{k_{\sigma}}  \label{eq:Z(0)}
\end{equation}
which is $\prec 1$ since outgoing edges are non-dominant; 
\begin{equation}   {\mathit{ (external\ rate)}}  
\qquad  k^{ext}_{{\cal G}_p} \sim \frac{1}{\tau_{{\cal G}_p}} \ \times\ 
Z(0)_{{\cal G}_p}   \label{eq:ext-rate}
\end{equation}
The {\bf resonance regime} is defined by $Z(0)_{{\cal G}_p}\sim \bar{\eps}_{{\cal G}_p}$; 
in that regime, we cannot decide whether ${\cal G}_p$ is autocatalytic or not. Banning
this regime,
 the {\em Lyapunov\ exponent of} ${\cal G}_p$ is defined as
\begin{equation}   \lambda_{{\cal G}_p} \sim
\begin{cases} 0, \qquad   Z(0)_{{\cal G}_p} \succ \bar{\eps}_{{\cal G}_p} \\
 \bar{\eps}_{{\cal G}_p}/\tau_{{\cal G}_p} , \qquad   Z(0)_{{\cal G}_p} \prec \bar{\eps}_{{\cal G}_p} 
 \end{cases}     \label{eq:lambdaGp}
\end{equation}

When $Z(0)_{{\cal G}_p} \succ \bar{\eps}_{{\cal G}_p}$ {\em (non-autocatalytic case)}, the Lyapunov exponent of $G^{int}_p$ is actually $\prec 0$, but this may change when ingoing edges are taken into account at step $(i+1)$. 
However, when $Z(0)_{{\cal G}_p} \prec \bar{\eps}_{{\cal G}_p}$ {\em (autocatalytic case)}, the exponent  is $\succ 0$, 
and adding ingoing edges can only increase it. 

\Medskip The renormalized rates of the compound vertex $G_{i,p}$ are now

\begin{equation}  {\mathit{(outgoing\ \ rates)}}   \qquad
 k_{G_{i,p}\to\sigma'} \sim \frac{1}{\tau_{{\cal G}_p}} \, \times\, \max_{\sigma \in V(G_p)}  \frac{k_{\sigma\to\sigma'}}{k_{\sigma}}  \label{eq:ren-rates}
\end{equation}

\begin{equation}  {\mathit{(deficiency\ rate) }}\qquad   \kappa_{G_{i,p}}  \sim \bar{\eps}_{{\cal G}_p}/\tau_{{\cal G}_p} \label{eq:ren-def-rate}
\end{equation}
\noindent giving rise to a self-edge $G_{i,p}\to G_{i,p}$, interpreted as a doubling reaction 
$G_{i,p}\overset{\kappa_{G_{i,p}}}{\to} G_{i,p}+G_{i,p}$;  note that,  in the autocatalytic case,  $\lambda_{{\cal G}_p}\sim
\kappa_{G_{i,p}} \succ k_{G_{i,p}}$ is dominant;

\begin{equation} {\mathit{(ingoing\ rates)}}  \qquad  k_{\sigma' \to 
G_{i,p}} \sim  \max_{\sigma \in V(G_p)}  k_{\sigma'\to\sigma}.  \label{eq:ren-ingoing-rate}
\end{equation}

\Medskip Finally, we call {\em weight of $G_{i,p}$} the factor
\begin{equation} {\mathit{(weight)}}\qquad Z^{-1}_{G_p} \sim \Big( 
\max(Z(0)_{{\cal G}_{i,p}}, \bar{\eps}_{{\cal G}_p}) \Big)^{-1}   \label{eq:ren-weight}
\end{equation}
and let
\begin{equation} Z(\eps,\alpha)_{{\cal G}_p}:= Z(0)_{{\cal G}_p} - \bar{\eps}_{{\cal G}_p} + \alpha \tau_{{\cal G}_p}, \qquad \eps,\alpha\ge 0  \label{eq:Zepsalpha}
\end{equation}
Note that $Z(\eps,0)_{{\cal G}_p} \sim Z(0)_{{\cal G}_p}$  in the free regime. When $\alpha\succ Z_{G_p}/ \tau_{{\cal G}_p}$, $Z(\bar{\eps}_{{\cal G}_p},\alpha)_{{\cal G}_p}
\sim Z(0,\alpha)_{{\cal G}_p} \sim \alpha\tau_{{\cal G}_p}$ simply.

\Bigskip Looking at the new graph, we see that, by construction, 
\BEQ
 k_{G_{i,p}} = 
\sum_{\sigma'\not= G_{i,p}} k_{G_{i,p}\to \sigma'} \sim k^{ext}_{{\cal G}_p},
 \EEQ
  from which we get
the renormalized deficiency weight, 
\begin{equation} 
\eps_{G_{i,p}}\sim \kappa_{G_{i,p}}/k_{G_{i,p}} \sim \frac{\bar{\eps}_{{\cal G}_p}}{Z(0)_{{\cal G}_p}},  \label{eq:ren-def-weight}
\end{equation}
and the new transition weights,
\begin{equation}
w(\alpha)_{G_{i,p}\to \sigma'} \sim \frac{k_{G_{i,p}\to\sigma'}}{k_{G_{i,p}}+\alpha} 
\label{eq:new-transition-weights}
\end{equation}
In the non-autocatalytic case,
\begin{equation}
 w(\alpha)_{G_{i,p}\to \sigma'} \sim \begin{cases} w(0)_{G_{i,p}\to \sigma'}  \sim
Z(0)_{{\cal G}_p}^{-1}\ \times\  \max_{\sigma\in V(G_p)} \frac{k_{\sigma\to\sigma'}}{k_{\sigma}}, 
\qquad \alpha\preceq k_{G_{i,p}}   \\   \frac{k_{G_{i,p}\to \sigma'}}{\alpha} \sim (\alpha\tau_{{\cal G}_p})^{-1} 
\ \times\  \max_{\sigma\in V(G_p)} \frac{k_{\sigma\to\sigma'}}{k_{\sigma}}, \qquad \alpha\succeq  
 k_{G_{i,p}}  
\end{cases}.
\end{equation} 
In the autocatalytic case, on the other hand, we only consider  $\alpha\succeq 
\lambda_{{\cal G}_p} \succ k_{G_{i,p}}$; then 
\begin{equation}
 w(\alpha)_{G_{i,p}\to \sigma'} \sim  \frac{k_{G_{i,p}\to \sigma'}}{\alpha}, \qquad \alpha\succeq  
\lambda_{{\cal G}_p}.  \label{eq:w(alpha)-free-autocata-regimes}
\end{equation} 
The {\em free regime} is obtained in the non-autocatalytic case for $\alpha \prec  k_{G_{i,p}}$. 
The {\em autocatalytic regime} is obtained in the autocatalytic case for $\alpha \sim  \lambda_{{\cal G}_p}$. In both cases,  
\BEQ w(\alpha)_{G_{i,p}\to \sigma'} \sim w(\lambda_{{\cal G}_p})_{G_{i,p}\to \sigma'} \sim 
Z_{G_p}^{-1} \ \times\  \max_{\sigma\in V(G_p)} \frac{k_{\sigma\to\sigma'}}{k_{\sigma}} \qquad {\mathit{(free\ and\ autocatalytic\ regimes)}}  \label{eq:walphaGsigma'Z}
\EEQ 
This defines a threshold value, $\alpha_{thr} \sim k_{G_{i,p}}$ (non-autocatalytic case), 
$\alpha_{thr}\sim \lambda_{{\cal G}_p}$ (autocatalytic case); when $\alpha\succ \alpha_{thr}$, defining the {\em degraded regime}, 
the denominator in the expression (\ref{eq:new-transition-weights}) behaves like $\alpha$.  {\em Note 
also that, in all three regimes, the prefactor in $w(\alpha)_{G_{i,p}\to\sigma'}$ is 
$Z(0,\alpha)_{{\cal G}_p}^{-1}$.}  This is to be remembered when considering (\ref{eq:pisigma1}),
(\ref{eq:pisigma2}) below.

\Medskip  Since $k_{G_{i,p}} \prec k_{p,min}$, the next-step cut-off scale will be $<n(i)$, allowing
downward induction on $n$.  Note also that, when $G_{i,p}$ is autocatalytic, $\kappa_{G_{i,p}}
\succ k_{G_{i,p}}$ is the highest reaction scale with reactant $G_{i,p}$, see below (\ref{eq:dominant edges}); therefore, if $\alpha\succeq \kappa_{G_{i,p}}$,  
\BEQ w(\alpha)_{G_{i,p}\to \sigma'} \sim 
\frac{k_{G_{i,p}\to\sigma'}}{\alpha} \prec 1 \label{eq:outgoing-autocata-small-weight}
\EEQ   
Thus we have this essential fact: {\em edges outgoing from an autocatalytic vertex are small}. 

\Medskip\noindent {\bf Final step.} The next renormalization steps $(i=2,3,\ldots)$ are exactly similar.
 Renormalization stops when there are no more non-trivial maximal
dominant SCCs. Each step involves a non-trivial nesting step $\Sigma^{(i)}\to \Sigma^{(i+1)}$, hence the total number of steps, $i_{max}$ is less than the number of species. The step $i$ effective
graph $G^{(i)}=(V(G^{(i)}),E(G^{(i))})$ is obtained from the bare graph $G^{(0)}$ by 
successive merging/rewiring steps. 

\bigskip\noindent {\bf Step $i$ Lyapunov data.} In order to follow the evolution of the system as slower and slower transitions are incorporated, we consider the sequence of cut-off graphs
$G_{cut}(i)$, $i=0,\ldots,i_{max}$, and compute their Lyapunov data (for $i<i_{max}$, these
may be understood as "transient" Lyapunov data, though we do not discuss dynamics here). By definition,  $E(G_{cut}(i))$ is
the set of bare edges (edges of $G^{(0)}$ with their rates) involved in the formation of the
effective graph $G^{(i)}_{cut}$ cut-off just before scale $n(i+1)$, and $V(G_{cut}(i))\subset\Sigma$ the subset of sources and targets of edges in  $E(G_{cut}(i))$.   By construction, $G_{cut}(i_{max})=G^{(0)}$, and $\lambda^*(G_{cut}(i))\le 
\lambda^*(G_{cut}{(i+1)})$, since cutting edges reduces the growth rate.

\Medskip We wish to approximate the Lyapunov weights of  $G_{cut}(i)$; when the latter is not strongly connected, the Lyapunov vector is not uniquely determined, which leads us to the following "initial condition dependent"
construction.   Fix $\sigma_0\in V(G_{cut}(i))$; we let $(G_{cut}(i))_{\sigma_0} \subset 
G_{cut}(i)$ be the subgraph with vertex subset  $V((G_{cut}(i))_{\sigma_0})= \{\sigma\in V(G_{cut}(i))\ |\ \sigma$   accessible from $\sigma_0\}$, where "$\sigma$ accessible from $\sigma_0$" means : $\sigma$ is connected to $\sigma_0$ by some path; if we assume that 
the initial condition is $X_{\sigma}(0) = \delta_{\sigma,\sigma_0}$, only this subgraph
can be reached.  The {\bf $\sigma_0$-SCCs} of $G_{cut}(i)$     are the maximal dominant 
SCCs $G_q, q=1,2,\ldots$ of $G_{cut}(i)$ included in $(G_{cut}(i))_{\sigma_0}$. 

\Medskip {\bf Lyapunov exponent of $(G_{cut}(i))_{\sigma_0}$.}  Each $G_q$ has a threshold rate $\alpha_q$, which
 is (by definition) $0$ if $G_q$ is not autocatalytic, otherwise gives the order of magnitude
 of the Lyapunov exponent $\lambda^*(G_q)$, with logarithm equal to the deficiency scale of $G_q$. The $G_q$ are connected between themselves in
 various ways, but a non-autocatalytic $G_q$ has no outgoing edge in    $(G_{cut}(i))_{\sigma_0}$  (otherwise $G_q$ would not be maximal).   The {\bf threshold scale} 
$\lfloor \log_b \alpha \rfloor$ of 
$(G_{cut}(i))_{\sigma_0}$ (logarithm of the {\bf threshold rate}) is the maximal deficiency scale (if any), $-\infty$ else ($\alpha=0$).  If $\alpha>0$,
\BEQ  \lambda^*((G_{cut}(i))_{\sigma_0}) \sim \alpha \EEQ
The growth rate of the graph $G_{cut}(i))$ started from $\sigma_0$ will
be $\sim \alpha$.  Otherwise (banning resonance cases), $(G_{cut}(i))_{\sigma_0}$ is not
autocatalytic. 

\Medskip {\bf Cores.} Cores are maximal dominant $\sigma_0$-SCCs  maximizing
the set $\{\alpha_q,q=1,2,\ldots\}$. If $\alpha=0$, all $G_q,q=1,2,\ldots$ are cores, and
$(G_{cut}(i))_{\sigma_0}$ is not autocatalytic. In the contrary case ($(G_{cut}(i))_{\sigma_0}$
autocatalytic), we ban the {\bf resonance regime} defined by the case when there exist
$\alpha_{q},\alpha_{q'}$ with $q\not=q'$ such that $\alpha_q\sim\alpha_{q'}\sim \alpha$ are
both maximal. Thus (by reindexing), we may assume that $\alpha\sim \alpha_1\succ \alpha_q, q\not =1$, and 
$G_1$ is the only core. In
the non-autocatalytic case, $\alpha_1=\cdots=\alpha_q=\alpha=0$, so that all $G_q$ are cores.

\Medskip {\bf Hierarchical formulas for Lyapunov vector/weights of $(G_{cut}(i))_{\sigma_0}$.}
We approximate the Lyapunov eigenvector of $(G_{cut}(i))_{\sigma_0}$
by $v_{\sigma}\sim (k_{\sigma} + \alpha_{\sigma_0})^{-1} \pi_{\sigma}$, where
\begin{equation}    \pi_{\sigma} \sim  \prod_{\sigma'\supsetneq \sigma} 
Z^{-1}_{\sigma'} , 
\qquad \sigma\subset G_q\   {\mathrm{core}}   \label{eq:pisigma1}
\end{equation}
\begin{eqnarray}    
&& \pi_{\sigma} \sim\Big(\max_{\gamma_: G_1\to G^{(i)}(\sigma)} w(\alpha_1)_{\gamma} \Big) \times \nonumber\\
 && \qquad
 \times \prod_{\sigma\subsetneq\sigma'\subset \bar{\sigma}} (Z(0,\alpha_1)_{\sigma'} )^{-1} \qquad
{\mathrm{else}}    \label{eq:pisigma2}
\end{eqnarray}
if $\sigma\in V((G_{cut}(i))_{\sigma_0})$, where: 
\begin{itemize}
\item[(i)] in (\ref{eq:pisigma1}),
$"\sigma\subset G_q\   {\mathrm{core}}$"
means: $\exists q, \ \sigma\in V(G_q)$ . Then   the product $\prod_{\sigma'\supsetneq 
\sigma}(\cdots)$ is over the chain of merged vertices containing $\sigma$ (if any);

\item[(ii)] in (\ref{eq:pisigma2}), it is assumed that $\sigma$ is not in a core, but
there exists a core $G_q$ and a  path $\gamma: G_q\to G^{(i)}(\sigma)$  (otherwise 
$\pi_{\sigma}=0$); by construction,
$G_q$ is autocatalytic, so that $q=1$.  Then  $\bar{\sigma}$ is the maximum compound
vertex containing $\sigma$, and $\prod_{\sigma\subsetneq\sigma'\subset \bar{\sigma}} (\cdots)$ is the product
over  the chain of merged vertices  $\sigma'$ contained in  $\bar{\sigma}$. 
\end{itemize}

\Medskip {\em Remark.} When $\sigma$ is not nested in a core, see (\ref{eq:pisigma2}),  paths $\gamma$  cannot, by construction, form an  $\alpha$-dominant cycle.  Though there is no upper  bound over the length of paths $\gamma$, which
may be arbitrarily large if there are cycles, the maximum order of magnitude is attained over the 
finite subset of simple excursions (paths with no self-intersection). Moreover, we can write a generalization of  Dijkstra's algorithm for finding the shortest paths between nodes in a weighted graph, which returns the set of maximal weight $\gamma$'s
in the form of  a DAG {\em (directed acyclic graph)}  rooted in $G_1$ 
(see \S \ref{subsection:DAG}).

\Medskip {\bf Hierarchical formula for adjoint Lyapunov eigenvector.} The estimate for $v^{\dagger}$ is simpler. Reversing edges, let $\bar{\Sigma}_{\sigma_0}\subset \Sigma$
be the subset of bare vertices from which a $\sigma_0$-core is accessible. Then
\begin{equation}
v^{\dagger}_{\sigma} \sim \max_{G_q\, {\mathrm{core}},\, \gamma^{\dagger}:G^{(i)}(\sigma)\to G_q} 
w(\alpha_1)_{\gamma^{\dagger}}  \label{eq:vdaggersigma}
\end{equation}
In (\ref{eq:vdaggersigma}), it is assumed that there exists a core $G_q$ and a reverse path   $\gamma^{\dagger} : G^{(i)}(\sigma) \to G_q$ (otherwise $v^{\dagger}_{\sigma}=0$). 
In particular, $v^{\dagger}_{\sigma}\sim 1$ if $\sigma$ is a core.  The maximum order
of magnitude, as in the case of the hierarchical formula for $\pi$, is attained over the
finite subset of simple excursions, and may be explicited in terms of {\em reverse} DAGs
with edges oriented towards the root $G_1$ (instead of away from $G_1$).

\Medskip The key formulas (\ref{eq:pisigma1})--(\ref{eq:vdaggersigma}) constitute what 
we call the {\bf hierarchical formulas}.


\subsection{Two examples}  \label{subsection:examples}


We illustrate the above estimates on two examples, which are simple enough to be dealt with by
hand.


\subsubsection{Example 1}  \label{subsubsection:example1}


We consider the following reaction network connecting species indexed by 1,2,3:
$1\overset{k_{1\to 2}}{\to} 2, 2\overset{k_{2\to 1}}{\to } 1, 2\overset{k_{2\to 3}}{\to} 3$, and 
a doubling reaction $1\overset{\kappa_1}{\to} 1+1$.  Let $\eps:= \eps_1 = \frac{\kappa_1}{k_{1\to 2}}$. We assume the following three conditions:
\begin{itemize}
\item[(i)] $\eps\ll 1$;
\item[(ii)] $k_{2\to 3} \ll \min(k_{1\to 2},k_{2\to 1})$;
\item[(iii)] $\eps \gg \frac{k_{2\to 3}}{k_{2\to 1}}$.
\end{itemize}
 The model is easily seen to be equivalent when $k_{1\to 2}\gg k_{2\to 1}$ to the cycle
$1\rightleftarrows 2$, $1\overset{\nu_+}{\to} 2+2$, $2\overset{\xi_2}{\to} 2_{ext}$ of 
\S \ref{subsection:cycle-elementary} up to the change of parameters  $(k_{1\to 2}, 
\kappa_1)\leftrightarrow (k_{1\to 2}+2\nu_+,\nu_+)\sim (k_{1\to 2},\nu_+)$ with external
structure relabeling $3\leftrightarrow 2_{ext}$, $k_{2\to 3}\leftrightarrow \xi_2$; resp., 
to the autocatalytically extended Michaelis-Menten model of \S \ref{subsection:multi-scale-method}
when $k_{1\to 2}\ll k_{2\to 1}$, with $(1,2,3)\leftrightarrow (S,ES,P), (k_{1\to 2},k_{2\to 1},
k_{2\to 3},\kappa_1)\leftrightarrow (k_+, k_-,k_2,\nu_+)$.

\Medskip (i) implies that $1$ is not an autocatalytic species.  Let $n_{i\to j} = \lfloor \log_b(k_{i\to j}) \rfloor$. (ii) implies  $\{n^{(0)}_1,n^{(0)}_2\} = 
\{n_{1\to 2},n_{2\to 1}\}$ and $n^{(0)}_3= n_{2\to 3}$. Let $k_{min}=\min(k_{1\to 2},k_{2\to 1})$, 
so that  $n^{(0)}_2=\lfloor \log_b(k_{min})\rfloor$. Then the step 1 cut-off scale is  $n(1) = n^{(0)}_2$;
the unique step 1 maximal dominant SCC is ${\cal G}_1 = \{1,2\}$, and $\bar{\eps}_{{\cal G}_1} \sim \eps$,
$Z(0)_{{\cal G}_1}  = \frac{k_{2\to 3}}{k_{2\to 1}}$. By (iii), $G_1$ is autocatalytic. The step 1
graph $G^{(1)}$ contains two reactions $G_1\overset{k_{G_1\to 3}}{\to}  3$ and an equivalent
doubling reaction (self-edge)   $G_1 \overset{\kappa_{G_1}}{\to} G_1 + G_1$, with  $k_{G_1\to 3} \sim 
k_{min} \, \times\, \frac{k_{2\to 3}}{k_{2\to 1}}$, and $\kappa_{G_1} \sim \lambda_{{\cal G}_1} \sim \eps k_{min}$. 
Since $G_1$ is autocatalytic, $w(\lambda_{{\cal G}_1})_{G_1\to 3} \sim \frac{k_{G_1\to 3}}{\lambda_{G_1}} \sim 
\eps^{-1} \frac{k_{2\to 3}}{k_{2\to 1}} \ll 1$.

\Medskip This example is simple enough to be solved directly by formal diagonalization. Namely, 
$A = \left[\begin{array}{ccc} -k_{1\to 2}+\kappa_1 & k_{2\to 1} & 0 \\ k_{1\to 2} & -k_{2\to 1} -k_{2\to 3} 
 & 0 \\ 0 & k_{2\to 3} & 0 \end{array}\right]$. We let $A' = (a_{i,j})_{1\le i,j\le 2}$; the non-zero eigenvalues  of $A$ may be found by solving $\det (A'-\lambda^* \Id)=0$. Here (using (i)-(iii)) $-\Tr(A') \sim \max(k_{1\to 2},
 k_{2\to 1})$, and  $\kappa_1 k_{2\to 3} \ll  k_{1\to 2} k_{2\to 3}  \ll \kappa_1 k_{2\to 1}$, whence
 $-\det(A') = \kappa_1(k_{2\to 1} + k_{2\to 3}) - k_{1\to 2} k_{2\to 3} \sim \kappa_1 k_{2\to 1} \sim 
 \eps k_{1\to 2} k_{2\to 1} \ll (\Tr(A'))^2$. Hence: $\lambda^* =\half(\Tr(A') + \sqrt{(\Tr(A'))^2 -4
 \det(A')}) \sim \frac{|\det(A')|}{|\Tr(A')|} \sim \eps k_{min}$. As expected from our results, $\lambda^*\sim \lambda_{{\cal G}_1}$.  Solving for $(A-\lambda^*)v=0$, we find
 $k_{1\to 2} v_1\sim k_{2\to 1} v_2, k_{2\to 3} v_2 \sim \lambda^* v_3$, whence $v\sim \left[\begin{array}{c}
 k_{2\to 1}/k_{1\to 2} \\ 1 \\ k_{2\to 3}/\eps k_{min} \end{array}\right]$ up to normalization. 
 
\Medskip Consider now the equivalent discrete-time eigenvalue problem;  $\pi^* = 
\left[\begin{array}{c} (k_{1\to 2}-\kappa_1 + \lambda^*)v_1 \\ (k_{2\to 1} + k_{2\to 3}+ \lambda^*)v_2 \\
\lambda^* v_3 \end{array}\right] \sim  \left[\begin{array}{c} k_{1\to 2}v_1 \\ k_{2\to 1}v_2 \\
\lambda^* v_3 \end{array}\right] \sim \left[\begin{array}{c} k_{2\to 1} \\ k_{2\to 1} \\ k_{2\to 3}
\end{array}\right]$ 
 is an (unnormalized) solution of the eigenvalue problem $(\, ^t {\cal W}(\lambda^*) -\Id)\pi^*=0$,
  with $^t{\cal W}(\lambda^*) = \left[\begin{array}{ccc} 0 & \frac{k_{2\to 1}}{k_{2\to 1}+k_{2\to 3} + 
\lambda^*} & 0 \\ \frac{k_{1\to 2}}{k_{1\to 2} - \kappa_1 + \lambda^*} & 0 & 0 \\
 0 & \frac{k_{2\to 3}}{k_{2\to 1} + k_{2\to 3} + \lambda^*} & 0 \end{array}\right]$.   This is in agreement (up to normalization) with our hierarchical formulas (\ref{eq:pisigma1}), 
 (\ref{eq:pisigma2}), for $\pi$, namely,
\BEQ \pi_1,\pi_2 \sim Z_{G_1}^{-1} \sim \bar{\eps}_{{\cal G}_1}^{-1} \sim \eps^{-1}, \qquad
\pi_3  w(\lambda_{{\cal G}_1})_{G_1\to 3} \sim \eps^{-1} \frac{k_{2\to 3}}{k_{2\to 1}}
\EEQ
Considering the adjoint problem, the solution of $({\cal W}(\lambda^*)-\Id)v^{\dagger,*}=0$
satisfies as expected $v^{\dagger,*}_{1,2}\sim 1, v^{\dagger,*}_3=0$.  
 
\Medskip {\em Green kernel.}  \qquad Let $G^{\lambda}(\sigma,\sigma') := \sum_{\gamma:\sigma\to\sigma'} 
w(\lambda)_{\sigma\to\sigma'} \in \bar{\R}_+$, \\ $G^{\lambda}_n(\sigma,\sigma') := \sum_{\ell \le n}
\sum_{\gamma:\sigma\to\sigma'\ |\ \ell(\gamma)=\ell} w(\lambda)_{\sigma\to\sigma'}$ and 
$G^*_n(\cdot,\cdot):=G^{\lambda^*}_n(\cdot,\cdot)$.  Resumming over iterated loops $1\to 2\to 1$, we get the geometric series  $G^{\lambda}(2,3) = \frac{1}{1-r(\lambda)} \, \times w(\lambda)_{2\to 3}$, with $r(\lambda)= w(\lambda)_{2\to 1} w(\lambda)_{1\to 2}$. 
Expanding the equation $\det(A'-\lambda^*\Id)=0$, one sees that $\lambda^*$ is defined precisely by $r=1$,
so that the series diverges. The truncated Green kernel is thus easily seen to scale like $N$, 
\BEQ \frac{G^*_{2N+1}(2,1)}{N+1} = w(\lambda^*)_{2\to 1} \sim 1,   \qquad
\frac{G^*_{2N}(2,2)}{N+1} = 1, \qquad \frac{G^*_{2N+1}(2,3) }{N+1} = w(\lambda^*)_{2\to 3} \sim k_{2\to 3}/ k_{2\to 1}
\EEQ
which is proportional to $\pi^*$ to leading order.

\begin{figure}[!t]
\centering
\includegraphics[scale = 0.8]{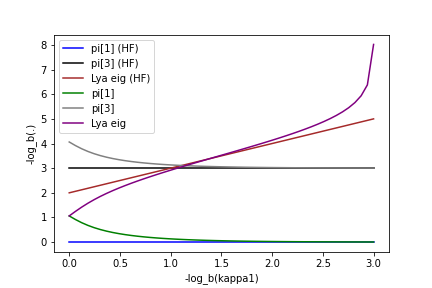}
\caption{Numerical checks for  Example 1. Legend: (HF) = hierarchical formulas,
see text.}
\label{fig:num_check_ex1_var}
\end{figure}

\Medskip Hierarchical formulas are reported on Fig. \ref{fig:num_check_ex1_var} in the form of 
log-log plots, and compared with numerical simulations. Parameters are $k_{1\to 2}=1, k_{2\to 1} = b^{-2}, k_{2\to 3} = b^{-5}$, while
$\kappa_1$ is allowed to vary in the {\em domain} $(b^{-3},1)$. When $\kappa_1\sim 
\frac{k_{2\to 3}}{k_{2\to 1}} \, \times\, k_{1\to 2} = b^{-3}$, the Lyapunov eigenvalue 
(in purple) becomes negative $(-\log_b \lambda^*\to +\infty)$, as predicted above; then  the
Lyapunov eigenvector is along the vector
$(0 \ 0 \ 1)^t$. In the autocatalytic range, hierarchical formulas yields $-\log_b \lambda$
(brown), and $-\log_b \pi_1$ (blue), resp. $-\log_b \pi_3$ (black), as compared to the
numerical formulas in purple, resp. green, resp. grey. Numerical values and hierarchical formulas for  $-\log_b\pi^*_2$, not reported, are both very close to 0. As a general rule,
hierarchical formulas are a good approximation {\em inside} a domain, far from its walls. The
word {\em domain} refers to \S \ref{subsection:aim}, the reader may easily compute the
coefficients $c_e^{T,i}, c_e, C_{\sigma,e}$ where $(\T,i)$ points here to the domain 
$\{k_{1\to 2}\succ k_{2\to 1} \succ k_{2\to 3},  k_{1\to 2} \succ \kappa_1 \succ \frac{k_{2\to 3}}{k_{2\to 1}} \, \times\, k_{1\to 2} \}$.

\subsubsection{Example 2. Two coupled cycles}  \label{subsubsection:example2}


 Next, we consider
the following multi-scale reaction network with 4 species indexed by $1,2,\bar{1},\bar{2}$,
\begin{eqnarray*}
(k_{max})\qquad  1\to 2 \\
(\bar{k}_{max}) \qquad \bar{1}\to \bar{2} \\
(\bar{\nu}_{+,1}) \qquad \bar{1} \to \bar{1} +\bar{1} \\ 
--------- \\
(k_{min}) \qquad 2\to 1 \\
--------- \\
(\bar{k}_{min})\qquad \bar{2}\to \bar{1} \\
(\nu_{+,1}) \qquad 1 \to 1+1 \\
(k_+) \qquad 2\to 2+2, 1\to \bar{1} \\
(k_-) \qquad \bar{1}\to 1, \bar{2}\to \bar{2} + \bar{2} 
\end{eqnarray*}

\noindent with $k_{max}\succ \bar{k}_{max}\succ \bar{\nu}_{+,1}\succ k_{min}\succ \bar{k}_{min}\succ 
\nu_{+,1}\succ k_+\succ k_-$, under the following assumptions:
\begin{itemize}
\item[(i)] $\frac{k_+}{k_{min}} \succeq \frac{\nu_{+,1}}{k_{max}}$;
\item[(ii)] $\frac{\bar{\nu}_{+,1}}{\bar{k}_{max}} \succ \frac{k_+}{\bar{k}_{min}}$.
\end{itemize}

We illustrate our computations by choosing $k_{max} = b^0, \bar{k}_{max} = b^{-2}, \bar{\nu}_{+,1}= b^{-3}, 
k_{min} = b^{-5}$, $\bar{k}_{min} = b^{-6}, \nu_{+,1} = b^{-7}, k_+ = b^{-12}$
and $k_-= b^{-16}$. We get the following multi-scale representation, see \S \ref{subsection:cycle-elementary} (self-edges have been removed for simplicity),

\begin{center}
\begin{tikzpicture}

\draw[dashed](-0.55,-2.4) rectangle(2.55,0.5);
\draw[dashed](4-0.55,-3.4) rectangle(4+2.55,-0.2);
\draw(-1,-1) node {${\cal G}_1$}; 
\draw(3,-3) node {${\cal G}_2$};

 \draw(2,1) node {$1$};
\draw[->, ultra thick](2,0)--(0,0); \draw(1,0.3) node {$k_{max}$};
\draw(0,1) node {$2$};

\draw(4,1) node {$\bar{1}$};
\draw[->, ultra thick](4,-1)--(6,-1); \draw(5,-1+0.3) node {$\bar{k}_{max}$};
\draw(6,1) node {$\bar{2}$};

\draw(8,0) node {$n_1^{(0)}$};
\draw[dashed](-1,-0.5)--(8,-0.5); 

\draw[dashed](-1,-2.5)--(8,-2.5);  \draw(8,-2) node {$n_3^{(0)}$};
\draw(8,-3) node {$n_4^{(0)}$};

 \draw[->, ultra thick](0,-2)--(2,-2);
\draw(1,-1-0.7) node {$k_{min}$}; 

 \draw[->, ultra thick](6,-3)--(4,-3);
\draw(5,-2-0.7) node {$\bar{k}_{min}$};
\draw(8,-1) node {$n_2^{(0)}$};
\draw[dashed](-1,-1.5)--(8,-1.5);

\draw[->](2,-4)--(4,-4); \draw(3,-2-1.7) node {$k_+$};
\draw(8,-4) node {$n_5^{(0)}$};
\draw[dashed](-1,-3.5)--(8,-3.5);

\draw[<-](2,-5)--(4,-5); \draw(3,-2-2.7) node {$k_-$};
\draw(8,-5) node {$n_6^{(0)}$};
\draw[dashed](-1,-4.5)--(8,-4.5);

\draw[dotted](0,0)--(0,-2); 
\draw[dotted](2,0)--(2,-5); 
\draw[dotted](4,-1)--(4,-5); 
\draw[dotted](6,-1)--(6,-3);

\end{tikzpicture}
\end{center}

\medskip Following our analysis, we have $\eps_1 \sim \frac{\nu_{+,1}}{k_{max}}= b^{-7}, \eps_2 \sim
\frac{k_+}{k_{min}} = b^{-7}, \eps_{\bar{1}} \sim \frac{\bar{\nu}_{+,1}}{\bar{k}_{max}} = 
b^{-1}, \eps_{\bar{2}} \sim \frac{k_-}{\bar{k}_{min}} = b^{-10}$. Assumption (i) is equivalent to 
$\eps_2\succeq \eps_1$; assumption (ii), which amounts to $\lambda_{{\cal G}_2}\succ 
\lambda_{{\cal G}_1}$, see below, implies in particular  $\eps_{\bar{1}}\succeq \eps_{\bar{2}}$. The original scales (from top to bottom) are 
$n^{(0)}_{1,\cdots,8} =n_{max}$,  $\bar{n}_{max}$, $\bar{n}_{+,1}$, $n_{min}$, $\bar{n}_{min}$, $n_{+,1}$, 
$n_+$, $n_-$. 

\Medskip  (1) At cut-off scale $n(1) = n_{min}$, the cycle
$1 \underset{k_{min}}{\overset{k_{max}}{\rightleftarrows}} 2$ defines a dominant SCC 
$G_1 = \{1 \rightleftarrows 2\}$ with renormalized weight $Z(0)_{{\cal G}_1} \sim \frac{k_+}{k_{max}} \sim b^{-12}$ and
deficiency weight $\bar{\eps}_{G_1}  \sim \max(\eps_1,\eps_2)$. Assumption (i) means $\eps_2\succeq 
\eps_1$, therefore $\bar{\eps}_{G_1} \sim \eps_2 \sim b^{-7} \succ Z(0)_{{\cal G}_1}$. Thus $G_1$ is autocatalytic, and $\lambda_{{\cal G}_1} \sim \bar{\eps}_{G_1} k_{min} \sim k_+ \sim 
b^{-12}$, $Z_{G_1}\equiv Z^{(1)}_{G_1}\sim \bar{\eps}_{G_1}$.   The new effective graph has vertices $\{G_1,\bar{1},\bar{2}\}$ and 
edges $k_{\bar{1}\to G_1}
\sim k_{\bar{1}\to 1} \sim k_-\sim b^{-16}$, $k_{G_1\to \bar{1}} \sim \frac{k_+}{k_{max}} k_{min} \sim b^{-17}$, both $\prec \bar{k}_{min}$. 
Choosing $\sigma_0=1,2$, $G_{cut,\sigma_0}(1)=\{1\rightleftarrows 2\}$ has   Lyapunov eigenvector  components 
$v^{(1)}_1 \sim k_{max}^{-1} Z_{G_1}^{-1} \sim k_{min}/(k_+ k_{max}) \sim b^{7}, 
v^{(1)}_2\sim k_{min}^{-1} Z_{G_1}^{-1}\sim k_+^{-1}\sim b^{12}$;  dividing by $b^{12}$ for normalization, $v^{(1)} \sim  \left(\begin{array}{c} b^{-5} \\ 1 \\ 
0 \\ 0 \end{array}\right)$, while $v^{\dagger,(1)}\sim \left(\begin{array}{c} 1 \\ 1 \\ 0 \\
0 \end{array}\right)$. 

\Medskip  (2) At cut-off scale $n(2) = 
\bar{n}_{min}$, the cycle
$\bar{1} \underset{\bar{k}_{min}}{\overset{\bar{k}_{max}}{\rightleftarrows}} \bar{2}$ defines a dominant SCC 
$G_2 = \{\bar{1},\bar{2}\}$ with renormalized weight $Z(0)_{{\cal G}_2} \sim \frac{k_{\bar{1}\to G_1}}{\bar{k}_{max}} \sim \frac{k_-}{\bar{k}_{max}} \sim b^{-14}$ and
deficiency weight $\bar{\eps}_{G_2}  \sim \max(\eps_{\bar{1}},\eps_{\bar{2}})\sim  \eps_{\bar{1}}$ 
as a consequence of Assumption (ii), whence $\bar{\eps}_{G_2}\sim  b^{-1} \succ Z(0)_{{\cal G}_2}$. Thus $G_2$ is autocatalytic, and $\lambda_{{\cal G}_2} \sim \bar{\eps}_{G_2} \bar{k}_{min} \sim \frac{\bar{\nu}_{+,1}}{\bar{k}_{max}} \bar{k}_{min} \sim 
b^{-7}$, $Z_{G_2}\sim \bar{\eps}_{G_2}$.  The new effective graph has vertices $\{G_1,G_2\}$
and edges $k_{G_1\to G_2} \sim k_{G_1\to \bar{1}}$, $k_{G_2\to G_1} \sim \frac{k_{\bar{1}\to G_1}}{k_{\bar{1}}} \bar{k}_{min} \sim \frac{k_-}{\bar{k}_{max}} \bar{k}_{min} \sim b^{-20}$.
Assumption (ii) is equivalent to  $\lambda_{{\cal G}_2}\succ \lambda_{{\cal G}_1}$. Hence the new root is $G_2$, and  now 
$Z_{G_1}\equiv Z^{(2)}_{G_1} \sim Z(\lambda_{{\cal G}_2})_{{\cal G}_1} \sim \frac{\lambda_{{\cal G}_2}}{k_{min}} \sim b^{-2}$, $w(\lambda_{{\cal G}_2})_{G_2\to G_1} \sim \frac{k_{G_2\to G_1}}{\lambda_{{\cal G}_2}}\sim 
\frac{k_-}{\bar{\nu}_{+,1}} \sim b^{-13}$.  
Then hierarchical formulas imply (using standard basis along $1,2,\bar{1},\bar{2}$) 
\BEQ \pi \propto \left(\begin{array}{c} 
 w(\lambda_{{\cal G}_2})_{G_2\to G_1} (Z^{(2)}_{G_1})^{-1} \\ 
 w(\lambda_{{\cal G}_2})_{G_2\to G_1} (Z^{(2)}_{G_1})^{-1} \\
 Z_{G_2}^{-1} \\
 Z_{G_2}^{-1} \end{array}\right)  \sim \frac{\bar{k}_{max}}{\bar{\nu}_{+,1}} 
 \left(\begin{array}{c} (k_- k_{min})/(\bar{\nu}_{+,1} \bar{k}_{min}) \\
    (k_- k_{min})/(\bar{\nu}_{+,1} \bar{k}_{min}) \\ 1 \\ 1 \end{array}\right),
 \label{eq:ex2-pi}
    \EEQ
     from which 
$v \sim 
\left(\begin{array}{c} k_{max}^{-1} \pi^*_1 \\ 
k_{min}^{-1}\pi^*_2 \\
 \bar{k}_{max}^{-1} \pi^*_{\bar{1}}\\
 \bar{k}_{min}^{-1}  \pi^*_{\bar{2}} \end{array}\right) \sim \lambda_{{\cal G}_2}^{-1}  
 \left(\begin{array}{c} b^{-18} \\ b^{-13} \\ b^{-4} \\ 1 \end{array}\right)$.  The symbol $\propto$ means: equivalent $(\sim)$ up to normalization. As follows from the Remark below, we see that
 $\pi\sim   \left(\begin{array}{c} (k_- k_{min})/(\bar{\nu}_{+,1} \bar{k}_{min}) \\
    (k_- k_{min})/(\bar{\nu}_{+,1} \bar{k}_{min}) \\ 1 \\ 1 \end{array}\right) \sim \left(\begin{array}{c} b^{-12}\\
b^{-12} \\ 1 \\ 1 \end{array}\right)$ is correctly
    normalized as a probability measure $(\sum_{\sigma} \pi_{\sigma}\sim 1$).  Looking
 at the adjoint problem instead, and using $w(\lambda_{{\cal G}_2})_{G_1\to G_2} \sim 
 \frac{k_{G_1\to \bar{1}}}{\lambda_{{\cal G}_2}} \sim \frac{k_+ k_{min}/k_{max}}{\bar{\nu}_{+,1} 
 \bar{k}_{min}/\bar{k}_{max}} \sim b^{-10}$, we get
\BEQ v^{\dagger} \propto  \left(\begin{array}{c} w(\lambda_{{\cal G}_2})_{G_1\to G_2} \\ w(\lambda_{{\cal G}_2})_{G_1\to G_2} \\ 1 \\ 1 \end{array}\right) \sim \left(\begin{array}{c} b^{-10}\\
b^{-10} \\ 1 \\ 1 \end{array}\right)
\EEQ
which is already correctly normalized in such a way that $\langle v^{\dagger},\pi\rangle\sim 1$. 

\Bigskip Lyapunov data are not easy to compute directly from the characteristic equation
$\det (A-\lambda^*\Id)=0$, since the latter does not factorize in this case, and   $A$ is now  a $4\times 4$ matrix. On the other hand, the Green kernel approach, although a bit lengthy, works. 
We sketch the computations. Let $r(\lambda):=w(\lambda)_{2\to 1}w(\lambda)_{1\to 2} = 
\frac{k_{min}}{k_{min}+\lambda} \, \frac{k_{max}}{k_{max}-\nu_{+,1} + k_+ + \lambda}$, 
     $\bar{r}(\lambda):=w(\lambda)_{\bar{2}\to \bar{1}}w(\lambda)_{\bar{1}\to \bar{2}} 
     = \frac{\bar{k}_{min}}{\bar{k}_{min}+\lambda} \, \frac{\bar{k}_{max}}{\bar{k}_{max}-\bar{\nu}_{+,1} + k_- + \lambda}$. The
threshold equation is $R(\lambda):=\frac{1}{1-r(\lambda)} w(\lambda)_{1\to \bar{1}} \frac{1}{1-\bar{r}(\lambda)}
w(\lambda)_{\bar{1}\to 1} \equiv 1$. Then $R(\lambda^*)=1$ implies  
\BEQ (\frac{\lambda^*}{k_{min}}-\frac{\nu_{+,1}}{k_{max}})(\frac{\lambda^*}{\bar{k}_{min}}-\frac{\bar{\nu}_{+,1}}{\bar{k}_{max}})  \sim \frac{k_+ k_-}{k_{max}\bar{k}_{max}}. 
\label{eq:thr_ex2}
\EEQ
 Solving instead 
$(\frac{\lambda}{k_{min}}-\frac{\nu_{+,1}}{k_{max}})(\frac{\lambda}{\bar{k}_{min}}-
\frac{\bar{\nu}_{+,1}}{\bar{k}_{max}}) =0$ yields $\lambda\sim \lambda_{{\cal G}_1}$
 or $\lambda_{{\cal G}_2}$. 
The product is $<0$ if $\lambda\sim\lambda_{{\cal G}_1}$, so $\lambda^*\sim \lambda_{{\cal G}_2}$; this is consistent with (\ref{eq:thr_ex2}) since, for $\lambda^*\sim\lambda_{{\cal G}_2}$,
the l.-h.s. is at most $\sim \frac{\lambda_{{\cal G}_2}}{k_{min}} \, \times\, \frac{\lambda_{{\cal G}_2}}{\bar{k}_{min}} \succ \frac{\lambda_{{\cal G}_1}}{k_{min}} \, \times\,  \frac{\lambda_{{\cal G}_2}}{\bar{k}_{min}} \sim \frac{k_+}{k_{min}} \, \times\, 
\frac{\bar{\nu}_{+,1}}{\bar{k}_{max}} \succ \frac{k_+ k_-}{k_{max}\bar{k}_{max}}$. 
Then we get:   $1-r(\lambda^*) \sim 
\frac{\lambda^*}{k_{min}} \sim \frac{\bar{\nu}_{+,1}}{\bar{k}_{max}} \frac{\bar{k}_{min}}{k_{min}}$, and (plugging in (\ref{eq:thr_ex2})) $1-\bar{r}(\lambda^*)\sim  \frac{k_+ k_-}{k_{max}\bar{k}_{max}} (\frac{\lambda^*}{k_{min}})^{-1} \sim 
\frac{k_+ k_- k_{min}}{k_{max} \bar{\nu}_{+,1} \bar{k}_{min}}$ (with our choice of scales,
$1-\bar{r}(\lambda^*) \sim b^{-24}$ and $1-r(\lambda^*)\sim b^{-2}$). 
The truncated geometric series $\sum_{k=0}^n (\bar{r}(\lambda^*))^k $ is $\sim  (1-\bar{r}(\lambda^*))^{-1}$ for $n\gg (1-\bar{r}(\lambda^*))^{-1}$. 
Similarly, $\sum_{k=0}^n (r(\lambda^*))^k \sim (1-r(\lambda^*))^{-1} $ for $n\gg (1-r(\lambda^*))^{-1}$. Thus   
\BEA  
G^*_N(\bar{1},1) \approx G^*_N(\bar{1},2) &\approx& w(\lambda^*)_{\bar{1}\to 1} \, \times\, (1-r(\lambda^*))^{-1}
\, \times\,  G^*_N(\bar{1},\bar{1})  \nonumber\\
&\sim & \frac{k_-}{\bar{k}_{max}} \,\times \, (\frac{\bar{\nu}_{+,1}}{\bar{k}_{max}} \frac{\bar{k}_{min}}{k_{min}})^{-1} \, \times\,  G^*_N(\bar{1},\bar{1}) 
\approx  \frac{k_- k_{min}}{\bar{\nu}_{+,1} \bar{k}_{min}} G^*_N(\bar{1},\bar{1}). \nonumber\\
\EEA
The vector $(G^*_N(\bar{1},\sigma))_{\sigma=1,2,\bar{1},\bar{2}}$ is again proportional to $\pi$ to leading order.  

\Medskip {\bf Remark.} Note that $\bar{1},\bar{2}$ are always {\bf $\pi$-dominant} under the above hypotheses, in the sense that $\pi_{\bar{1}}\sim 
\pi_{\bar{2}}\succ \pi_1,\pi_2$. Namely, (ii) implies: $\frac{\pi_1}{\pi_{\bar{1}}}
\sim \frac{k_- k_{min}}{\bar{\nu}_{+,1} \bar{k}_{min}} \prec \frac{k_-}{k_+} \, \times \, 
\frac{k_{min}}{k_{max}}\prec 1$. The variant discussed in the next paragraph allows for
more possibilities.


\subsubsection{Variant of Example 2, comparison with numerics} \label{subsubsection:example2-var}

 
We consider here the "reverse wiring case" where edges connect $2$ and $\bar{2}$, instead of $1$ and $\bar{1}$, i.e. $k_+=k_{2\to \bar{2}}, k_- = k_{\bar{2}\to 2}$, so that 
$Z(0)_{{\cal G}_1} \sim \frac{k_+}{k_{min}}, Z(0)_{{\cal G}_2}\sim \frac{k_-}{\bar{k}_{min}}$. Leave $\bar{\eps}_{G_1},\bar{\eps}_{G_2}$ unspecified, but keep the hypotheses $\bar{\eps}_{G_1}\succ 
Z(0)_{{\cal G}_1}, \bar{\eps}_{G_2}\succ Z(0)_{{\cal G}_2}$, so that the two subgraphs are
autocatalytic, and $\lambda_{{\cal G}_1}\sim \bar{\eps}_{G_1} k_{min}, \lambda_{{\cal G}_2}\sim \bar{\eps}_{G_2} \bar{k}_{min}$. Finally, assume $\lambda_{{\cal G}_2}\succ \lambda_{{\cal G}_1}$ as above. The above three hypotheses are equivalent to $\lambda_{{\cal G}_2}\succ\lambda_{{\cal G}_1}\succ k_+$.  Then we get successively: $k_{G_1\to G_2} \sim k_{G_1\to \bar{2}}\sim k_+$; 
$k_{G_2\to G_1} = k_{G_2\to 2} \sim k_-$;   $Z_{G_1} \sim \frac{\lambda_{{\cal G}_2}}{k_{min}}, 
Z_{G_2} \sim \frac{\lambda_{{\cal G}_2}}{\bar{k}_{min}}$; and $w(\lambda_{{\cal G}_2})_{G_2\to G_1} \sim \frac{k_{G_2\to G_1}}{\lambda_{{\cal G}_2}} \sim \frac{k_-}{\lambda_{{\cal G}_2}}$,
  $w(\lambda_{{\cal G}_2})_{G_1\to G_2} \sim \frac{k_{G_1\to G_2}}{\lambda_{{\cal G}_2}} 
  \sim \frac{k_+}{\lambda_{{\cal G}_2}}$.
Thus (using the general Ansatz of (\ref{eq:ex2-pi}) and plugging in the corrected
values) 
\BEQ \frac{\pi_2}{\pi_{\bar{2}}} \sim \frac{(k_-/ \lambda_{{\cal G}_2})(k_{min}/\lambda_{{\cal G}_2})}{\bar{k}_{min}/\lambda_{{\cal G}_2}}  \sim \frac{k_-}{\lambda_{{\cal G}_2}} \times \frac{k_{min}}{\bar{k}_{min}}  \label{eq:kkklambda}
\EEQ
 is not necessarily $\prec 1$ here. 

\begin{figure}[!t]
\centering
\includegraphics[scale=0.8]{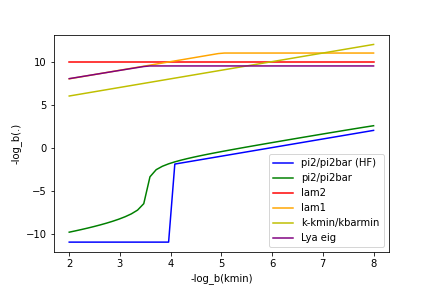}
\caption{Numerical checks for variant of Example 2. Legend: (HF) = hierarchical formulas,
see text.}
\label{fig:num_check_ex2_var}
\end{figure}

\Medskip We report on Fig. \ref{fig:num_check_ex2_var} log-log plots of the
above computations, with parameters $b=5$ (results hardly depend on $b$ at all, as long as 
$b$ is not too close to 1);  $k_{max} = 1, \bar{k}_{max} = b^{-2}, \bar{k}_{min} = 
b^{-9}, k_+ = b^{-12}, k_- = b^{-13}$; $\nu_{+,1} = b^{-6}$, $\nu_{+,2} = b^{-11}$, 
$\bar{\nu}_{+,1} = b^{-3}$, $\bar{\nu}_{+,2} = b^{-10}$; and $k_{min}$ is left free
to vary from $b^{-8}$ to $b^{-2}$. For the computations, we used the formulas $\bar{\eps}_{G_1} = \max(\frac{\nu_{+,1}}{k_{max}}, \frac{\nu_{+,2}}{k_{min}})$, $\bar{\eps}_{G_2} = \max(\frac{\bar{\nu}_{+,1}}{\bar{k}_{max}}, \frac{\bar{\nu}_{+,2}}{\bar{k}_{min}}), \lambda_{{\cal G}_1}= \bar{\eps}_{G_1} k_{min}, \lambda_{{\cal G}_2} = \bar{\eps}_{G_2} \bar{k}_{min}$ found
above.  One has $\bar{\eps}_{G_1}\succ k_+/k_{min}$ in the whole
range, and $\bar{\eps}_{G_2}\sim \frac{\bar{\nu}_{+,1}}{\bar{k}_{max}} \sim b^{-3} \succ k_-/\bar{k}_{min}$, therefore $G_1$ and $G_2$ are autocatalytic. Then $\lambda_{{\cal G}_2}\succ 
\lambda_{{\cal G}_1}$ only if $k_{min}\prec b^{-4}$; the orange curve $-\log_b(\lambda_{{\cal G}_1}(k_{min}))$ and the red curve $-\log_b(\lambda_{{\cal G}_2}(k_{min}))$, actually cross
when $-\log_b(k_{min})=4$. The purple curve, representing $-\log_b(\lambda^*)$, illustrates
the fact that $\lambda^* \sim \max(\lambda_{{\cal G}_1},\lambda_{{\cal G}_2})$.  Then 
the yellow curve represents the log-ratio $-\log_b(\frac{k_- k_{min}}{\bar{k}_{min}})$; it
crosses the red curve $(-\log_b (\lambda_{{\cal G}_2}))$ at $-\log_b k_{min}=6$. This signals
(see (\ref{eq:kkklambda})) a transition from the case when $\bar{1},\bar{2}$ are $\pi$-dominant 
to the case when $1,2$ are $\pi$-dominant, which is indeed apparent on the green (numerical computations) or blue  curve (hierarchical formulas) reporting
$-\log_b(\pi_2/\pi_{\bar{2}})$ . 

\Medskip Around $k_{min}\sim b^{-4}$, both green and blue curves undergo a sharp transition.
When $k_{min}\succ b^{-4}$,  $\lambda_{{\cal G}_1}\succ \lambda_{{\cal G}_2}$,  and 
the formulas for $\pi$ change,
\BEQ \pi_{1,2} \sim Z_{G_1}^{-1}, \ \pi_{\bar{1},\bar{2}} \sim w(\lambda_{{\cal G}_1})_{G_1\to  G_2} \, \times\,  Z_{G_2}^{-1} \sim Z_{G_1}^{-1} \frac{k_+}{k_{max}} \, \times\,  Z_{G_2}^{-1}
\EEQ
with $Z_{G_i} \sim \bar{\eps}_{G_i}$, $i=1,2$, therefore $\pi_2/\pi_{\bar{2}} \sim 
 \frac{k_{max}}{k_+} \bar{\eps}_{G_2}$, which is the value used for the left part of 
 the blue curve. The green curve transitions smoothly between the two formulas, whereas 
 the blue curve is discontinuous at the transition.

\section{Proofs}  \label{section:proofs}


The heuristics developed in \S \ref{subsection:heuristics} was based on a perturbative
expansion of the resolvent formula in terms of paths. The renormalization process relies on
partial resummations of the sum over paths. Typically, the resummation involves  "pseudo-inverses"
$({\cal W}-\Id)^{-1}$ or $(\, ^t{\cal W}-\Id)^{-1}$, where ${\cal W}$, resp. $^t {\cal W}$ is a Markov matrix, resp. its adjoint, with
kernel generated by the constant vector ${\bf 1}$, resp. stationary measure $\pi$. Pseudo-inverses are defined on the
orthogonal of the (co)kernel. The procedure, which is perhaps not totally standard because ${\cal W}$ is 
not symmetric, involves both the kernel and the cokernel, and can be found in Appendix 
\ref{section:app3-pert}.

\Medskip We first (\S \ref{subsection:cut-off-setting}, \S \ref{subsection:cut-off-Lya-estimates}) discuss the theory for cut-off graphs; see Section \ref{section:notations}. 
Namely, we start from a graph $G=(V,E)$, split $V$ into   ${\cal V}^{int} \uplus {\cal V}^{ext}$ and
consider the cut-off graphs ${\cal G}^{int}, {\cal G}$ obtained by stopping the Markov chain
upon hitting ${\cal V}^{ext}$, and estimate the Lyapunov data of ${\cal G}^{int}$, and 
exit probabilities of ${\cal G}$.

\subsection{Cut-off graphs: the setting} \label{subsection:cut-off-setting}

\noindent{\bf General setting.}  We consider in this subsection (see Section \ref{section:notations} for definitions and notations):

\begin{itemize}
\item[\textbullet] a (non-empty, connected) graph $G=(V,E)$, with adjoint Markov/defective Markov  generators $\tilde{A}$, resp. $A$; and transition matrix $\tilde{\cal W}$, resp. $\cal W$;
\item[\textbullet] a splitting of $V$ into ${\cal V}^{int} \uplus {\cal V}^{ext}$ ({\em internal}, resp. 
{\em external} vertices), with ${\cal V}^{int},{\cal V}^{ext}\not=\emptyset$, defining an {\em internal subgraph} ${\cal G}^{int} = ({\cal V}^{int},{\cal E}^{int})$, with
${\cal E}^{int} = \{(\sigma,\sigma')\in E\ |\ \sigma,\sigma'\in {\cal V}^{int}\}$ ({\em internal} edges), and degradation rates $\beta_{\sigma}^{int}=\beta_{\sigma} + \sum_{\sigma'\in {\cal V}^{ext}}
k_{\sigma\to\sigma'}$;
\item[\textbullet]   
$\tilde{A}^{int}$, resp. $A^{int}$,  $\tilde{\cal W}^{int}$, 
resp. ${\cal W}^{int}$ 
the (non-defective, resp. defective) adjoint Markov generators/transition matrices of ${\cal G}^{int}$; also, $\tilde{\pi}^{int}$,
the probability measure associated
to the Markov generator ${\cal W}^{int}$;
\item[\textbullet] finally, the extended cut-off graph ${\cal G}$, obtained by adding to the internal edges in ${\cal G}^{int}$ the
 set ${\cal E}^{out}$ of outgoing edges $k_{\sigma\to\sigma'},\sigma\in \Sigma^{int},\sigma'\in 
 \Sigma^{ext}$.
\end{itemize}

\Medskip {\em We assume that the cut-off graph ${\cal G}$ satisfies properties (i), (ii) introduced
in Section \ref{section:notations}: 
\\   (i) ${\cal G}^{int}$ is connected by  dominant edges in ${\cal E}^{int}$;  \\ 
(ii) every outgoing edge $(\sigma,\sigma'_{ext})\in {\cal E}^{out}$ is dominated
by some internal edge $(\sigma,\sigma')\in {\cal E}^{int}$. }

\Medskip We denote by $\lambda^*$ the Lyapunov exponent of $A^{int}$, and by $v^*$, resp. $v^{\dagger}$, 
associated right, resp. left eigenvectors, which are unique up to normalization.

\Medskip Fix $\alpha\in\R$. The matrix $A(\alpha)$ has a block decomposition along the splitting,  
\BEQ A(\alpha) = \left(\begin{array}{cc} A^{int}-\alpha\Id &  A^{ext\to int}
 \\ A^{int\to ext} 
& A^{ext}-\alpha\Id \end{array}\right),   \label{eq:A-block}
\EEQ
where (letting $A = (A_{\sigma',\sigma})_{\sigma',\sigma\in V}$): 
\BEQ A^{ext\to int} := (k_{\sigma_{ext}\to \sigma'})_{\sigma',\sigma_{ext}}, \qquad
 A^{int\to ext} := (k_{\sigma\to\sigma'_{ext}})_{\sigma'_{ext},\sigma} 
\EEQ
Here  $\sigma,\sigma'$, resp. $\sigma_{ext},\sigma'_{ext}$, range in ${\cal V}^{int}$, resp. ${\cal V}^{ext}$, and 
$A^{ext} = (k_{\sigma_{ext}\to \sigma'_{ext}})_{\sigma'_{ext},\sigma_{ext}\in {\cal V}^{ext}}$. 
Similarly, 
\BEQ {\cal W}(\alpha)  = \left(\begin{array}{cc} {\cal W}^{int}(\alpha)  &  {\cal W}^{int\to ext}(\alpha) 
 \\ {\cal W}^{ext\to int} (\alpha) 
& {\cal W}^{ext}(\alpha)  \end{array}\right), 
\EEQ
where (letting ${\cal W}(\alpha)  = (w(\alpha) _{\sigma\to\sigma'})_{\sigma,\sigma'\in {\cal V}}$): 
\BEQ {\cal W}^{int\to ext}(\alpha)  := (w(\alpha) _{\sigma\to \sigma'_{ext}})_{\sigma,\sigma'_{ext}}, \qquad
 {\cal W}^{ext\to int}(\alpha)  := (w(\alpha) _{\sigma_{ext}\to\sigma'})_{\sigma_{ext},\sigma'} 
\EEQ

\Bigskip We must first introduce some definitions.

\Medskip {\bf External transition rates} of ${\cal G}^{int}$ are the
total transition rates to ${\cal V}^{ext}$ from a given vertex 
in ${\cal V}^{int}$,
\BEQ k^{ext}_{\sigma} := \sum_{\sigma'_{ext}\in {\cal V}^{ext}}
k_{\sigma\to\sigma'_{ext}}.  \label{eq:kext}
\EEQ
Letting
\BEQ k_{\sigma}^{int} = \sum_{\sigma'\in {\cal V}^{int}} k_{\sigma\to\sigma'} \EEQ
be the total outgoing rate of the subgraph ${\cal G}^{int}$,  
we have the obvious identities relating the  outgoing rates of ${\cal G}^{int}$ and ${\cal G}$
and the deficiency rates:
\BEQ  k_{\sigma} = a_{\sigma} + \kappa_{\sigma} = k^{int}_{\sigma} + k_{\sigma}^{ext} , \qquad \sigma \in {\cal V}^{int}.
\EEQ

\Medskip {\bf Influx rates.}  Let 
\BEQ k_{\sigma_{ext}\to {\cal G}^{int}}:= \sum_{\sigma'\in {\cal V}^{int}} k_{\sigma_{ext}\to \sigma'}, \qquad
\sigma_{ext} \in {\cal V}^{ext}  \label{eq:influx-rate}
 \EEQ
be the total influx into ${\cal G}^{int}$ from the external vertex $\sigma_{ext}$.

\Medskip {\bf Exit probabilities.} The matrix $\tilde{\cal W}$, resp. ${\cal W}(\alpha)$, also describes a discrete-time Markov chain (resp. defective Markov chain) on ${\cal V} = {\cal V}^{int} \uplus 
{\cal V}^{ext}$ with set of absorbing states ${\cal V}^{ext}$.    Call $f(\sigma^*_{ext}) = 
( f^{\sigma}(\sigma^*_{ext}))_{\sigma\in {\cal V}}$
the unique solution of 
the linear  equation inside ${\cal V}^{int}$
\BEQ (-\Id + {\cal W}(\alpha) )f(\sigma^*_{ext})=0 \label{eq:exit-proba-eq}
\EEQ
 with
boundary conditions 
\BEQ f^{\sigma_{ext}}(\sigma^*_{ext})= \del_{\sigma_{ext},\sigma^*_{ext}}, \qquad \sigma_{ext} \in 
{\cal V}^{ext}.
\EEQ 
By analogy with standard Markov theory, the
vector $(f^{\sigma}(\sigma^*_{ext}))_{\sigma\in {\cal V}^{int}}$ may be interpreted as
the exit probability at $\sigma^*_{ext}\in {\cal V}^{ext}$ for the chain started at $\sigma\in{\cal V}^{int}$.

\Bigskip {\bf Some identities.}  From Section \ref{section:notations}, the following relations hold:
\BEQ   \tilde{\pi}^{int} \tilde{\cal W}^{int} = \tilde{\pi}, \qquad \tilde{A}^{int} \tilde{v}^{int} =0,
\label{eq:stationary}
 \EEQ
where $\tilde{v}_{\sigma}^{int}:= \frac{\tilde{\pi}^{int}_{\sigma}}{|\tilde{A}^{int}_{\sigma,\sigma}|}$, 
$\sigma\in {\cal V}^{int}$ is the stationary distribution of the time-continuous Markov chain associated
to $\tilde{A}^{int}$ with the usual normalization;

\BEQ \tilde{A}^{int}_{\sigma',\sigma}=A^{int}_{\sigma',\sigma}\qquad 
(\sigma'\not=\sigma), \qquad \tilde{A}^{int}_{\sigma,\sigma} = 
-\sum_{\sigma'\in {\cal V}^{int}, \sigma'\not=\sigma} A^{int}_{\sigma',\sigma} = -k^{int}_{\sigma};
\EEQ

\BEQ |\tilde{A}_{\sigma,\sigma}| = k_{\sigma} = |A_{\sigma,\sigma}| + \kappa_{\sigma} = 
|A_{\sigma,\sigma}|\, (1+\eps_{\sigma})
\EEQ

\BEQ |\tilde{A}^{int}_{\sigma,\sigma}| = |\tilde{A}_{\sigma,\sigma}| -
k^{ext}_{\sigma} = k_{\sigma}-k^{ext}_{\sigma} = k_{\sigma}^{int}.
\EEQ 
In particular,
\BEQ |\tilde{A}^{int}_{\sigma,\sigma}| = |A_{\sigma,\sigma}| + \kappa_{\sigma}-k_{\sigma}^{ext}
= |A_{\sigma,\sigma}|(1+\eps_{\sigma})-k_{\sigma}^{ext} \sim |A_{\sigma,\sigma}|  \label{eq:AAkk}
\EEQ
since $0\le \eps_{\sigma}\preceq 1$. 

\Medskip {\em Rewriting of the exit probability equation in terms
of the non-defective finite-time Markov generator.}
 We write the latter equation in terms of $\tilde{\cal W}^{int}$; if $\sigma\in {\cal V}^{int}, 
 \sigma'\in {\cal V}$,
\BEA  w(\alpha)_{\sigma\to \sigma'}  &=& \frac{k_{\sigma\to\sigma'}}{|A_{\sigma,\sigma}|+\alpha}  = \frac{k_{\sigma}}{|A_{\sigma,\sigma}|} \, \times\,  \frac{|A_{\sigma,\sigma}|}{|A_{\sigma,\sigma}|+\alpha} \, \times\, 
\frac{|\tilde{A}^{int}_{\sigma,\sigma}|}{k_{\sigma}} \, \times\, 
\frac{k_{\sigma\to\sigma'}}{|\tilde{A}^{int}_{\sigma,\sigma}|}
\nonumber\\
& \sim & \Big(1+ \eps_{\sigma} - \frac{k^{ext}_{\sigma} + \alpha}{k_{\sigma}}\Big)\, 
\tilde{w}^{int}_{\sigma \to \sigma'}.  \label{eq:5.14}
\EEA
  Now, the $\sigma$-coordinate
 equation $f^{\sigma}(\sigma^*_{ext})=({\cal W}(\alpha) f(\sigma^*_{ext}))^{\sigma}$ 
for exit probabilities may be rewritten
\BEA f^{\sigma}(\sigma^*_{ext}) &=& \Big( \sum_{\sigma'\in {\cal V}^{int}} w(\alpha)_{\sigma,\sigma'} f^{\sigma'}(\sigma^*_{ext}) \Big) + w(\alpha)_{\sigma\to\sigma^*_{ext}} \nonumber\\
&\sim &   \Big(1+ \eps_{\sigma} - \frac{k^{ext}_{\sigma}+\alpha}{k_{\sigma}}\Big)\, \sum_{\sigma'\in {\cal V}^{int}}  
\tilde{w}^{int}_{\sigma\to\sigma'} f^{\sigma'}(\sigma^*_{ext}) \  +  \frac{k_{\sigma\to\sigma^*_{ext}}}{k_{\sigma}+\alpha}   \label{eq:exit-proba-sim0}
\EEA
Conversely,
\BEQ ((\tilde{\cal W}^{int}-\Id)f(\sigma^*_{ext}))_{\sigma} \sim 
 \Big(- \eps_{\sigma} + \frac{k^{ext}_{\sigma}+\alpha}{k_{\sigma}}\Big)\,
 f^{\sigma}(\sigma^*_{ext}) \ - \frac{k_{\sigma\to\sigma^*_{ext}}}{k_{\sigma}+\alpha}  \label{eq:exit-proba-sim}
\EEQ

\Medskip {\em General identity.} For any function $f:{\cal V}\to\R$, we deduce from 
(\ref{eq:exit-proba-sim0}) by taking linear combinations over $\sigma^*_{ext}\in {\cal V}^{ext}$
the following general identity,
\BEQ (({\cal W}(\alpha)-\Id)f)^{\sigma} \sim ((\tilde{\cal W}^{int}-\Id)f)^{\sigma} + 
(\eps_{\sigma} - \frac{k^{ext}_{\sigma}+\alpha}{k_{\sigma}}) \sum_{\sigma'\in {\cal V}^{int}} 
\tilde{w}^{int}_{\sigma\to\sigma'} f^{\sigma'} \  +  \langle w(\alpha)_{\sigma\to\cdot}, f^{\cdot}\rangle_{ext} \label{eq:exit-proba-sim1}
\EEQ


\Bigskip {\bf Time-continuous adjoint problem.} In an analogous way, we call $v(\sigma^*_{ext}) = 
(v_{\sigma}(\sigma^*_{ext}))_{\sigma\in {\cal V}}$ the unique solution of the  adjoint,
time-continuous problem  inside ${\cal V}^{int}$
\BEQ A(\alpha)v(\sigma^*_{ext})=0 \label{eq:continuous-exit-proba-eq} \EEQ
with boundary conditions 
\BEQ v_{\sigma_{ext}}(\sigma^*_{ext}) = \del_{\sigma_{ext},\sigma^*_{ext}}, \qquad \sigma_{ext}\in {\cal V}^{ext}.
\EEQ
The above system of equations   inside ${\cal V}^{int}$ is equivalent to 
\BEQ (\tilde{A}^{int}v(\sigma^*_{ext}))_{\sigma} = \Big(-\kappa_{\sigma} + k^{ext}_{\sigma} + \alpha\Big) v_{\sigma}(\sigma^*_{ext}) - k_{\sigma^*_{ext}\to\sigma}  \label{eq:tcap}
\EEQ
and features the  influx rates $k_{\sigma^*_{ext}\to\sigma}$ instead of the outflux weights 
$\frac{k_{\sigma\to\sigma^*_{ext}}}{k_{\sigma}}$.



\subsection{Cut-off graphs: Lyapunov data and exit probability  estimates}  \label{subsection:cut-off-Lya-estimates}


The setting here is as in \S \ref{subsection:cut-off-setting}. We essentially prove a priori estimates
for the Lyapunov data. The last two paragraphs (\S \ref{subsubsection:exit-proba}, 
\S \ref{subsubsection:time-continuous-adjoint-pb}) help get the whole picture, in connection
with the probabilistic interpretation, the last one was actually used to find the renormalized generator, motivating Definition \ref{def:Aren};  Lemmas \ref{lem:exit-proba} 
and \ref{lem:time-cont-adj-pb} are, however, not used in the proof for the renormalization
algorithm.


\subsubsection{A priori estimates for Lyapunov exponent}  \label{subsubsection:apriori-estimates}


Though our main arguments are perturbative in essence, we prove here a priori estimates for
$\lambda^*$ in the autocatalytic case,  based on a technical result proved in Appendix. Let 
\BEQ \tau'_{{\cal G}}:= \Big(\min_{\sigma\in {\cal V}^{int}} k_{\sigma}\Big)^{-1}, \qquad \eps':= \max_{\sigma \in {\cal V}^{int}} \eps_{\sigma}, \qquad 
Z'(0)_{{\cal G}} := \max_{\sigma\in {\cal V}^{int}} (\frac{k^{ext}_{\sigma}}{k_{\sigma}})
\EEQ
Note that, in \S \ref{subsubsection:perturbation-Lyapunov} below, $\tau'_{{\cal G}}, \eps', 
Z'(0)_{{\cal G}}$ are replaced by  quantities $\tau_{{\cal G}}$, $\bar{\eps}_{\cal G}$,
$Z(0)_{{\cal G}}$ which are equivalent (i.e. have same scales) under hypotheses (i), (ii), see
\S \ref{subsection:cut-off-setting}. 
 We assume
\BEQ {\mathrm{(Autocata')}} :  \qquad \eps'\succ Z'(0)_{{\cal G}}  \label{eq:eps'>Z'} \EEQ

\Medskip We first consider the lower bound, $\lambda^*\ge \sup\{\alpha\ge 0 \ |\ \max_{\sigma^*}h(m,\alpha|\sigma^*)>1\}$, where $h(k,\alpha|\sigma^*):= \frac{k}{k+\alpha} \frac{|\tilde{A}_{\sigma^*,\sigma^*}|}{|\tilde{A}_{\sigma^*,\sigma^*}| - \kappa_{\sigma^*}+\alpha}$ $(k>0)$.
As proved in Appendix 1, this implies 
\BEQ \lambda^* \succeq y^2/x
\label{eq:apriorilowerbound-lambda}
\EEQ
where $m:= \min_{\sigma\in {\cal V}^{int}} |\tilde{A}_{\sigma,\sigma}|$ and   $x\equiv x(m|\sigma^*):= |A_{\sigma^*,\sigma^*}|+m$, $y^2\equiv y(m|\sigma^*)^2 = 4m(d|\tilde{A}_{\sigma^*,\sigma^*}| + \kappa_{\sigma^*}- k^{ext}_{\sigma^*})$, where $d:=\min_{\sigma\in {\cal V}^{int}} \Big(\frac{\kappa_{\sigma}-k^{ext}_{\sigma}}{|\tilde{A}_{\sigma,\sigma}|} \Big)$, $|d|\preceq 1$. We choose $\sigma^*$ such 
that $\eps'=\eps_{\sigma^*}$. Then $m\sim 1/\tau'_{{\cal G}}$, $x=x(m|\sigma^*)\sim |A_{\sigma^*,\sigma^*}|$,\  $\kappa_{\sigma^*} - k^{ext}_{\sigma^*}  \sim (\eps' - \frac{k^{ext}_{\sigma^*}}{k_{\sigma^*}})\,   |A_{\sigma^*,\sigma^*}| \sim \eps'  |A_{\sigma^*,\sigma^*}|$ as follows from the autocatalysis
assumption (\ref{eq:eps'>Z'}), whence $y^2 \sim 
 \eps' a_{\sigma^*}/\tau'_{{\cal G}}$ and 
\BEQ \lambda^*\succeq  \frac{\eps'}{\tau'_{{\cal G}}}. 
\EEQ 
As it turns out, this lower bound coincides with the estimate obtained in the perturbative
regime (see Lemma \ref{lem:Lya(1)}), which we prove only when $\eps\prec 1$ (equivalently,
$\eps'\prec 1$). 

\Medskip Next, we consider the upper bound, again from Appendix 1,
\BEQ \lambda^* \preceq  y^2/x
\label{eq:aprioriupperbound-lambda}
\EEQ
where $M:= \max_{\sigma\in {\cal V}^{int}} |\tilde{A}_{\sigma,\sigma}|$ and   $x\equiv x(M|\sigma^*):= |A_{\sigma^*,\sigma^*}| + M$, $y^2 \equiv y(M|\sigma^*)^2 = 4M(D |\tilde{A}_{\sigma^*,\sigma^*}| + \kappa_{\sigma^*}- k^{ext}_{\sigma^*})$, where $D:=\max_{\sigma\in {\cal V}^{int}} \Big(\frac{\kappa_{\sigma}-k^{ext}_{\sigma}}{|\tilde{A}_{\sigma,\sigma}|} \Big)$, $|D|\preceq 1$. Choosing $\sigma^*$ such 
that $k_{\sigma^*} \sim |\tilde{A}_{\sigma^*,\sigma^*}| \sim 1/\tau'_{{\cal G}}$, we get
\BEQ \lambda^*\preceq  \frac{1}{\tau'_{{\cal G}}} . \EEQ 
 
\Medskip When $\eps'\sim 1$ (a regime not covered by the perturbative arguments of \S 
\ref{subsubsection:perturbation-Lyapunov}), we get $\lambda^* \sim \frac{1}{\tau'_{{\cal G}}}. $
Otherwise $(\eps'\prec 1)$, the upper bound proves insufficient, but our arguments below show
that $\lambda^*\sim \frac{\eps'}{\tau'_{{\cal G}}}$.


\subsubsection{Perturbation theory arguments for Lyapunov exponent and weight}  \label{subsubsection:perturbation-Lyapunov}


We look here for approximate expansions of the form $\lambda^* = \lambda_{(0)} + \lambda_{(1)} + 
\cdots$ (Lyapunov exponent),   $v^* = v_{(0)} +  v_{(1)} + \cdots$ (right Lyapunov
eigenvector), $\pi^* = \pi_{(0)} +  \pi_{(1)} + \cdots$ (Lyapunov weights), and  $v^{\dagger} = v^{\dagger}_{(0)} + v^{\dagger}_{(1)} + \cdots$ (right Lyapunov
eigenvector), where terms of order (1) are smaller than the main terms of order (0) by a small  parameter $O(1/b)$ which converges to 0 in the limit $b\to\infty$.  

\Medskip We may perturb either $A$ (continuous-time adjoint generator) or ${\cal W}$ (discrete-time
generator). Ultimately we are interested in the top eigenvalue of $A$, so perturbing $A$ seems more
direct and easier. It is, however, easier to bound the pseudo-inverses appearing in perturbation theory
in the discrete-time case, using an old Markov chain argument attributed to Doeblin. The two expansions
are different; however, they give
the same value for the leading order of the Lyapunov exponent. We start with evaluating $\lambda_{(1)}$
by a perturbation argument
based on $A$, mainly as an appetizer, because the result is so simple. We then abandon this path and 
develop a chain of perturbation arguments based on ${\cal W}$, which lays the ground for the
rigorous fixed-point argument detailed in Appendix.

\Medskip {\bf First perturbation argument.} We rewrite the eigenvalue problem $A^{int}v^*=\lambda^* v^*$ in the form
\BEQ ({\cal H}_{(0)} + {\cal H}_{(1)})v^* = \lambda^* v^*, \qquad \lambda^* =\lambda_{(0)} + \lambda_{(1)}+\cdots, \qquad v^* =
v_{(0)} + v_{(1)} + \cdots
\EEQ
where ${\cal H}_{(0)} := \tilde{A}^{int}$, ${\cal H}_{(1)} := A^{int}-\tilde{A}^{int}$. Since the Lyapunov
eigenvalue of $\tilde{A}^{int}$ is $0$, we may actually set $\lambda_{(0)}=0$. We let $v_{(0)}:=\tilde{v}^{int} \in 
\R_+^{{\cal V}^{int}}$
(see (\ref{eq:stationary}) above) be the right Lyapunov eigenvector of $\tilde{A}^{int}$, 
normalized by imposing $||\tilde{v}^{int}||_1 := \sum_{\sigma\in {\cal V}^{int}} 
\tilde{v}^{int}_{\sigma}=1$,  and $u_{(0)}:= {\bf 1}$
(vector with constant coefficients $1$) the associated left Lyapunov eigenvector. The normalization is chosen in such a way
that $\langle u_{(0)} , v_{(0)} \rangle:=\sum_{\sigma\in {\cal V}^{int}} (u_{(0)})^{\sigma} (v_{(0)})_{\sigma} = 1$.   Elementary non-Hermitian Rayleigh
perturbation arguments (see Appendix, Section \ref{section:app3-pert}) yield
\BEQ \lambda_{(1)} = \langle u_{(0)}, {\cal H}_{(1)} v_{(0)}\rangle
\EEQ
The first-order perturbation ${\cal H}_{(1)}$ is a diagonal matrix, with $({\cal H}_{(1)})_{\sigma,\sigma} = 
\kappa_{\sigma} - k^{ext}_{\sigma}$. 
Using the relations proved in \S \ref{subsection:cut-off-setting}, we get
\BEQ \lambda_{(1)} =  \frac{ \sum_{\sigma\in {\cal V}^{int}} \frac{\kappa_{\sigma}-k^{ext}_{\sigma}}{k_{\sigma} -
k^{ext}_{\sigma}} \, \tilde{\pi}^{int}_{\sigma}}{\sum_{\sigma\in {\cal V}^{int}} \frac{\tilde{\pi}^{int}_{\sigma}}{k_{\sigma} - k^{ext}_{\sigma}} }.
\EEQ

By assumption, $k_{\sigma}\succ k^{ext}_{\sigma}$. Replacing $k_{\sigma}-k^{ext}_{\sigma}$ by
$k_{\sigma}$ in the denominators, we have proved the following:

\begin{Lemma}[first-order formula for Lyapunov eigenvalue]  \label{lem:Lya(1)}
Let 
\BEQ \tau_{{\cal G}}:= \sum_{\sigma\in {\cal V}^{int}} \frac{\tilde{\pi}^{int}_{\sigma}}{k_{\sigma}};  \label{eq:tauGint}
\EEQ
\BEQ \bar{\eps}_{{\cal G}}:=  \sum_{\sigma\in {\cal V}^{int}} \tilde{\pi}^{int}_{\sigma} 
 \eps_{\sigma}, \qquad Z(\eps,\alpha)_{{\cal G}}:= \sum_{\sigma\in {\cal V}^{int}} \tilde{\pi}^{int}_{\sigma}  \Big( \frac{k^{ext}_{\sigma}+\alpha}{k_{\sigma}} 
 -\eps_{\sigma} \Big)  \label{eq:epsZ(0)}
 \EEQ
Then  
\BEQ \lambda_{(1)} \sim   -\frac{Z(\eps,0)_{{\cal G}}}{\tau_{{\cal G}}}  \label{eq:lambda1-gen-bis}
\EEQ
\end{Lemma}

\noindent We sometimes use the short notation $\langle f_{\cdot}\rangle_{int}$ for the averaged quantity $\sum_{\sigma\in {\cal V}^{int}} \tilde{\pi}^{int}_{\sigma} f_{\sigma}$, so that 
\BEQ \tau_{{\cal G}} = \langle \frac{1}{k_{\cdot}}\rangle_{int}, \qquad 
\bar{\eps}_{{\cal G}} = \langle
\eps_{\cdot}\rangle_{int},\qquad Z(\eps,\alpha)_{{\cal G}} = \langle \frac{k^{ext}_{\cdot}+\alpha}{k_{\cdot}} - \eps_{\cdot}\rangle_{int}
\EEQ

\noindent {\bf Remark.} By definition, $Z(0)_{{\cal G}}:=Z(0,0)_{{\cal G}} = 
 \sum_{\sigma\in {\cal V}^{int}} \tilde{\pi}^{int}_{\sigma}   \frac{k^{ext}_{\sigma}}{k_{\sigma}}$
 and 
\BEQ Z(\eps,\alpha)_{{\cal G}} = Z(0)_{{\cal G}} + \tau_{{\cal G}}\alpha - \bar{\eps}_{{\cal G}}. \label{eq:Z(eps,alpha)}
\EEQ
Thus (\ref{eq:lambda1-gen}) may be rewritten: 
\BEQ \lambda_{(1)} \sim   \frac{1}{\tau_{{\cal G}}} \Big( \bar{\eps}_{{\cal G}} - Z(0)_{{\cal G}}\Big).   \label{eq:lambda1-gen}
\EEQ
  Note also that 
\BEQ Z(\eps,\bar{\eps}_{{\cal G}}/\tau_{\cal G})  = Z(0)_{{\cal G}} 
\EEQ
while 
\BEQ Z(-\eps,0)_{\cal G} =  Z(\eps, 2\, \bar{\eps}_{{\cal G}}/\tau_{\cal G}) =  Z(0)_{\cal G} + \bar{\eps}_{\cal G} \sim Z_G,   \label{eq:Z-eps} 
\EEQ
see (\ref{eq:ren-weight}), 
a sum of two positive terms, bounds $|Z(\eps,0)_{\cal G}|$ from above.
Note that, by construction, $0\le \bar{\eps}_{\cal G}\le 1$ since all $\eps_{\sigma}\in [0,1]$, and 
$Z(0)_{\cal G}\prec 1$, thus $\lambda_{(1)}\succ 0$ when $\bar{\eps}_{\cal G}\sim 1$ 
(equivalently, when $Z(-\eps,0)_{\cal G}\sim 1$), a special case
discussed in \S \ref{subsubsection:apriori-estimates} which we exclude henceforth. The two
quantities $Z(-\eps,0)_{\cal G}$ and $Z_G$ play the same role in this section; we use rather
$Z(-\eps,0)_{\cal G}$, because renormalization is performed only later on (\S \ref{subsection:renormalization}). 

\Medskip If the first term in the r.-h.s. of (\ref{eq:lambda1-gen})
is of strictly higher scale than the second one, 
then $\lambda_{(1)}>0$, and 
\BEQ \lambda_{(1)} \sim  \frac{\bar{\eps}_{\cal G}}{\tau_{{\cal G}}}   \label{eq:lambda1-eps-tauGint}
\EEQ
If, conversely,  the {\em second} term in the r.-h.s. of (\ref{eq:lambda1-gen})
is of strictly higher scale than the {\em first} one, then $\lambda_{(1)}<0$. If both
terms are of the same scale, then the sign of $\lambda_{(1)}$ cannot be determined from these
estimates.

\Medskip {\bf A key fact.} Under properties (i), (ii) (see \S \ref{subsection:cut-off-setting}),  {\bf all coefficients of $\tilde{\pi}^{int}$ are of order 1.}  Namely, up to the 'averaging procedure' described in the proof of Corollary
\ref{cor:eta}, one may assume that all $\tilde{w}^{int}_{\sigma\to\sigma'}\sim 1$, $\sigma\not=
\sigma'\in {\cal V}^{int}$. Let then $\underline{\sigma}\not= \overline{\sigma}$ s.t. 
$\tilde{\pi}^{int}_{\underline{\sigma}} = \min_{\sigma}(\tilde{\pi}^{int}_{\sigma}), \ 
\tilde{\pi}^{int}_{\overline{\sigma}} = \max_{\sigma}(\tilde{\pi}^{int}_{\sigma})$. The stationary
equation $\tilde{\pi}^{int}_{\underline{\sigma}} = \sum_{\sigma'\not=\underline{\sigma}} \tilde{\pi}^{int}_{\sigma'} \tilde{w}^{int}_{\sigma'\to\underline{\sigma}} \sim \max_{\sigma'\not = \underline{\sigma}} (\tilde{\pi}^{int}_{\sigma'})$ implies $\tilde{\pi}^{int}_{\underline{\sigma}}
\sim \tilde{\pi}^{int}_{\overline{\sigma}}$, from which 
\BEQ \tilde{\pi}^{int}_{\sigma}\sim 1 \label{eq:tildepiintsigma=1}
\EEQ
 for 
all $\sigma$.

\Medskip Then
\BEQ \tau_{{\cal G}}\sim \Big(\min_{\sigma\in  {\cal V}^{int}}(k_{\sigma}) \Big)^{-1}, \EEQ
\BEQ \eps\sim \max_{\sigma\in {\cal V}^{int}}(\eps_{\sigma}), \qquad Z(0)_{{\cal G}} \sim \max_{\sigma\in {\cal V}^{int}} (\frac{k^{ext}_{\sigma}}{k_{\sigma}})
\EEQ


\Bigskip  We set apart an easy subcase, namely, when $\eps'\sim 1\succ Z(0)_{{\cal G}}$, see \S \ref{subsubsection:apriori-estimates}. This
autocatalytic subcase was already  discussed in \S \ref{subsubsection:apriori-estimates};
it was proved there that $\lambda^*\sim 1/\tau_{{\cal G}}$ is less than or comparable to
all internal kinetic rates.  Thus $w(\lambda^*)_{\sigma\to\sigma'}
\sim \tilde{w}_{\sigma\to\sigma'}$. Therefore, solving the equations $\sum_{\sigma'\not=\sigma}
\pi^*_{\sigma'} w(\lambda^*)_{\sigma'\to\sigma} = \pi^*_{\sigma}$ (see just below) for the 
Lyapunov weights yields $\pi^*_{\sigma}\sim 1$ as above (\ref{eq:tildepiintsigma=1}). This makes
it possible to assume in the following that $\bar{\eps}_{\cal G}j,\prec 1$.



\Bigskip {\bf Second perturbation argument.} From now on, we develop perturbation arguments based on
the discrete-time Markov chain;  an application of the implicit function
theorem (see Appendix) makes the argument below rigorous.  Recall from (\ref{eq:pi*}) that
 $Av^* = \lambda^* v^*$ if and only if $\sum_{\sigma\not=\sigma'}
\pi^*_{\sigma} w(\lambda^*)_{\sigma\to\sigma'} = \pi^*_{\sigma'}$, where 
\BEQ  w(\lambda^*)_{\sigma\to\sigma'} = 
\frac{k_{\sigma\to\sigma'}}{|A_{\sigma,\sigma}|+\lambda^*}, \qquad 
\pi^*_{\sigma} := (|A_{\sigma,\sigma}|+\lambda^*) v^*_{\sigma}. \label{eq:2nd-pert-wstarpistar}
\EEQ
  Recall $|A_{\sigma,\sigma}| = k_{\sigma}-\kappa_{\sigma}$.  We rewrite ${\cal W}(\lambda^*)$ as a perturbation of  $\widetilde{\cal W}^{int}$ with coefficients $\tilde{w}^{int}_{\sigma\to\sigma'} = \frac{k_{\sigma\to\sigma'}}{k_{\sigma}^{int}} =
 \frac{k_{\sigma\to\sigma'}}{k_{\sigma} - k^{ext}_{\sigma}}$:
\BEQ w(\lambda^*)_{\sigma\to\sigma'} =  \frac{k_{\sigma\to\sigma'}}{(k_{\sigma}-k^{ext}_{\sigma}) + 
(k^{ext}_{\sigma}-\kappa_{\sigma} +\lambda^*)} \sim \tilde{w}^{int}_{\sigma\to\sigma'} ( 1 + \eps_{\sigma} - 
\frac{k^{ext}_{\sigma} + \lambda^*}{k_{\sigma}} )   \label{eq:w(lambda*)sim3}
\EEQ 
We rewrite the eigenvalue problem $ ^t{\cal W}(\lambda^*) \pi^* = \pi^* $ in the 
form 
\BEQ ({\cal H}_{(0)} + {\cal H}_{(1)}) (\pi_{(0)} + \pi_{(1)}) =0, 
 \label{eq:second-perturbation-argument}
\EEQ
 with
\BEQ  \pi^* \sim \pi_{(0)} + \pi_{(1)}, \qquad  \lambda^* \sim 
\lambda_{(1)};  \EEQ
\BEQ {\cal H}_{(0)} = \,   ^t\widetilde{\cal W}^{int}-\Id, \qquad 
\pi_{(0)} = \tilde{\pi}^{int} \EEQ
\BEQ ({\cal H}_{(1)})_{\sigma,\sigma'} \sim \tilde{w}^{int}_{\sigma'\to\sigma} (\eps_{\sigma'} - 
\frac{k^{ext}_{\sigma'} + \lambda_{(1)} }{k_{\sigma'}}) ,
\label{eq:H1}
 \EEQ
 which is first order in $\lambda_{(1)}$.
This
defines  a fixed-point equation for $\lambda^*$ which has a unique solution with our hypotheses; 
see Appendix, \S \ref{app: fixed1}.  Solving to first order, we neglect ${\cal H}_{(1)}\pi_{(1)}$
in  (\ref{eq:second-perturbation-argument}). Since 
\BEQ {\cal H}_{(0)} \pi_{(0)}=0; \qquad 
\langle {\bf 1}, {\cal H}_{(0)}(\pi_{(0)}+\pi_{(1)})\rangle = \langle \, ^t {\cal H}_{(0)}
{\bf 1}, \pi_{(0)}+\pi_{(1)}\rangle =0, 
\EEQ
we find first $\langle {\bf 1}, {\cal H}_{(1)}\pi_{(0)}\rangle =0$, from which
\BEQ \lambda_{(1)} \sim - \frac{\langle {\bf 1}, {\cal H}_{(1)}(\lambda_{(1)}=0) \pi_{(0)}\rangle}
{\langle {\bf 1}, \partial_{\lambda_{(1)}}{\cal H}_{(1)}(\lambda_{(1)}=0) \pi_{(0)}\rangle}
\label{eq:lambda1}
\EEQ  
Using the explicit expression (\ref{eq:H1}) for ${\cal H}_{(1)}$, we obtain
\BEQ \lambda_{(1)} \sim \frac{\sum_{\sigma',\sigma} \tilde{w}^{int}_{\sigma'\to\sigma} (\eps_{\sigma'} - \frac{k^{ext}_{\sigma'}}{k_{\sigma'}}) \tilde{\pi}^{int}_{\sigma'}}
 { \sum_{\sigma',\sigma} 
 \tilde{w}^{int}_{\sigma'\to\sigma}  \frac{\tilde{\pi}^{int}_{\sigma'}}{k_{\sigma'}} } \nonumber\\
\EEQ
and get back (\ref{eq:lambda1-gen})
using probability conservation $\sum_{\sigma} \tilde{w}^{int}_{\sigma'\to\sigma}=1$. 
Note the related identity from (\ref{eq:w(lambda*)sim3}), valid for any small $\lambda$, 
\BEQ \langle ({\cal W}(0) - {\cal W}(\lambda)) {\bf 1}, \tilde{\pi}^{int}\rangle  \sim 
\lambda \sum_{\sigma,\sigma'} \frac{\tilde{w}^{int}_{\sigma\to\sigma'}}{k_{\sigma}} 
\tilde{\pi}^{int}_{\sigma} = \lambda \sum_{\sigma'} \frac{\tilde{\pi}^{int}_{\sigma}}{k_{\sigma}}
\sim \lambda \tau_{\cal G}   \label{eq:related-id1}  \EEQ
In particular,
\BEQ \langle ({\cal W}(0) - {\cal W}(\lambda_{(1)})) {\bf 1}, \tilde{\pi}^{int}\rangle  \sim 
\bar{\eps}_{{\cal G}}-Z(0)_{\cal G}.  \label{eq:related-id2}
\EEQ

\Bigskip  Next, we get (see Appendix \ref{section:app3-pert})
\BEQ \pi_{(1)} = - ( \, ^t \tilde{\cal W}^{int}-\Id)^{-1} {\cal H}_{(1)} \pi_{(0)} 
\label{eq:pi1}
\EEQ
which makes sense because  $\langle {\bf 1}, {\cal H}_{(1)}\pi_{(0)}\rangle =0$.
 We now prove that, under suitable assumptions,  $\pi_{(1)}$ is  small compared to $\pi_{(0)} = \tilde{\pi}^{int}$. For this, we must bound the pseudo-inverse 
\BEQ -( \, ^t \tilde{\cal W}^{int}-\Id)^{-1} = \sum_{n\ge 0} (\, ^t\tilde{\cal W}^{int})^n.
\label{eq:pseudo-inverse-series}
\EEQ
The sum converges on the subspace $\{u= (u_{\sigma})_{\sigma\in {\cal V}^{int}} \ |\ 
\langle {\bf 1}, u\rangle =0$. 
 We use a simple, classical argument 
attributed to W. Doeblin. We say that the vector  $u=(u_{\sigma})_{\sigma\in {\cal V}^{int}}$
has zero average if 
\BEQ u^{av}:= \frac{1}{|{\cal V}^{int}|} \sum_{\sigma\in {\cal V}^{int}} u_{\sigma} \EEQ
vanishes.  Note that the zero-average space $\{u\ |\ u^{av}=0\}$ is preserved 
by the matrix $^t \tilde{\cal W}^{int}$; it is also preserved by the pseudo-inverse
$( \, ^t \tilde{\cal W}^{int}-\Id)^{-1}$. Dualizing, we say that the vector $f=(f^{\sigma})_{\sigma\in {\cal V}^{int}}$ as zero average w.r. to
$\tilde{\pi}^{int}$ if 
\BEQ f_{av}:= \sum_{\sigma\in {\cal V}^{int}} \tilde{\pi}^{int}_{\sigma} f^{\sigma} \EEQ
vanishes;  the  space $\{f\ |\ f_{av}=0\}$ is preserved 
by the matrix $\tilde{\cal W}^{int}$.   The relevant norms for (zero-average or not) vectors on ${\cal V}^{int}$ are the $L^1$ and $L^{\infty}$ norms,
\BEQ ||u||_1 := \sum_{\sigma\in {\cal V}^{int}} |u_{\sigma}|, \qquad  
   ||f||_{\infty} := \sup_{\sigma\in {\cal V}^{int}} |f^{\sigma}|
\EEQ
if $u = (u_{\sigma})_{\sigma\in {\cal V}^{int}}, \ f=(f^{\sigma})_{\sigma\in {\cal V}^{int}}$. 
The semi-norm $||f||^*_{\infty} :=  \sup_{u:{\cal V}^{int}\to\R \ |\  u^{av}=0, ||u||_1=1}  
|\langle f,u\rangle|$ will also be used. The two norms $||\cdot||_{\infty}, ||\cdot||^*_{\infty}$
are equivalent when restricted to $\{f\ |\ f_{av}=0\}$ since 
\BEQ ||f||^*_{\infty}\le ||f||_{\infty} = 
\max_{\sigma} |\langle f,\del^{\sigma}\rangle| = \max_{\sigma} |\langle f,\del^{\sigma}-\tilde{\pi}^{int}\rangle| \le ||f||^*_{\infty} ||\del^{\sigma}-\tilde{\pi}^{int}||_1 \le 2 ||f||^*_{\infty};
\label{eq:finfty*2finfty}
\EEQ 
we have used: $(\del^{\sigma}-\tilde{\pi}^{int})^{av}=0$ and $\langle f,\tilde{\pi}^{int}\rangle=f_{av}=0$. 

\Medskip
 Define 
\BEQ \rho\equiv \rho(\tilde{\cal W}^{int}):= \sum_{\sigma'} \min_{\sigma} (\tilde{w}^{int}_{\sigma\to\sigma'})   \label{eq:rho} \EEQ
The Lemma below may be stated by saying that the {\em spectral gap} of $\tilde{\cal W}^{int}$ is
$\ge \rho$.
We assume that $\rho>0$. Later on, we shall see how this condition is automatically ensured by
a suitable averaging procedure.

\begin{Lemma}[$L^1$-contraction bound for $^t\tilde{\cal W}^{int}$]  \label{lem:Doeb1}
For every zero-average $u=(u_{\sigma})_{\sigma\in {\cal V}^{int}}$, 
\BEQ ||\, ^t\tilde{\cal W}^{int} u||_1 \le (1-\rho) ||u||_1.
\EEQ 
\end{Lemma}

\noindent  


\begin{Lemma}[$L^{\infty}$-contraction bound for $\tilde{\cal W}^{int}$]   \label{lem:Doeb2}
For every  $f=(f^{\sigma})_{\sigma\in {\cal V}^{int}}$ such that $f_{av}=0$, 
\BEQ ||\tilde{\cal W}^{int} f ||^*_{\infty} \le (1-\rho) \, ||f||^*_{\infty}.
\EEQ 
\end{Lemma}

The proofs of Lemmas \ref{lem:Doeb1} and \ref{lem:Doeb2} may be found in Appendix.  We now apply 
the $L^1$-contraction bound for $^t \tilde{\cal W}^{int}$ to bound $(\pi_{(1)})_{\sigma}$, using
the definition (\ref{eq:pi1}). The spectral gap $\rho$ is as in (\ref{eq:rho}).  Recall
from (\ref{eq:Z-eps}) that $Z(-\eps,0)_{\cal G}$ bounds $|Z(\eps,0)_{\cal G}|$ from above.

\begin{Lemma}  \label{lem:eta}
Let
\BEQ \eta := \max_{\sigma,\sigma'\in {\cal V}^{int}} \frac{Z(-\eps,0)_{\cal G}}{ \tilde{w}^{int}_{\sigma\to\sigma'}} 
\EEQ
Then, for every $\sigma'\in {\cal V}^{int}$, 
\BEQ |(\pi_{(1)})_{\sigma}| \preceq \frac{\eta}{\rho} \tilde{\pi}^{int}_{\sigma}. \EEQ
\end{Lemma}

\noindent{\bf Proof.} Expanding in series as in (\ref{eq:pseudo-inverse-series}) and using Lemma 
\ref{lem:Doeb1} iteratively, we get 
\BEQ ||(\Id-\, ^t\tilde{\cal W}^{int})^{-1} {\cal H}_{(1)}\pi_{(0)} ||_1 \le \rho^{-1} ||{\cal H}_{(1)}\pi_{(0)}||_1, \label{eq:1/rho}
\EEQ
 whence 
\BEA \sum_{\sigma'} |(\pi_{(1)})_{\sigma'}| &\le& \frac{1}{\rho} \sum_{\sigma'} |({\cal H}_{(1)} \pi_{(0)})_{\sigma'}| \preceq \frac{1}{\rho} S_1,  \label{eq:S1S2}
\EEA
with 
\BEQ S_1:=  \sum_{\sigma,\sigma'} \tilde{w}^{int}_{\sigma\to\sigma'} |\eps_{\sigma} - \frac{k^{ext}_{\sigma} + \lambda_{(1)}}{k_{\sigma}}| \, \tilde{\pi}^{int}_{\sigma} \preceq \sum_{\sigma} Z(-\eps,0)_{\cal G} \tilde{\pi}^{int}_{\sigma} 
\EEQ
 Fix $\sigma'$,  then
\BEQ \rho\,  |(\pi_{(1)})_{\sigma'}| \preceq S_1 \le \eta \sum_{\sigma} \tilde{\pi}^{int}_{\sigma}
\tilde{w}^{int}_{\sigma\to \sigma'} = \eta \tilde{\pi}^{int}_{\sigma'}. 
\EEQ
\hfill \eop


\begin{Corollary} \label{cor:eta}
Assume $\bar{\eps}_{\cal G}\prec 1$. Then 
\BEQ |(\pi_{(1)})_{\sigma'}| \preceq Z(-\eps,0)_{\cal G} \, \tilde{\pi}^{int}_{\sigma'}, \qquad \sigma'\in {\cal V}^{int} \EEQ
with $Z(-\eps,0)_{\cal G}\prec 1$, so that $\pi_{(1)}$ is  small compared to $\pi_{(0)}$: for all $\sigma'\in {\cal V}^{int}$,  $(\pi_{(0)} + \pi_{(1)})_{\sigma'} \sim \tilde{\pi}^{int}_{\sigma'}$. 
\end{Corollary}

{\bf Proof.}  In general, $\rho$ may be zero (in particular, when ${\cal G}^{int}$ 
is not aperiodic) or too small; we now define a suitable 
'averaging procedure' ensuring that $\rho>0$ whenever ${\cal G}^{int}$ is irreducible, and $\rho\sim 1$ 
when, furthermore, ${\cal G}^{int}$ is dominant.   First, one substitutes the
averaged, aperiodic $\tilde{\cal W}^{int}_{av}:= \frac{1}{p} \sum_{\theta=1}^p (\tilde{\cal W}^{int})^{\theta}$ to
$\tilde{\cal W}^{int}$, where $p=$ gcd$\{n\ge 1\ |\ ((\tilde{\cal W}^{int})^n)_{\sigma,\sigma}>0\}$ is the
minimal period of $\tilde{\cal W}^{int}$. Then, it is well-known that for $q$ large enough, $\tilde{\cal W}^{int}_* := (\tilde{\cal W}^{int}_{av})^q$
is such that all its coefficients $(\tilde{\cal W}^{int}_*)_{\sigma,\sigma'}$ are
$>0$.  If ${\cal G}^{int}$ is dominant, then one can consider $q$ minimal such that   $\tilde{\cal W}^{int}_*$ connects each pair $(\sigma,\sigma')$ by dominant edges, so that $(\tilde{\cal W}^{int}_*)
_{\sigma,\sigma'}\sim 1$;  as a consequence, $\rho_*:= \rho(\tilde{\cal W}^{int}_*)\sim 1$. Working with 
$\tilde{\cal W}^{int}_*$ instead of $\tilde{\cal W}^{int}$, we get instead of (\ref{eq:1/rho})

\BEQ ||(\Id-\, ^t{\cal W}^{int}_*)^{-1} u ||_1 \le  \frac{1}{\rho_*}  ||u ||_1  \label{eq:inverse-tW-bound}
\EEQ
for every zero average $u = (u_{\sigma})_{\sigma\in {\cal V}^{int}}$, in particular for $u=  {\cal H}_{(1)}\pi_{(0)}$. 
Now, for every vector  $v$,  
\BEA || (\Id - \, ^t {\cal W}^{int}_*)v||_1 &\le & \sum_{\theta=0}^{q-1} || (\, ^t \tilde{\cal W}_{av}^{int})^{\theta}
(\Id - \, ^t \tilde{\cal W}^{int}_{av}) v||_1 \le q\,  || (\Id-\, ^t \tilde{\cal W}^{int}_{av}) v ||_1 
\nonumber\\
&\le & \frac{q}{p} \sum_{\theta=1}^{p} || (\Id - (\, ^t \tilde{\cal W}^{int})^{\theta} ) v ||_1 
\le \frac{q(p+1)}{2}\,  ||   (\Id - \, ^t {\cal W}^{int}) v ||_1 \nonumber\\
\EEA
from which, 
\BEQ ||((\Id-\, ^t{\cal W}^{int})^{-1} {\cal H}_{(1)}\pi_{(0)} ||_1 \le (\rho^*)^{-1} \,  || {\cal H}_{(1)}\pi_{(0)} ||_1 
\EEQ
where $\rho^*:= \frac{\rho_*}{q(p+1)/2}$. 
Redefining $\eta$ (see Lemma \ref{lem:eta}) as $\eta^*:= \max_{\sigma,\sigma'\in {\cal V}^{int}} \Big(\frac{ \eps_{\sigma} + (k^{ext}_{\sigma}/k_{\sigma})}{(\tilde{w}^{int}_*)_{\sigma\to\sigma'}}
\Big) \sim Z(-\eps,0)_{\cal G}$ and going through
the details of the proof of the Lemma, one gets $|(\pi_{(1)})_{\sigma'}|\preceq \frac{\eta^*}{\rho^*} \tilde{\pi}^{int}_{\sigma'}$. 
\hfill \eop


\subsubsection{Perturbation theory arguments for adjoint Lyapunov vector}  \label{subsubsection:perturbation-adjoint-Lyapunov}


Proceeding as in the second perturbation argument above, see in particular (\ref{eq:second-perturbation-argument})--(\ref{eq:pi1}), with  ${\cal H}_{(0,1)}$ replaced
 by their adjoints, and a Lyapunov vector $v^{\dagger}$ chosen as a perturbation of
  $v^{\dagger}_{(0)}={\bf 1}$, we prove a result similar to 
Corollary \ref{cor:eta}, with the same hypotheses.

\begin{Corollary}[adjoint Lyapunov vector] \label{cor:vdagger}
Assume $\bar{\eps}_{\cal G}\prec 1$. Then 
 \BEQ ||v^{\dagger}_{(1)}||_{\infty} \preceq Z(-\eps,0)_{\cal G} \prec 1 \EEQ
so that $v^{\dagger}_{(1)}$ is small compared to $v^{\dagger}_{(0)}$. 
\end{Corollary}

We briefly sketch the proof;  $v^{\dagger}= v^{\dagger}_{(0)}+v^{\dagger}_{(1)}$ is 
constructed as a perturbation of the constant vector $v^{\dagger}_{(0)}={\bf 1}$, with  
\BEQ v^{\dagger}_{(1)} = - (\tilde{\cal W}^{int}-\Id)^{-1} \  ^t{\cal H}_{(1)} v^{\dagger}_{(0)} 
\label{eq:vdagger1}
\EEQ
The proof uses Lemma \ref{lem:Doeb2} instead of \ref{lem:Doeb1}. Instead of (\ref{eq:1/rho}), 
(\ref{eq:S1S2}), we get:  $||(\Id-\, \tilde{\cal W}^{int})^{-1} \  ^t {\cal H}_{(1)}
v^{\dagger}_{(0)} ||^*_{\infty} \le \frac{1}{\rho^*} \, ||\, ^t {\cal H}^{(1)} v^{\dagger}_{(0)}||^*_{\infty}$. Now, $||\,  ^t {\cal H}^{(1)} v^{\dagger}_{(0)}||^*_{\infty} = \sup_{u \ |\ u^{av}=0, \, ||u||_1=1}
|\langle \, ^t{\cal H}_{(1)}v^{\dagger}_{(0)},u\rangle|$. For such $u$,  
\BEA
|\langle \, ^t{\cal H}_{(1)}v^{\dagger}_{(0)},u\rangle| &=& |\langle v^{\dagger}_{(0)}, 
{\cal H}_{(1)} u\rangle| = |\sum_{\sigma'} {\cal H}_{(1)}u_{\sigma'}| \nonumber\\
& \le & \sum_{\sigma,\sigma'}  \tilde{w}^{int}_{\sigma\to\sigma'}
 |\eps_{\sigma}-\frac{k^{ext}_{\sigma} + \lambda_{(1)}}{k_{\sigma}}| \, |u_{\sigma}|  \le | |\eps_{\cdot}-\frac{k^{ext}_{\cdot} + \lambda_{(1)}}{k_{\cdot}}||_{\infty} \preceq \eta^*
 \nonumber\\
\EEA
 whence (using (\ref{eq:finfty*2finfty}))
$ ||v^{\dagger}-{\bf 1}||_{\infty} \le   2\frac{\eta^*}{\rho^*}. $   \hfill \eop


\subsubsection{Exit probabilities}  \label{subsubsection:exit-proba}


We estimate here exit probabilities by solving   equation (\ref{eq:exit-proba-eq}) using $-( \,  \tilde{\cal W}^{int}-\Id)^{-1} = \sum_{n\ge 0} (\, \tilde{\cal W}^{int})^n$ and the $L^{\infty}$-contraction bond, see Lemma \ref{lem:Doeb2}. The value of $\alpha$ is important; it must be 
chosen larger than the Lyapunov exponent $\lambda^*$ for the operator ${\cal W}(\alpha)-\Id$ 
to be invertible, and computations produce the quantity $Z(\eps,\alpha)_{\cal G}>0$ in 
the denominator. On the other hand, error terms are of order (at most) $Z(-\eps,0)_{\cal G}$; 
see (\ref{eq:Z-eps}). 
We choose $\alpha$ such that $Z(\eps,\alpha)_{\cal G}$ is of the same order as $Z(-\eps,0)_{\cal G}$.
In the non-autocatalytic case, this holds true if $\tau_{\cal G}|\alpha|\prec Z(0)_{\cal G}$, including
the trivial case $\alpha=0$. In the autocatalytic case, we need $\tau_{\cal G}\alpha>0$ to be comparable to
$Z(-\eps,0)_{\cal G}\sim |Z(\eps,0)_{\cal G}|$, and also large enough so
that $Z(\eps,\alpha)_{\cal G}=Z(\eps,0)_{\cal G} + \tau_{\cal G}\alpha$ is positive and 
bounded below, e.g.
$\tau_{\cal G}\alpha\ge 2 \bar{\eps}_{\cal G}$.

\begin{Lemma}[exit probabilities]  \label{lem:exit-proba}
 Choose $\alpha$ such that $Z(\eps,\alpha)_{\cal G} \sim Z(-\eps,0)_{\cal G}$.  Let $f\equiv f(\sigma^*_{ext})$ be the exit 
probabilities, solution of  (\ref{eq:exit-proba-eq}). Then
 
\BEQ f^{\sigma}(\sigma^*_{ext}) \equiv  w(\alpha)_{G\to\sigma^*_{ext}} \ \times\      (1+O(Z(-\eps,0)_{{\cal G}}))
\EEQ
where by definition,
\BEQ w(\alpha)_{G\to\sigma^*_{ext}}  := \frac{ 
\sum_{\sigma} 
\tilde{\pi}^{int}_{\sigma} \frac{k_{\sigma\to \sigma^*_{ext}}}{k_{\sigma}}  }
 { Z(\eps,\alpha)_{{\cal G}} 
 }   \label{eq:lem:exit-proba}
\EEQ
In particular, 
\BEQ w(\alpha)_{G\to\sigma^*_{ext}}  \sim \frac{1}{Z(\eps,\alpha)_{{\cal G}}} \ \times \ \max_{\sigma} (\frac{k_{\sigma\to\sigma^*_{ext}}}{k_{\sigma}})
\EEQ
\end{Lemma}

{\bf Proof.}  Let $g= (g_{\sigma})_{\sigma\in {\cal V}^{int}}$ be the r.-h.s. of  (\ref{eq:exit-proba-sim}), so $f =  -(\tilde{\cal W} - \Id)^{-1}g$.
As proved above (see previous paragraph),
\BEQ ||P_{\perp}f||_{\infty} \preceq  ||P_{\perp} g||_{\infty} \preceq ||g||_{\infty}. \EEQ

\Medskip 
Let $w_{\to \sigma^*_{ext}} := \max_{\sigma} (\frac{k_{\sigma\to\sigma^*_{ext}}}{k_{\sigma}})$. Replacing
$g$ by its value, one gets
\BEA ||g||_{\infty} &\preceq  &
Z(-\eps,0)_{{\cal G}}\ 
 ||f||_{\infty} + w_{\to \sigma^*_{ext}} 
\nonumber\\
&\preceq & Z(-\eps,0)_{{\cal G}} \ \Big( ||P_{\perp}f||_{\infty} + |f_{av}|\Big) 
+ w_{\to \sigma^*_{ext}}
\EEA
from which
\BEQ ||P_{\perp}f ||_{\infty} \preceq   Z(-\eps,0)_{{\cal G}}\  |f_{av}| +  w_{\to \sigma^*_{ext}}. \label{eq:Pavfinfty}
\EEQ
We now rewrite
 identity (\ref{eq:exit-proba-sim}) as
\BEQ 
(1 - \eps_{\sigma} + \frac{k^{ext}_{\sigma}+\alpha}{k_{\sigma}})  f^{\sigma}= \sum_{\sigma'\in {\cal V}} \tilde{\cal W}^{int}_{\sigma,\sigma'} f^{\sigma'} + \frac{k_{\sigma\to\sigma^*_{ext}}}{k_{\sigma}+\alpha},  \label{eq:227}
\EEQ
 multiply both
sides by $\tilde{\pi}^{int}_{\sigma}$, and sum over $\sigma$. Because $\sum_{\sigma} \tilde{\pi}^{int}_{\sigma} \tilde{\cal W}_{\sigma,\sigma'} = \tilde{\pi}^{int}_{\sigma'}$, the first term in the r.-h.s. cancels
with the term with 1 in factor in the l.-h.s.
Since (by (\ref{eq:Pavfinfty})) $f^{\sigma} = (P_{\perp}f)^{\sigma} + f_{av} =  (1+O(Z(-\eps,0)_{{\cal G}})) f_{av} + 
 O (w_{\to \sigma^*_{ext}})$ 
is  constant to leading order, it comes out
\BEQ 
 \sum_{\sigma} \tilde{\pi}^{int}_{\sigma} (-\eps_{\sigma} + \frac{k^{ext}_{\sigma}+\alpha}{k_{\sigma}})   \Big\{ 
 (1+O(Z(\eps,\alpha)_{{\cal G}})  f_{av}  \ + O (w_{\to \sigma^*_{ext}}) \Big\}
  = \sum_{\sigma} \tilde{\pi}^{int}_{\sigma} \frac{k_{\sigma\to\sigma^*_{ext}}}{k_{\sigma} + \alpha} 
\EEQ

The main term on the l.-h.s. is $f_{av}$, times  $\sum_{\sigma} \tilde{\pi}^{int}_{\sigma} (-\eps_{\sigma} + \frac{k^{ext}_{\sigma}+\alpha}{k_{\sigma}})  = Z(\eps,\alpha)_{{\cal G}}$. 
Next, $w_{\to \sigma^*_{ext}} \sim \sum_{\sigma} \tilde{\pi}^{int}_{\sigma} \frac{k_{\sigma\to\sigma^*_{ext}}}{k_{\sigma}} $, whence $\sum_{\sigma} \tilde{\pi}_{\sigma} \frac{k_{\sigma\to\sigma^*_{ext}}}{k_{\sigma}+\alpha} 
\ -\   O (w_{\to \sigma^*_{ext}})  \sum_{\sigma} \tilde{\pi}^{int}_{\sigma} (-\eps_{\sigma} + \frac{k^{ext}_{\sigma} + \alpha^{int}}{k_{\sigma}} )
=  \sum_{\sigma} \tilde{\pi}^{int}_{\sigma} \frac{k_{\sigma\to\sigma^*_{ext}}}{k_{\sigma}}  (1+ O ( Z(-\eps,0)_{{\cal G}}))$, from which 
 one gets 
 \BEQ f_{av} =  \frac{ 
\sum_{\sigma} 
\tilde{\pi}_{\sigma} \frac{k_{\sigma\to \sigma^*_{ext}}}{k_{\sigma}}  }
 { Z(\eps,\alpha)_{{\cal G}}
 } \ \times\      (1+O(Z(-\eps,0)_{{\cal G}}))  \label{eq:fav-formula}.
\EEQ
 Using $|f^{\sigma} - f_{av}| = |(P_{\perp}f)^{\sigma}|  \le ||P_{\perp}f||_{\infty}$, and plugging the
 value (\ref{eq:fav-formula}) for $f_{av}$ into  the bound (\ref{eq:Pavfinfty}) yields finally  
 (\ref{eq:lem:exit-proba}). \hfill\eop


\subsubsection{Time-continuous adjoint problem}  \label{subsubsection:time-continuous-adjoint-pb}


We finally estimate the solution of the time-continuous adjoint problem introduced in  (\ref{eq:continuous-exit-proba-eq}).

\begin{Lemma}[time-continuous adjoint problem]  \label{lem:time-cont-adj-pb}
 Choose $\alpha$  such that $Z(\eps,\alpha)_{\cal G} \sim Z(-\eps,0)_{\cal G}$.   Let  $v\equiv v(\sigma^*_{ext})$ be the solution of 
the time-continuous adjoint problem (\ref{eq:continuous-exit-proba-eq}).
Then 
\BEQ v_{\sigma'}(\sigma^*_{ext})  \sim C(\sigma^*_{ext})\frac{\tilde{\pi}^{int}_{\sigma'}}{k_{\sigma'}} \ \times\      (1+O(Z(-\eps,0)_{{\cal G}}))
\EEQ
where by definition, 
\BEQ C(\sigma^*_{ext}):=   \frac{k_{\sigma^*_{ext}\to G}}{Z(\eps,\alpha)_{{\cal G}} } \label{eq:lem:time-cont-adj-pb}
  \EEQ
\end{Lemma} 

{\bf Proof.}
We use the identity $\tilde{A}^{int}_{\sigma',\sigma} = k_{\sigma\to\sigma'} - k^{int}_{\sigma} \del_{\sigma,\sigma'} =   (\tilde{w}^{int}_{\sigma\to\sigma'} - \del_{\sigma,\sigma'})\, k^{int}_{\sigma} = 
\, ((^t \tilde{\cal W}^{int}- \Id) K^{int})_{\sigma',\sigma}$, where $K^{int} = \diag(k^{int}_{\sigma})$ is diagonal, and 
combine arguments found in the proofs of Corollary \ref{cor:eta} and Lemma \ref{lem:exit-proba}. 
First, the equations $A(\alpha)v=0$ inside ${\cal V}^{int}$ can be rewritten as (see
(\ref{eq:tcap}))
\BEQ \tilde{A}^{int}v = u; \qquad u_{\sigma'} = (k^{ext}_{\sigma'} + \alpha - \kappa_{\sigma'}) v_{\sigma'} - k_{\sigma^*_{ext}\to \sigma'}, \ \sigma'\in {\cal V}^{int} \EEQ
or equivalently, $(\, ^t\tilde{\cal W}^{int} - \Id) K^{int} v = u$. The projection matrix $P^{\perp}: u \mapsto 
u- |{\cal V}^{int}| u^{av}\tilde{\pi}^{int}$ commutes with $\, ^t \tilde{\cal W}^{int}$. Thus
\BEQ (\, ^t\tilde{\cal W}^{int} - \Id) P^{\perp} (K^{int} v) = P^{\perp} u \EEQ
from which $||P^{\perp} (K^{int} v)||_1 \preceq  ||P^{\perp} u||_1$. Thus,
\BEQ k_{\sigma'}^{int} v_{\sigma'}= (K^{int}v)_{\sigma'} =  C \tilde{\pi}^{int}_{\sigma'}+ O(||P^{\perp}u||_1) \EEQ
with $C\equiv C(\sigma^*_{ext}):= (K^{int}v)^{av}$. 
Now, 
\BEA ||P^{\perp}u||_1 \preceq ||u||_1 &\le& \sum_{\sigma'} (\frac{k^{ext}_{\sigma'} + \alpha}{k_{\sigma'}} + \eps_{\sigma'}) 
k_{\sigma'} v_{\sigma'} + k_{\sigma^*_{ext}\to G}  \nonumber\\
& \preceq & Z(-\eps,0)_{{\cal G}}
\sum_{\sigma'} k_{\sigma'} v_{\sigma'} + k_{\sigma^*_{ext}\to G} \nonumber\\
& \preceq & Z(-\eps,0)_{{\cal G}} \ C
 + k_{\sigma^*_{ext}\to G}   \label{eq:Pavu-bound}
\EEA
from which we conclude that
\BEQ k^{int}_{\sigma'}v_{\sigma'} = C\tilde{\pi}^{int}_{\sigma'} \, \times \Big(1+O(Z(-\eps,0)_{{\cal G}})\Big) + O( k_{\sigma^*_{ext}\to G})
\EEQ
The same identity holds for $k_{\sigma'}v_{\sigma'}$ since $k_{\sigma'}/k^{int}_{\sigma'} = 1+O(Z(-\eps,0)_{{\cal G}})$. 

\Medskip To identify
the coefficient $C(\sigma^*_{ext})$, we rewrite the equations as
\BEQ \sum_{\sigma \in {\cal V}^{int}} \tilde{A}^{int}_{\sigma',\sigma} v_{\sigma} = 
(k^{ext}_{\sigma'} + \alpha - \kappa_{\sigma'}) v_{\sigma'} - k_{\sigma^*_{ext}\to \sigma'},
\EEQ
and sum over $\sigma'\in {\cal V}^{int}$ to cancel the l.-h.s. We obtain
\BEQ C\sum_{\sigma'} (\frac{k^{ext}_{\sigma'}+\alpha}{k_{\sigma'}}-\eps_{\sigma'}) \tilde{\pi}^{int}_{\sigma'} \, \times \Big(1+O(Z(-\eps,0)_{{\cal G}})\Big)    =  k_{\sigma^*_{ext}\to G}\, \Big(1+O(Z(-\eps,0)_{{\cal G}})\Big)
\EEQ
The main term on the l.-h.s. is $C \, \times\,  Z(\eps,\alpha)_{{\cal G}}$. 
Hence
\BEQ C \sim  \frac{k_{\sigma^*_{ext}\to G}}{Z(\eps,\alpha)_{{\cal G}}} \, \times\,  (1+ O ( Z(-\eps,0)_{{\cal G}}))  \label{eq:C-formula}
\EEQ

 Using $|k^{int}_{\sigma} v_{\sigma} - C |{\cal V}^{int}|\tilde{\pi}_{\sigma}^{int}| = |(P^{\perp}K^{int}v)_{\sigma}|  \le ||P^{\perp}(K^{int}v)||_1 
 \preceq  ||P^{\perp}u||_1$, and plugging the
 value (\ref{eq:C-formula}) for $C$ into  the bound (\ref{eq:Pavu-bound}) yields finally  
 (\ref{eq:lem:time-cont-adj-pb}).  \hfill\eop


\subsection{Renormalization}  \label{subsection:renormalization}


We change  notations and scope compared to the previous two subsections. The starting point is 
an "underlined"  graph $\underline{G}=(\underline{V},\underline{E})$. The splitting $\underline{V} = 
{\cal V}^{int} \uplus {\cal V}^{ext}$ defines a cut-off graph ${\cal G}^{int}$, extended into 
${\cal G}$, with the same properties as in \S \ref{subsection:cut-off-setting}.  The adjoint generator $\underline{A}$ of $\underline{G}$ has a block structure
\BEQ \underline{A} = \left(\begin{array}{cc} \underline{A}^{int} &  \underline{A}^{ext\to int}
 \\ \underline{A}^{int\to ext} 
& \underline{A}^{ext} \end{array}\right).   \label{eq:underlineA}
\EEQ
In general,  ${\cal G}^{int} = \uplus_{p=1,\ldots,q} {\cal G}^{int}_p$ is a disjoint union of
connected subgraphs with vertex subset ${\cal V}^{int}_p$ and deficiency weights $\eps_p$. 
The general purpose is to prove that the Lyapunov data of $\underline{A}$ are to leading order
equivalent to those of a renormalized  matrix $\underline{A}^{ren}$, which is the adjoint 
generator of a defective Markov chain on the merged vertex set $\{G_p,p=1,\ldots,q\}\uplus 
{\cal V}^{ext}$, see \S \ref{subsection:merging}.


\subsubsection{An introductory computation}  \label{subsubsection:an-introductory-computation}


\begin{Definition}[renormalized generator $\underline{A}^{ren}$] \label{def:Aren}
Let  $ \underline{A}^{ren}(\alpha) = \underline{A}^{ren}-\alpha \Id$  be the  $(q, |{\cal V}^{ext} |)$-block matrix with coefficients
\BEQ  \underline{A}^{ren}(\alpha)_{\sigma',\sigma} \equiv 
\left(\begin{array}{cc}    \diag(\underline{A}^{ren}(\alpha)_{G_p,G_p})_p &  \underline{A}^{ren}(\alpha)_{G_p,\sigma} \\  \underline{A}^{ren}(\alpha)_{\sigma',G_p} &  \underline{A}^{ren}(\alpha)_{\sigma',\sigma}
\end{array}\right)
= \left(\begin{array}{cc}  
- \diag(\frac{Z(\eps_p,\alpha)_{{\cal G}_p}}{\tau_{{\cal G}_p}})_p & k_{\sigma\to G_p}
 \\ k_{G_p\to \sigma'} 
& A^{ext}(\alpha)_{\sigma',\sigma} \end{array}\right),   \label{eq:A-ren}
\EEQ
where $A^{ext}(\alpha)=A^{ext}-\alpha \Id$ and :

\Medskip \textbullet\ for every $p=1,\ldots,q$,\   $\tau_{{\cal G}_p} = \sum_{\sigma\in {\cal V}^{int}_p} \frac{\tilde{\pi}^{int}_{p,\sigma}}{k_{\sigma}}, Z(\eps_p,\alpha)_{{\cal G}^{int}_p} = \sum_{\sigma\in {\cal V}^{int}} \tilde{\pi}^{int}_{p,\sigma} \Big(
\frac{k^{ext}_{\sigma} + \alpha}{k_{\sigma}} - \eps_{p,\sigma}\Big)$ are as in Lemma \ref{lem:Lya(1)}, and   
\BEQ k_{\sigma\to G_p} = \sum_{\sigma' \in {\cal V}^{int}_p} k_{\sigma\to \sigma'}, \qquad 
\sigma\in 
{\cal V}^{ext}\EEQ
as in the introduction to \S \ref{subsection:cut-off-setting}); 

\Medskip \textbullet\ (equivalent transition rates from $G$) as in the
heuristic \S \ref{subsection:heuristics}, we define
\BEQ k_{G_p \to\sigma'} := \frac{1}{\tau_{{\cal G}_p}} \sum_{\sigma\in {\cal V}^{int}_p} \tilde{\pi}^{int}_{p,\sigma}
\frac{k_{\sigma\to\sigma'}}{k_{\sigma}}, \qquad  
\sigma' \in 
{\cal V}^{ext} 
 \EEQ

Then $\underline{k}^{ren}_{\sigma\to\sigma'}$, resp. $\underline{k}^{ren}_{\sigma}$, are the transition rates, resp. outgoing rates, for $\underline{A}^{ren}$, namely,
\BEQ \underline{k}^{ren}_{\sigma\to\sigma'} = \underline{A}^{ren}_{\sigma',\sigma} \qquad 
(\sigma'\not=\sigma), \qquad  \underline{k}^{ren}_{\sigma} := \sum_{\sigma'\not=\sigma} 
\underline{k}^{ren}_{\sigma\to\sigma'}
\EEQ
while $\underline{k}^{ren}_{min}:= \min_{\sigma\not=\sigma'} \underline{k}^{ren}_{\sigma\to\sigma'}$
is a characteristic (minimum scale transition) rate of $\underline{A}^{ren}$, see (\ref{eq:char-rate}). 

\Medskip Finally, let  $\underline{\tilde{A}}^{ren}$ be the matrix with corrected diagonal coefficients,
\BEQ \underline{\tilde{A}}^{ren}_{\sigma',\sigma}= \underline{A}^{ren}_{\sigma',\sigma}
 \qquad (\sigma\not=\sigma'), \qquad   \underline{\tilde{A}}^{ren}_{\sigma,\sigma}= - 
 \sum_{\sigma'}  \underline{A}^{ren}_{\sigma',\sigma} 
\EEQ
\end{Definition}

\noindent The sum of the coefficients of each  column of $\underline{A}^{ren}(\alpha)$ 
vanishes when $\eps=\alpha=0$; this is obvious except for the merged vertices, for which the sum
\BEA && \sum_{\sigma'} \underline{A}^{ren}(\alpha)_{\sigma',G_p} = - \frac{Z(\eps,\alpha)_{{\cal G}_p}}{\tau_{{\cal G}_p}} + \sum_{\sigma'\in  {\cal V}^{ext}}  k_{G_p
\to \sigma'} \nonumber\\
&& = -\alpha +  \frac{1}{\tau_{{\cal G}_p}} \Big(\bar{\eps}_{{\cal G}_p}-Z(0)_{{\cal G}_p} +  
\sum_{\sigma\in {\cal V}^{int}} \tilde{\pi}_{p,\sigma}^{int} \  \frac{\sum_{\sigma'\in {\cal V}^{ext}} k_{\sigma\to\sigma'}}{k_{\sigma}} \Big) =    -\alpha +
\frac{\bar{\eps}_{{\cal G}_p}}{\tau_{{\cal G}_p}}  \nonumber\\
\label{eq:sumcoeff1stcolAren=0}
\EEA
vanishes if $\eps_p=0$ and $\alpha=0$.


\Bigskip {\em Time-continuous adjoint problem} (see \S \ref{subsection:cut-off-setting}). We
simply motivate heuristically the coefficients of $\underline{A}^{ren}$ by looking at the
time-continuous adjoint problem. For simplicity, we assume $q=1$ and write simply 
$G$ instead of $G_1$.  
Let $\underline{U}\in\R^{\underline{V}}$, and $\underline{U}=\left(\begin{array}{c} 
U^{int} \\ U^{ext} \end{array}\right)$ its decomposition along ${\cal V}^{int} \uplus {\cal V}^{ext} $.   We solve the system 
\BEQ \underline{A}(\alpha) \underline{U} =0   \label{eq:mathbbAintU=0}
\EEQ
on ${\cal V}^{int}$ with boundary condition $U^{ext}$ on ${\cal V}^{ext}$, with an appropriate
choice of $\alpha$ so that Lemma \ref{lem:time-cont-adj-pb} holds.
 Then the solution is (by the superposition principle)
\BEQ U^{int}_{\sigma} \sim {\cal C}(U^{ext}) \frac{\tilde{\pi}^{int}_{\sigma}}{k_{\sigma}}, \EEQ
where:  
\BEQ {\cal C}(U^{ext}):=\sum_{\sigma^*_{ext}\in {\cal V}^{ext}}
 U^{ext}_{\sigma^*_{ext}} C(\sigma^*_{ext}), \qquad
  C(\sigma^*_{ext}) = \frac{k_{\sigma^*_{ext}\to G}}{Z(\eps,\alpha)_{{\cal G}}}.
  \label{eq:CUext}
 \EEQ
Equivalently,
\BEQ  \sum_{\sigma \in 
{\cal V}^{int}} k_{\sigma\to\sigma'} U^{int}_{\sigma} \sim  k_{G\to\sigma'} \tau_{{\cal G}} \  {\cal C}(U^{ext}).  \label{eq:ktauUVkU} 
\EEQ
 
\Medskip  Introduce the scalar
\BEQ u_{G}:=  \tau_{{\cal G}}\  {\cal C}(U^{ext}) 
 \EEQ
Then (\ref{eq:mathbbAintU=0}) 
is equivalent to leading order to the system 
\BEQ \underline{A}^{ren}(\alpha) \underline{U}^{ren} =0, \qquad \underline{U}^{ren} := \left(\begin{array}{c} u_{G} \\ U^{ext} \end{array} \right).
\EEQ
 Namely (see (\ref{eq:A-ren})),  looking at the $G$-index line, 
 \BEQ \sum_{\sigma^*_{ext} \in {\cal V}^{ext}} k_{\sigma^*_{ext}\to G} \ 
\underline{U}^{ren}_{\sigma^*_{ext}} = \sum_{\sigma^*_{ext} \in {\cal V}^{ext}} k_{\sigma^*_{ext}\to G} \  U^{ext}_{\sigma^*_{ext}}   = Z(\eps,\alpha)_{{\cal G}} \ 
{\cal C}(U^{ext})    = - \underline{A}^{ren}_{G,G} \, u_{G} \EEQ
by (\ref{eq:CUext}). Looking now at the $\sigma^*_{ext}$-line $(\sigma^*_{ext}\in {\cal V}^{ext})$,  $k_{G\to \sigma^*_{ext}} u_G \sim \sum_{\sigma \in 
{\cal V}^{int}} k_{\sigma\to \sigma^*_{ext}} U^{int}_{\sigma}$ (see (\ref{eq:ktauUVkU})),  which adds up
to $\sum_{\sigma\in {\cal V}^{ext}} A^{ext}(\alpha)_{\sigma^*_{ext},\sigma} U^{ext}_{\sigma}$ to give back $\sum_{\sigma\in \underline{V}} \underline{A}_{\sigma^*_{ext},\sigma} (\alpha)
U_{\sigma}=0$.



\subsubsection{Table of regimes}  \label{subsubsection:table-of-regimes}


\noindent Let $G$ be one of the $G_p$.

\Bigskip {\bf Renormalized weight matrix and stationary measure.} Let 
\BEQ \underline{\cal W}^{ren}(\alpha) \equiv   \left(\begin{array}{cc}  0 & 
\underline{w}^{ren}(\alpha)_{G\to\sigma'} \\ \underline{w}^{ren}(\alpha)_{\sigma\to G} & 
\underline{w}^{ren}(\alpha)_{\sigma\to\sigma'} \end{array}\right) =
\left(\begin{array}{cc}  0 &  \frac{k_{G\to\sigma'}}{
Z(\eps,\alpha)_{{\cal G}}/\tau_{{\cal G}}} \\
 \frac{k_{\sigma\to G}}
{|\underline{A}_{\sigma,\sigma}|
 + \alpha}  & \frac{k_{\sigma\to\sigma'}}{|A_{\sigma,\sigma}|+\alpha} \end{array}\right)  
\EEQ
 be the weight
matrix associated to $\underline{A}^{ren}(\alpha)$. In particular, let $\underline{\lambda}^{ren}$ be the Lyapunov eigenvalue of $\underline{A}^{ren}$, and $\underline{\pi}^{ren}:=\underline{\pi}^{ren}(\underline{\lambda}^{ren})$ the
left eigenvector of $\underline{\cal W}^{ren}(\underline{\lambda}^{ren})$ with eigenvalue 1, normalized
by letting $\underline{\pi}^{ren}_G + \sum_{\sigma\in {\cal V}_{ext}} \underline{\pi}^{ren}_{\sigma}=1$.  

\Medskip Similarly, let 
\BEQ \underline{\tilde{\cal W}}^{ren} \equiv  
 \left(\begin{array}{cc}  0 & 
\underline{\tilde{w}}^{ren}_{G\to\sigma'} \\ \underline{\tilde{w}}^{ren}_{\sigma\to G} & 
\underline{\tilde{w}}^{ren}_{\sigma\to\sigma'} \end{array}\right) =
\left(\begin{array}{cc}  0 &  \frac{k_{G\to\sigma'}}{
Z(\eps,\tilde{\alpha})_{{\cal G}}/\tau_{{\cal G}} + \tilde{\eps}_{{\cal G}}} \\
 \frac{k_{\sigma\to G}}
{k_{\sigma}}  & \frac{k_{\sigma\to\sigma'}}{k_{\sigma}}  \end{array}\right)  
\label{eq:underlineWtilderen}
\EEQ
 be the weight
matrix associated to $\underline{\tilde{A}}^{ren}$. If $\underline{A}$ is irreducible, then 
the Markov chain with matrix  $\underline{\tilde{\cal W}}^{ren}$ has a unique stationary probability measure,
which we denote $\underline{\tilde{\pi}}^{ren}$.

\Medskip  By construction, the following identities hold,
 \BEQ \underline{w}^{ren}(\alpha)_{\sigma\to\sigma'} = w(\alpha)_{\sigma\to\sigma'} \ \  (\sigma,\sigma'\in {\cal V}^{ext}),  \qquad \underline{w}^{ren}(\alpha)_{\sigma\to G} = \sum_{\sigma'\in {\cal V}^{int}} w(\alpha)_{\sigma\to\sigma'} \ \  (\sigma\in {\cal V}^{ext})
 \label{eq:wrenext}
 \EEQ
  and similarly for $\underline{\tilde{w}}^{ren}$ instead of 
$\underline{w}^{ren}(\alpha)$.

\Medskip Thus the table of regimes introduced above in the case of a cycle (see \S \ref{subsection:cycle-elementary})
may be generalized with hardly any modifications as follows:

\Bigskip 
\begin{center}
\begin{tabular}{c|c|c|c} \label{table:regimes2}
Regime & free &  autocatalytic & degraded \\ \hline
Conditions & $Z(0)_{{\cal G}}\succ \alpha \tau_{{\cal G}},\bar{\eps}_{{\cal G}}$ & 
$\alpha \tau_{{\cal G}} \sim \bar{\eps}_{{\cal G}} \succ Z(0)_{{\cal G}}$ & 
$\alpha\tau_{{\cal G}} \succ Z(0)_{{\cal G}}, \bar{\eps}_{{\cal G}} $ \\
   \hline
Value of $Z(\eps,\alpha)_{{\cal G}}$ & $Z(0)_{{\cal G}}$ & $\alpha \tau_{{\cal G}}$ & $\alpha\tau_{{\cal G}}$ \\
Value of $\underline{w}^{ren}(\alpha)_{G\to \sigma_{ext}}$ & $Z(0)_{{\cal G}}^{-1} 
\langle \tilde{w}_{\cdot\to \sigma_{ext}}\rangle_{int}
 $ & $\frac{1}{\alpha\tau_{{\cal G}}} \langle \tilde{w}_{\cdot\to \sigma_{ext}}\rangle_{int}$
& $\frac{1}{\alpha\tau_{ {\cal G} }} \langle \tilde{w}_{\cdot\to \sigma_{ext}}\rangle_{int}$ \\ \hline
Value of $k_{G\to \sigma_{ext}}$ & $\frac{1}{\tau_{{\cal G}}} 
\langle \tilde{w}_{\cdot\to \sigma_{ext}}\rangle_{int}$ & 
$\frac{1}{\tau_{{\cal G}}} 
\langle \tilde{w}_{\cdot\to \sigma_{ext}}\rangle_{int}$
  & 
$\frac{1}{\tau_{{\cal G}}} 
\langle \tilde{w}_{\cdot\to \sigma_{ext}}\rangle_{int}$ 
\end{tabular}

\bigskip {\textsc{TABLE \ref{table:regimes} -- The three different regimes for a general
non-trivial maximal dominant SCC.}}
\end{center}
where:  "value of" means
"same order as" ($"\sim"$), and $\langle \tilde{w}_{\cdot\to \sigma_{ext}}\rangle_{int} := 
\sum_{\sigma} 
\tilde{\pi}^{int}_{\sigma} \tilde{w}_{\sigma\to \sigma_{ext}}
$ as in \S \ref{subsubsection:exit-proba}. Let $\alpha_{thr} := k_{{\cal G}}^{ext} \sim 
\frac{1}{\tau_{{\cal G}}} \, \times \, Z(0)_{{\cal G}}$  (if $G$ is not autocatalytic),
$\alpha_{thr} \sim \lambda_{{\cal G}} \sim \frac{1}{\tau_{{\cal G}}} \, \times\, \bar{\eps}_{{\cal G}}$ 
(if $G$ is autocatalytic), see discussion in \S \ref{subsection:description}. 
The three different Conditions in the Table cover all values of  $\alpha$, under the
restriction $\alpha\succeq \alpha_{thr}$ in the autocatalytic case. 
The autocatalytic Regime is obtained in the autocatalytic case for $\alpha$ such that $\alpha\sim\alpha_{thr}$. On the other hand, choosing $\alpha\succ \alpha_{thr}$, 
 leads
to the degraded Regime.


\subsubsection{One renormalization step}  \label{subsubsection:one-renormalization-step}


All computations in \S \ref{subsubsection:an-introductory-computation} and \S \ref{subsubsection:table-of-regimes} generalize straightforwardly to the case when ${\cal G}$ is a disjoint union $\uplus_p {\cal G}_p$ 
of subgraphs ${\cal G}_p$ with internal vertices ${\cal V}^{int}_p$ and deficiency weights $\eps_p$. 
{\em We assume that none of the ${\cal G}_p$ is autocatalytic.}
We shall use the obvious notations ${\cal V}^{int} = \uplus_p {\cal V}^{int}_p$;  
$({\bf 1}_{int,p})_{\sigma}=1$ if and only if $\sigma\in {\cal V}_p^{int}$; ${\bf 1}_{int} = 
\sum_p {\bf 1}_{int,p}$;  $\langle f,u\rangle_{int,p} := \sum_{\sigma \in {\cal V}_{int,p}} 
f_{\sigma} u_{\sigma}$;   $\langle f,u\rangle_{int} := \sum_{\sigma \in {\cal V}_{int}} 
f_{\sigma} u_{\sigma}$.  

\Medskip We rewrite the eigenvalue equations $\underline{A} \underline{v}^* = \underline{\lambda}^* \underline{v}^*, \  ^t \underline{A}\, 
\underline{v}^{*,\dagger} = \underline{\lambda}^* \underline{v}^{*,\dagger} $ in the equivalent form 
\BEQ \underline{\pi}^* \underline{\cal W}(\underline{\lambda}^*) = \underline{\pi}^*, \qquad
 \underline{\cal W}(\underline{\lambda}^*) \underline{v}^{*,\dagger} = \underline{v}^{*,\dagger}
\label{eq:underlinepi*}
\EEQ
 where
$\underline{w}(\underline{\lambda}^*)_{\sigma\to\sigma'} = \frac{k_{\sigma\to\sigma'}}{|\underline{A}_{\sigma,\sigma}|+\bar{\lambda}^*}$, $\underline{\pi}^*_{\sigma} = 
(|\underline{A}_{\sigma,\sigma}|+\underline{\lambda}^*) \underline{v}^*_{\sigma}$. 
{\em We assume here that $\underline{\pi}^{ren}_{G_p}\sim 1$ for all $p$, and that 
$\underline{\lambda}^{ren}\prec \underline{k}^{ren}_{min}$ (see Definition 
\ref{def:Aren}).} In the next subsection, this will hold as a consequence of an induction.   

\Medskip 
 Let first  $d_{\sigma}(\alpha) := \frac{k^{ext}_{\sigma} + \alpha}{k_{\sigma}} - \eps_{\sigma}$, $\sigma\in {\cal V}^{int}$, and $D(\alpha) = \diag(d_{\sigma}(\alpha))$, so that $D(\alpha) {\bf 1}_{int}=d(\alpha)$. Finally, let 
 $\Del(\alpha):=(Z(\eps,\alpha))^{-1}D(\alpha)$, short-hand for $\Del(\alpha) = (\Del_p(\alpha))_p, \ \Del_p(\alpha) = Z^{-1}(\eps_p,\alpha)D_p(\alpha)$, where $D_p := D\Big|_{{\cal V}^{int}_p}$.  By definition, $\langle \tilde{\pi}^{int}_p, d(\alpha)\rangle_{int,p} = Z(\eps_p,\alpha)_{{\cal G}_p}$. 
 Note that \BEQ \langle \tilde{\pi}^{int}_{p,\cdot}, w(\alpha)_{\cdot\to \sigma'} \rangle_{int,p}  \sim \tau_{{\cal G}_p^{int}} k_{G_p\to \sigma'} \sim Z(\eps_p,\alpha)_{{\cal G}_p} \,  \underline{w}^{ren}(\alpha)_{G_p\to\sigma'}, \qquad \sigma'\in {\cal V}^{ext}.
 \EEQ 
Also,  the following two identities hold,
\BEQ  (\underline{\cal W}^{int\to ext}(0) {\bf 1}_{ext})_{\sigma} = \frac{k^{ext}_{\sigma}}{k_{\sigma}} 
= d_{\sigma}(0) + \eps_{\sigma}, \qquad \sigma \in {\cal V}^{int}  \label{eq:id-intext}
\EEQ
\BEQ \frac{d}{d\alpha}Z(\eps_p,\alpha)_{{\cal G}_p} = \langle \tilde{\pi}^{int}_p, 
 \frac{1}{k_{\cdot}}\rangle_{int,p} \sim \tau_{{\cal G}_p}, \label{eq:id-dZ/dalpha}
\EEQ
 
\Medskip There are two small parameters here; the first one, relative to the
subgraphs $G_p$, is 
$Z(0)_{{\cal G}_p}  = \langle   \tilde{\pi}^{int}_p, \frac{k^{ext}_{\cdot}}{k_{\cdot}} \rangle_{int,p}$, which is $>0$ and $\prec 1$. Since we have assumed that none of the ${\cal G}_p$ is autocatalytic,  
$Z(\eps_p,0)\sim Z(0)_{{\cal G}_p}$ is $>0$ too. We let $\alpha$ be a free parameter in the following, and later on expand to order 1 around $0$. Then all components
of $d(\alpha)|_{{\cal V}^{int}_p}$, $D(\alpha)|_{{\cal V}^{int}_p}$ and  $\underline{{\cal W}}^{int\to ext}(\alpha)|_{{\cal V}^{int}_p}$, resp. $\Del(\alpha)$, are of order $Z(0)_{{\cal G}_p}$, resp. of order $1$. By abuse of notation, we write simply: 
$d(\alpha),D(\alpha),\underline{\cal W}^{int\to ext}(\alpha)=O(Z(0)_{\cal G})$; this kind of block  bounds will
appear several times later on. The second samall parameter, on the other hand, is relative to the renormalized graph,
\BEQ \underline{Z}^{ren}:= \underline{\lambda}^{ren} / \underline{k}^{ren}_{min} 
\prec 1 \label{eq:Zren} 
\EEQ

\Medskip {\em Step 1.}  Decompose a vector $f= \left(\begin{array}{c} f_{int} = (f_{int,p})_p \\ f_{ext} 
\end{array}\right)$ 
into the sum $f =  f_{\perp} + f_{av}  +  f_{ext}$, where  $f_{av}\equiv  \sum_p f_{av,p}  {\bf 1}_{int,p}$,  $f_{av,p} \equiv P_{av,p}f:=   \langle \tilde{\pi}^{int}_p, f_{int} \rangle_{int,p}$ and   $f_{\perp,p} \equiv  P_{\perp,p} f_{int,p}
:= f_{int,p} - f_{av,p} {\bf 1}_{int,p}$. If $p=1$,  any linear map $\cal M$ with  matrix $M = (M_{\sigma,\sigma'})_{\sigma,\sigma'\in 
{\cal V}}$ in the standard species basis  may be rewritten
as a 3-block matrix $\left[\begin{array}{ccc} M_{\perp,\perp} & M_{\perp,av} & 
M_{\perp,ext} \\ M_{av,\perp} & M_{av,av} & M_{av,ext} \\ M_{ext,\perp} & 
M_{ext, av} & M_{ext,ext} \end{array}\right]$  by projecting vectors along this decomposition; generalization
to $p>1$ is straightforward. Dualizing, a measure $\pi = \left(\begin{array}{c} \pi^{int} \\
\pi^{ext} \end{array}\right)$ decomposes into a sum $\pi^{\perp} + \pi^{av} + \pi^{ext}$ using
the projections $^t P_{av} = P^{av}$ and $^t P_{\perp} = P^{\perp}$. 
The above averaged identities for $d$ and $w(\alpha)_{\cdot\to\sigma'}$ translate into
\BEQ P_{av,p}d_p = P_{av,p} D_p {\bf 1}_{int,p} = Z(\eps_p,\alpha)_{{\cal G}_p};\ \  (P_{av,p} \underline{\cal W}^{int\to ext}(\alpha))_{\sigma'} \sim Z(\eps_p,\alpha)_{{\cal G}_p} \, \underline{w}^{ren}(\alpha)_{G_p\to\sigma'}, \sigma'\in 
{\cal V}^{ext}  \label{eq:PavDw}
\EEQ 
Let us project ${\cal M}= \underline{\cal W}(\alpha)-\Id$ along this new system of 
coordinates.

\Medskip  Looking first at the internal components of  $(\underline{\cal W}(\alpha)-\Id)f$, we get from (\ref{eq:exit-proba-sim1})
\BEQ ((\underline{\cal W}(\alpha)-\Id)f)^{\sigma} \sim  ((\tilde{\cal W}^{int}-\Id)f_{int})^{\sigma} 
- (D(\alpha) \tilde{\cal W}^{int} f_{int})^{\sigma}  + \langle w(\alpha)_{\sigma\to \cdot}, f_{ext}^{\cdot} \rangle_{ext}    , \qquad\sigma
\in {\cal V}^{int}
\EEQ
Since $P_{\perp}$, $P_{av}$ commute with $\tilde{\cal W}^{int}-\Id$, 
\BEQ  ((\underline{\cal W}(\alpha)-\Id)f)_{\perp} =  (\tilde{\cal W}^{int}-\Id)f_{\perp}
- P_{\perp} (D(\alpha) \tilde{\cal W}^{int} f_{\perp} +  f_{av} d(\alpha)) + P_{\perp} \underline{\cal W}^{int\to ext}(\alpha) f_{ext} 
\EEQ 
and (by (\ref{eq:PavDw}))
\BEA
 ((\underline{\cal W}(\alpha)-\Id)f)_{av} &=&  -P_{av} D(\alpha) \tilde{\cal W}^{int} f_{\perp} + Z(\eps,\alpha)_{\cal G} \, 
 (-f_{av} + \langle \underline{w}^{ren}(\alpha)_{G\to\cdot}, f_{ext}^{\cdot} \rangle_{ext})
 \nonumber\\
 &=&  -P_{av} D(\alpha) \tilde{\cal W}^{int}  f_{\perp} + Z(\eps,\alpha)_{\cal G} \ P_{av}  (\underline{\cal W}^{ren}(\alpha)-\Id) 
 \left(\begin{array}{c} f_{av} \\ f_{ext} \end{array}\right) \nonumber\\
 &=& Z(\eps,\alpha)_{\cal G} \ P_{av}  \ \Big(-\Del(\alpha) \tilde{\cal W}^{int} f_{\perp} + \underline{\cal W}^{ren}(\alpha)-\Id) 
 \left(\begin{array}{c} f_{av} \\ f_{ext} \end{array}\right) \Big) \nonumber\\
 \EEA 
 
\Medskip Let us now look at the external components: 
\BEQ   ((\underline{\cal W}(\alpha)-\Id)f)^{\sigma} =   ((\underline{\cal W}^{ext}(\alpha)-\Id)f_{ext})^{\sigma} 
+ \langle w(\alpha)_{\sigma\to \cdot}, f_{int}^{\cdot} \rangle_{int}    , \qquad\sigma
\in {\cal V}^{ext}
\EEQ
Further, for every $\sigma\in {\cal V}^{ext}$,   $((\underline{\cal W}^{ext}(\alpha)-\Id)f_{ext})^{\sigma} =   ((\underline{\cal W}^{ren}(\alpha)-\Id)f_{ext})^{\sigma}$ and   
$ \langle w(\alpha)_{\sigma\to \cdot}, f_{int}^{\cdot} \rangle_{int} =
 \underline{w}^{ren}(\alpha)_{\sigma\to G} f_{av} +  \langle w(\alpha)_{\sigma\to \cdot}, f_{\perp}^{\cdot} \rangle_{int}$ from (\ref{eq:wrenext}),  hence, summing:
 \BEQ  ((\underline{\cal W}(\alpha)-\Id)f)^{\sigma} =  ({\cal W}^{ext,\perp}(\alpha)f_{\perp})^{\sigma} +  ((\underline{\cal W}^{ren}(\alpha)-\Id)f)^{\sigma} , \qquad\sigma
\in {\cal V}^{ext}
\EEQ
where ${\cal W}^{ext,\perp}(\alpha)$ is the restriction of ${\cal W}^{ext\to int}(\alpha)$ to vectors $f$ of the form
$f=f_{\perp}$. 
All together, we have proved:
\BEQ \underline{\cal W}(\alpha)-\Id =  \left[\begin{array}{ccc} \Id & & \\ & Z(\eps,\alpha)_{\cal G} &  \\  & & \Id 
\end{array}\right] \ \left(\begin{array}{cc} P_{\perp} (\tilde{\cal W}^{int} - \Id) P_{\perp} &  0 \\ 
\left(\begin{array}{c} -P_{av}\Del(\alpha) \tilde{\cal W}^{int} P_{\perp} \\ 
{\cal W}^{ext,\perp}(\alpha) \end{array}\right) &  
\underline{\cal W}^{ren}(\alpha)-\Id  \end{array}\right) +    \underline{\cal W}_{(1)}(\alpha),
 \label{eq:WalphaZPperpW1} \EEQ
where
\BEQ   \underline{\cal W}_{(1)}(\alpha):= 
\left[\begin{array}{ccc} -P_{\perp}D(\alpha) \tilde{\cal W}^{int} P_{\perp}  & -P_{\perp}D(\alpha) P_{av} & P_{\perp} \underline{\cal W}^{int\to ext}(\alpha) \\
0 & 0 & 0  \\ 0 & 0 & 0 
\end{array} \right]   
\EEQ 
Since $\underline{\cal W}_{(1)}(\alpha)$ has non-zero entries only along the $\perp$-component, 
the left multiplying matrix in square brackets in (\ref{eq:WalphaZPperpW1}) can be put in factor 
of $\underline{\cal W}_{(1)}(\alpha)$ too. 
The second matrix in (\ref{eq:WalphaZPperpW1}) is block diagonal, with first block along the $\perp$ component, and second
block along the $(av,ext)$-components. All coefficients of $ \underline{\cal W}_{(1)}(\alpha)$ are of 
order $O(Z(0)_{\cal G})$ at most.


\Bigskip  {\em Step 2.}  Rewrite $\underline{\cal W}(\alpha)-\Id$ as  
$ {\cal Z}(\alpha) (  {\cal H}_{(0)} + {\cal H}_{(1)}(\alpha))$,
 where 
\BEQ {\cal Z}(\alpha):= \left[\begin{array}{ccc} \Id & & \\ & \diag(Z(\eps_p,\alpha)_{{\cal G}_p})_p &  \\  & & \Id 
\end{array}\right] 
\EEQ
is the {\bf renormalization matrix}; and
\BEA &&  
  {\cal H}_{(0)} :=  \left(\begin{array}{cc} P_{\perp} (\tilde{\cal W}^{int} - \Id) P_{\perp} & 0 \\
\left(\begin{array}{c} -P_{av}\Del(0) \tilde{\cal W}^{int} P_{\perp} \\ 
{\cal W}^{ext,\perp}(0) \end{array}\right) 
   &  
\underline{{\cal W}}^{ren}(\underline{\lambda}^{ren}) - \Id  \end{array}\right), \nonumber\\
&& \qquad
{\cal H}_{(1)} =    \left(\begin{array}{cc} 0 &  0 \\
\left(\begin{array}{c} -P_{av}(\Del(\alpha)-\Del(0)) \tilde{\cal W}^{int} P_{\perp} \\ 
{\cal W}^{ext,\perp}(\alpha)- {\cal W}^{ext,\perp}(0)\end{array}\right)
 &  
\underline{\cal W}^{ren}(\alpha) -  \underline{{\cal W}}^{ren}(\underline{\lambda}^{ren})  \end{array}\right) \, +\,  
 \underline{\cal W}_{(1)}(\alpha) \nonumber\\
\EEA
Compared to the arguments in \S \ref{subsubsection:perturbation-Lyapunov}, the leading order
matrix ${\cal H}_{(0)}$ has an $\alpha$-dependent prefactor. Thus (\ref{eq:lambda1}), (\ref{eq:pi1}) hold
provided one replaces ${\cal H}_{(0)}$ by ${\cal Z}(0){\cal H}_{(0)}$, $\bf 1$ by $v^{\dagger}_{(0)}$ (right null-eigenvector of ${\cal Z}(0){\cal H}_{(0)}$), and  ${\cal H}_{(1)}$ by 
\BEQ \del {\cal H}(\alpha):=({\cal Z}(\alpha)-{\cal Z}(0)) {\cal H}_{(0)} + {\cal Z}(\alpha){\cal H}_{(1)}(\alpha).  \label{eq:deltaH}
\EEQ 
Expanding $\del{\cal H}$ to order 1 in $\alpha$, one gets: 
\BEQ \del{\cal H}(0) = {\cal Z}(0) \Big( \left(\begin{array}{cc} 0 & 0 \\ 0 & 
\underline{\cal W}^{ren}(0)-\underline{\cal W}^{ren}(\underline{\lambda}^{ren}) \end{array}\right)
+ \underline{\cal W}_{(1)}(0) \Big) \label{eq:delH0}
\EEQ and 
\BEQ \del{\cal H}'(0) \sim {\cal Z}'(0) ({\cal H}_{(0)}
+{\cal H}_{(1)}(0)) + {\cal Z}(0) {\cal H}'_{(1)}(0). \label{eq:delpH0}
\EEQ

 
\Bigskip {\em Step 3. Order 0.}   By construction (up to normalization), the only right, resp. left, null eigenvector
of ${\cal Z}(0) {\cal H}_{(0)}$ is 
$v^{\dagger}_{(0)} := \left(\begin{array}{cc} 0 \\ 
{\bf 1} \end{array}\right)$, resp.
\BEQ \pi_{(0)}:=  z(0) {\cal Z}(0)^{-1} \left(\begin{array}{c} 
\underline{\pi}^{\perp}  \\ \underline{\pi}^{ren} \end{array}\right) 
=  z(0) \left(\begin{array}{c} 
 \, \underline{\pi}^{\perp} \\ 
 Z(\eps_p,0)_{{\cal G}_p}^{-1}  \,  \underline{\pi}^{ren}_{G_p} \\  \underline{\pi}^{ren}_{ext} 
\end{array}\right)
\EEQ
where 
\BEQ \underline{\pi}^{\perp} := - \Big(P^{\perp}(\, ^t\tilde{\cal W}^{int}-\Id)P^{\perp}\Big)^{-1}
\Big(-P^{\perp} \, ^t \tilde{\cal W}^{int} \Del(0) P^{av} +\,  ^t {\cal W}^{ext,\perp}(0) \Big) \underline{\pi}^{ren}.  \label{eq:piperp}
\EEQ  
We choose the normalization constant
\BEQ z(0):= \Big(  \sum_p \frac{\underline{\pi}^{ren}_{G_p}}{Z(\eps_p,0)_{{\cal G}_p}}  +    \langle {\bf 1}_{ext},
\underline{\pi}^{ren}_{ext} \rangle \Big)^{-1} 
\EEQ
so that $\langle v^{\dagger}_{(0)},\pi_{(0)}\rangle=1$. Since $\underline{\tilde{\pi}}^{ren}_{G_p}\sim 1$
by assumption, $z(0)\sim \min_p(Z(0)_{{\cal G}_p})$.  Note that 
\BEQ \pi_{(0)} \sim z(0)\  \Big\{ \sum_p Z(0)^{-1}_{{\cal G}_p} \,  \underline{\pi}^{ren}_{G_p} \ \del^{G_p}  +
 \underline{\pi}^{ren}_{ext} + O({\bf 1}^{\perp}) \Big\} \label{eq:pi_(0)}  
\EEQ 
where $O({\bf 1}^{\perp})$ is a sum of  error terms of order $O(1)$ along the $(\perp,p)$ components.


\Bigskip {\em Step 4. Lyapunov eigenvalue.} 
From (\ref{eq:lambda1}), 
\BEQ \lambda_{(1)} = -\frac{\langle \del{\cal H}(0) v^{\dagger}_{(0)},\pi_{(0)}\rangle}
{\langle \del{\cal H}'(0) v^{\dagger}_{(0)},\pi_{(0)}\rangle} 
=  - \frac{\langle {\cal Z}(0) {\cal H}_{(1)}(0) v^{\dagger}_{(0)}, \pi_{(0)}\rangle}{\langle ({\cal Z}'(0) ({\cal H}_{(0)} + {\cal H}_{(1)}(0))+   {\cal Z}(0) {\cal H}'_{(1)}(0) )   v^{\dagger}_{(0)}, \pi_{(0)}\rangle }  \label{eq:ren-lambda1}.
\EEQ
  The numerator is 
 $\sim  z(0)  \langle  {\cal H}_{(1)}(0) \left(\begin{array}{c} 0 \\ {\bf 1} 
 \end{array}\right), \left(\begin{array}{c} \underline{\pi}^{\perp} \\ \underline{\pi}^{ren}
 \end{array}\right)\rangle = z(0) \Big\{ \langle (\underline{\cal W}^{ren}(0)-\underline{\cal W}^{ren}(
 \underline{\lambda}^{ren})) {\bf 1}, \underline{\pi}^{ren}\rangle + \langle 
 \underline{\cal W}_{(1)}(0) {\bf 1}, \underline{\pi}^{\perp}\rangle \Big\}$. The absolute value of the
 first term between curly brackets, see (\ref{eq:related-id1}), is $\preceq \underline{\lambda}^{ren} 
 /\underline{k}^{ren}_{min}$. 
The second term is equal to  $\langle P_{\perp}(- d + \underline{\cal W}^{int\to ext}(0) {\bf 1}_{ext}),
 \underline{\pi}^{\perp} \rangle  = \langle P_{\perp} \eps, \underline{\pi}^{\perp}\rangle $ by (\ref{eq:id-intext}), and is therefore of order $O(Z(0)_{\cal G})$. 
 
\Medskip Now, the denominator is the sum of two terms; by (\ref{eq:id-dZ/dalpha}) and (\ref{eq:w(lambda*)sim3}), 
\BEA && \langle ({\cal Z}'(0) ({\cal H}_{(0)} + {\cal H}_{(1)}(0)) v^{\dagger}_{(0)}, \pi_{(0)}\rangle
 \sim  z(0)  \sum_p \tau_{{\cal G}_p} \langle P_{av,p} {\cal H}_{(1)}(0)   \left(\begin{array}{c} 0 \\ {\bf 1} 
 \end{array}\right), {\cal Z}(0)^{-1} \left(\begin{array}{c} \underline{\pi}^{\perp} \\ \underline{\pi}^{ren}
 \end{array}\right)\rangle   \nonumber\\
 &&\qquad\sim  z(0) \sum_p \tau_{{\cal G}_p} \langle P_{av,p} (\underline{\cal W}^{ren}(0)-
 \underline{\cal W}^{ren}(\underline{\lambda}^{ren})) {\bf 1}, {\cal Z}(0)^{-1} \left(
 \begin{array}{c} \underline{\pi}^{\perp} \\ \underline{\pi}^{ren} \end{array}\right) \rangle
 \nonumber\\ 
 &&\qquad\sim  z(0) \sum_p \tau_{{\cal G}_p}  \frac{\underline{\lambda}^{ren}}{\underline{k}^{ren}_{G_p}}
 Z(0)_{{\cal G}_p}^{-1} \underline{\pi}^{ren}_{G_p} \nonumber\\  \label{eq:Z'HHvpi}
\EEA
Since $ \underline{\pi}^{ren}_{G_p} \sim 1$,
 $\underline{k}^{ren}_{G_p} \sim \frac{1}{\tau_{{\cal G}_p}}
 \, \times\, Z(0)_{{\cal G}_p}$, we get: 
 $\sum_p \tau_{{\cal G}_p}  \frac{\underline{\lambda}^{ren}}{\underline{k}^{ren}_{G_p}}
 Z(0)_{{\cal G}_p}^{-1} \underline{\pi}^{ren}_{G_p} \sim \underline{\lambda}^{ren}
  \sum_p 
 (\underline{k}^{ren}_{{\cal G}_p})^{-2}$; now, $(\underline{k}^{ren}_{{\cal G}_p})^{-2} \preceq
  (\underline{k}^{ren}_{min})^{-2}$
   and  $\underline{\lambda}^{ren} \preceq \underline{k}^{ren}_{min}$ (by hypothesis), whence
 $\underline{\lambda}^{ren} \sum_p 
 (\underline{k}^{ren}_{{\cal G}_p})^{-2} \preceq (\underline{k}^{ren}_{min})^{-1}$.
Let us now deal with the second term in the denominator:
\BEQ \langle {\cal Z}(0) {\cal H}'_{(1)}(0)    v^{\dagger}_{(0)}, \pi_{(0)}\rangle 
\sim z(0) \, \Big\{ \langle (\underline{\cal W}^{ren})'(0) {\bf 1}, \underline{\pi}^{ren}\rangle
+ \frac{d}{d\alpha}\Big|_{\alpha=0} \langle P_{\perp}(-d(\alpha) + \underline{\cal W}^{int\to ext}(\alpha) {\bf 1}_{ext}), \underline{\pi}^{\perp} \rangle \Big\}
\EEQ
is a sum of two terms.  By (\ref{eq:w(lambda*)sim3}),  $ \langle (\underline{\cal W}^{ren})'(0) {\bf 1}, \underline{\pi}^{ren}\rangle 
 \sim  \langle \frac{1}{\underline{k}^{ren}_{\cdot}}, \underline{\pi}^{ren}\rangle \sim 
 (\underline{k}^{ren}_{min})^{-1}$. Then  $\frac{d}{d\alpha}\Big|_{\alpha=0} \langle P_{\perp}(-d(\alpha) + \underline{\cal W}^{int\to ext}(\alpha) {\bf 1}_{ext}) \sim -\frac{1}{k_{\sigma}} + 
\frac{d}{d\alpha}\Big|_{\alpha=0} (\frac{k^{ext}_{\sigma}}{k_{\sigma}+\alpha}) \sim -\frac{1}{k_{\sigma}} $ and $|\langle\frac{1}{k_{\cdot}},\underline{\pi}^{\perp}\rangle|\preceq \tau_{{\cal G}_p}$, which
is negligible compared to the other term. In total, we have proved that the denominator is 
$\sim z(0) (\underline{k}^{ren}_{min})^{-1}$, whence 
\BEQ |\lambda_{(1)}| \preceq |\underline{\lambda}^{ren}| + O(Z(0)_{\cal G}) \, \underline{k}^{ren}_{min}.
\EEQ
Note that, under our hypotheses, $|\lambda_{(1)}| \prec \underline{k}^{ren}_{min}$.


\Bigskip {\em Step 5. Lyapunov adjoint eigenvector.}
To first order, 
 \BEA v^{\dagger}_{(1)} &\sim &   - \ ({\cal Z}(0) {\cal H}_{(0)})^{-1} \ 
 \del{\cal H}(0) v^{\dagger}_{(0)}  \nonumber\\
 &=&  - ({\cal H}_{(0)})^{-1}\  \Big\{ (\underline{\cal W}^{ren}(0)-\underline{\cal W}^{ren}(\underline{\lambda}^{ren})){\bf 1}_{ren} + P_{\perp} (-d(0) + \underline{\cal W}^{int\to ext}(0) 
 {\bf 1}_{ext})) \Big\} \nonumber\\
 &\sim& - 
\left(\begin{array}{cc} P_{\perp} (\tilde{\cal W}^{int} - \Id) P_{\perp} & 0 \\
\left(\begin{array}{c} -P_{av}\Del(0) \tilde{\cal W}^{int} P_{\perp} \\ 
{\cal W}^{ext,\perp}(0) \end{array}\right) 
   &  
\underline{{\cal W}}^{ren}(\underline{\lambda}^{ren}) - \Id  \end{array}\right)^{-1}
\Big(  (\frac{\underline{\lambda}^{ren}}{\underline{k}^{ren}_{\cdot}})_{ren} + P_{\perp}\eps
\Big) \nonumber\\
& \equiv &   \left[ \begin{array}{c}  v^{\dagger}_{(1),\perp} \\  v^{\dagger}_{(1),ren} \end{array}
\right] 
\EEA
where we have used once again (\ref{eq:w(lambda*)sim3}) and (\ref{eq:id-intext}), and 
\BEQ v^{\dagger}_{(1),\perp} = - \Big( P_{\perp} ({\tilde{\cal W}}^{int} -\Id) P_{\perp} \Big)^{-1}   P_{\perp}\eps
\EEQ

\BEQ  v^{\dagger}_{(1),ren} =  
  \Big(\underline{\cal W}^{ren}(\underline{\lambda}^{ren}) - \Id\Big)^{-1}  \,  \Big\{ \frac{\underline{\lambda}^{ren}}{\underline{k}^{ren}_{\cdot}}  - (-P_{av}\Del(0) 
  \tilde{\cal W}^{int} +   {\cal W}^{ext,\perp}(0)) v^{\dagger}_{(1),\perp} \Big\}
\EEQ

\Medskip Thus $\perp$-component $v^{\dagger}_{(1),\perp}$ 
is of order $O(Z(0)_{\cal G})$, and   the  $(av,ext)$-components of $v^{\dagger}_{(1)}$  of 
order $O(\underline{Z}^{ren}) + O(Z(0)_{\cal G})$. 


\Bigskip {\em Step 6. Lyapunov weights.} 
To first order,
\BEQ  \pi_{(1)} \sim - (\, ^t {\cal H}_{(0)}  {\cal Z}(0))^{-1}  \ ^t \del{\cal H}(0) \pi_{(0)} 
\sim   - z(0) {\cal Z}(0)^{-1}   (\, ^t {\cal H}_{(0)})^{-1} u,   
\EEQ
 where $u^{\perp} = -P^{\perp} \, ^t \tilde{\cal W}^{int} D \underline{\pi}^{\perp}$ and
\BEQ  u^{ren} = (\, ^t\underline{\cal W}^{ren}(0)-\, ^t\underline{\cal W}^{ren}(\underline{\lambda}^{ren}))   \underline{\pi}^{ren}
+ \left(\begin{array}{c} -P^{av} D \\ ^t \underline{\cal W}^{int\to ext}(0)
\end{array}\right) \underline{\pi}^{\perp}.
\EEQ
Then
\BEQ (^t {\cal H}_{(0)})^{-1} u = \left(\begin{array}{cc} P^{\perp} (\, ^t \tilde{\cal W}^{int}-\Id)
P^{\perp} & (-P^{\perp}\, ^t \tilde{\cal W}^{int} \Del(0) P^{av} \ \   ^t {\cal W}^{ext,\perp}(0) )
\\ 
0 &  ^t \underline{\cal W}^{ren}(\underline{\lambda}^{ren})-\Id \end{array}\right)^{-1} u \nonumber\\
= \left[ \begin{array}{c} \pi_{(1)}^{\perp} \\ \pi_{(1)}^{ren} \end{array} \right] 
\label{eq:fixed1-tH0inv}
\EEQ
where
\BEQ \pi_{(1)}^{ren} = ( ^t \underline{\cal W}^{ren}(\underline{\lambda}^{ren})-\Id)^{-1} 
u^{ren}  \label{eq:5.143}  \EEQ
\BEQ  \pi_{(1)}^{\perp} = (P^{\perp} (\, ^t \tilde{\cal W}^{int}-\Id) P^{\perp})^{-1} \Big(u^{ren} 
- (-P^{\perp} \, ^t \tilde{\cal W}^{int} \Del(0) u^{av} + \,  ^t {\cal W}^{ext,\perp}(0) u^{ext})
\Big)  \label{eq:5.144}
\EEQ
From (\ref{eq:5.143}) and (\ref{eq:5.144}) it is apparent that $|| (^t {\cal H}_{(0)})^{-1} u||_{\infty} = O(1) ||u||_{\infty}$ for any $u$, i.e. that $|||(\, ^t {\cal H}_{(0)})^{-1}|||=O(1)$  for the matrix norm associated to the $L^{\infty}$-norm on the components. 
Next,  
\BEQ \Big((\, ^t\underline{\cal W}^{ren}(0)-\, ^t\underline{\cal W}^{ren}(\underline{\lambda}^{ren}))
\underline{\pi}^{ren}\Big)_{\sigma} = \sum_{\sigma'} \underline{\pi}^{ren}_{\sigma'} (\underline{w}^{ren}_{\sigma'\to\sigma}(0) - \underline{w}^{ren}(\underline{\lambda}^{ren})_{\sigma'\to\sigma}) \nonumber\\
\EEQ
From (\ref{eq:w(lambda*)sim3}), $\underline{w}^{ren}_{\sigma'\to\sigma}(0) - \underline{w}^{ren}(\underline{\lambda}^{ren})_{\sigma'\to\sigma} \sim \underline{w}^{ren}_{\sigma'\to\sigma} 
\, \times\, \frac{\underline{\lambda}^{ren}}{\underline{k}^{ren}_{\sigma'}}$. Now, 
$ \frac{\underline{\lambda}^{ren}}{\underline{k}^{ren}_{\sigma'}}\le \underline{Z}^{ren}$ is $o(1)$, and  
$\sum_{\sigma'} \underline{\pi}^{ren}_{\sigma'} \underline{w}^{ren}_{\sigma'\to\sigma} = (\, ^t \underline{\cal W}^{ren} \underline{\pi}^{ren})_{\sigma} = \underline{\pi}^{ren}_{\sigma} $. Thus
$||u^{ren}||_{\infty} = O(\underline{Z}^{ren}) + O(Z(0)_{\cal G})$. 

\Medskip
As a result, components of $\pi_{(1)}$ are of
  order $O(\underline{Z}^{ren}) + O(Z(0)_{\cal G})  = o(1)$ compared to the components
of $\pi_{(0)}$.


\subsubsection{Full induction, first part: 'stopped' graph}  \label{subsubsection:full-non-auto}


\noindent The construction in the previous paragraph may be iterated. Referring to 
the coalescence tree described in \S \ref{subsection:merging}, one gets the branches from the leaves, terminating at the first autocatalytic vertices (called: {\em terminal vertices}) met on the way down to the root. Namely,
we have explicitly assumed in \S \ref{subsubsection:one-renormalization-step} that none of the subgraphs of the previous steps is autocatalytic.  This means that renormalized couplings to or from autocatalytic vertices are replaced by degradation rates, i.e. $\sigma\overset{k}{\to}\sigma'$ is replaced by $\sigma\overset{k}{\to}\emptyset$ if $\sigma$ or $\sigma'$ is autocatalytic. The resulting graph, $G^{stop}$, is called the {\bf 'stopped' graph}; after merging, it is simply made up of isolated terminal vertices with degradation rates, which
are the images of the connected components $(G^{stop}_q)_{q=1,2,\ldots}$ of $G^{stop}$. Cutting edges outgoing from terminal vertices ensures that the only hierarchical formula to be
checked for 
the Lyapunov weights of the 'stopped' graph is (\ref{eq:pisigma1}). Cutting edges
going into terminal vertices ensures that (\ref{eq:vdaggersigma}) simply states: 
$v^{\dagger}_{\sigma}\sim 1$. Note that there can be no edge  in $G$ outgoing form $G^{stop}_q$ if $G^{stop}_q$ is not autocatalytic (otherwise one of them would be dominant after renormalization, which would contradict our stopping rule), so that non-autocatalytic components $G^{stop}_q$ are also terminal in the whole graph $G$.

\Medskip  We now prove (\ref{eq:pisigma1}), (\ref{eq:vdaggersigma}) for the  stopped graph 
$G^{stop}$. We may restrict to one of the components, $G^{stop}_q$, and let $i_{max,q}$ be
the last merging step for $G^{stop}_q$. Iterate the renormalization steps.
At step $i=1$,  ${\cal G} = G^{(0)}_{\searrow n(1)}$, $\underline{G} = G^{(0)}$. Then, 
${\bf 1}^{int,p} = {\bf 1}_{V(G_p)}$, ${\bf 1}^{int} = {\bf 1}_{V(G)}$,  $\tilde{\pi}^{int,p} 
= \tilde{\pi}^{int}_{{\cal G}_p}$, and $\underline{\pi}^{ren} = \pi^{*,(1)}$ are the Lyapunov
weights of the renormalized network with adjoint generator $\underline{A}^{ren}$. Then, assuming $\pi^{*,(1)}_{G_p}\sim 1$ for all $p$,  we have proved:
\BEQ \pi^*  \sim z^{(1)}(0) \  \Omega^{(1)} \pi^{*,(1)}  \EEQ
where $\Omega^{(1)}:\Sigma^{(1)}\to \Sigma^{(0)}=\Sigma$ is the {\bf weight renormalization matrix} defined by
\BEQ \Omega^{(1)}_{G_p,\sigma} = Z(0)_{{\cal G}_p}^{-1}, \qquad \sigma\in V({\cal G}_p) \EEQ
\BEQ \Omega^{(1)}_{\sigma',\sigma'} = 1,\qquad \sigma'\in \Sigma^{(1)}\setminus\{G_p,p\ge 1\}
\simeq \Sigma^{(0)}\setminus\{V({\cal G}_p),p\ge 1\} 
\EEQ
 all other coefficients vanishing.  
Under similar hypothesis, we get 
\BEQ \pi^{*,(i)} \sim z^{(i+1)}(0)    \  \Omega^{(i+1)} \pi^{*,(i+1)}, \qquad i\ge 1  \EEQ
Solving by induction on the number $i'_{max}$ of renormalization steps ($i'_{max}\le i_{max,q}$
because of the 'stopping' procedure), and using the results 
of \S \ref{subsection:cut-off-Lya-estimates} to initialize the induction ($i'_{max}=1$) -- a 
case when it has been proved that $\pi^*\sim {\bf 1}$ so that the main assumption holds --,    we get the hierarchical formula (\ref{eq:pisigma1}) for the Lyapunov weights of the 'stopped' graph $G^{stop}$, up to normalization by the product of the prefactors $\prod_{i=1}^{i'_{max}} z^{(i)}(0)$. The normalization in this paragraph ensures that $\sum_{\sigma} \pi^{*,(i)}_{\sigma}\sim 1$.

\Medskip Proceeding similarly for the adjoint Lyapunov eigenvectors $v^{\dagger}$, we have proved:
\BEQ v^{*\dagger} \sim \Omega^{\dagger,(1)} v^{*\dagger,(1)} \EEQ
where  
\BEQ \Omega^{\dagger,(1)}_{G_p,\sigma} = 1, \qquad \sigma\in V({\cal G}_p) \EEQ
\BEQ \Omega^{\dagger,(1)}_{\sigma',\sigma'} = 1,\qquad \sigma'\in \Sigma^{(1)}\setminus\{G_p,p\ge 1\}
\simeq \Sigma^{(0)}\setminus\{V({\cal G}_p),p\ge 1\} 
\EEQ
Thus  we
simply get by induction
\BEQ  v^{*\dagger,(i)} \sim {\bf 1}. \EEQ
This is in agreement with the hierarchical formula (\ref{eq:pisigma1}) for the adjoint Lyapunov
eigenvector. With this normalization, $\langle v^{*\dagger,(i)},\pi^{*,(i)}\rangle \sim 1$.

\Bigskip There is a way to directly interpret formulas (\ref{eq:pisigma1}) for 
$\pi^*_{\sigma}$ with $\sigma\in V(G^{stop})$ in terms of sums over paths, which we shall
now spell out. We call $\sigma$  {\bf $\pi$-dominant} if  $\pi^*_{\sigma}\succeq \pi^*_{\sigma'}$ for all $\sigma'$ in the same connected component of $G^{stop}$. 

\Bigskip {\bf Leading  paths (non autocatalytic case).} Let $\sigma,\sigma'$ lie in the same {\em non-autocatalytic} connected
component $G^{stop}_{q}$ 
with threshold rate $\alpha_q=0$  of $G^{stop}$.   By hypothesis, $\sigma,\sigma'$ merge
 at some step $i_{\sigma,\sigma'}$, i.e. there exists $i=i_{\sigma,\sigma'}\le i_{max,q}$ such that $G^{(i-1)}(\sigma)
 \not=G^{(i-1)}(\sigma'), G^{(i)}(\sigma)=G^{(i)}(\sigma')\equiv G_{\sigma,\sigma'} \subset
 G^{stop}$.   Thus there exists a simple path $\gamma:\sigma\to\sigma'$ descending from $\sigma$ within $G_{\sigma,\sigma'}$,
 in the terminology of \S \ref{subsection:merging} (see Fig. \ref{fig:merging}.1, \ref{fig:merging}.2 for an illustration, and discussion at the end of the subsection) such that: (i) paths {\em circulating within} $G_i$'s
 are made of dominant edges; (ii) {\em leaving} edges $\sigma_1\to \sigma_2$, which leave
 $G_{i_1} = G^{(i_1)}(\sigma_1)$ from $\sigma_1$, have small weight, but the associated edge $G_{i_1}\to G^{(i_1)}(\sigma_2)$
 in $G^{(i_1)}$ is dominant; (iii) {\em reentering} edges $\sigma_1\to\sigma_2$ are dominant.  
 Such a path is called {\bf  simple leading path} from $\sigma$ to $\sigma'$; in the example of 
 \S \ref{subsection:merging}, see Fig. \ref{fig:merging}.2, both $(0\to 1\to 2\to 6)$ and 
 $(0\to 1\to 2\to 3\to 4\to 5\to 6)$ are simple leading paths from $0$ to $6$, provided 
 the $G_i$ are step $i$ maximal dominant SCCs, and edges $(0\to 1\to  2\to 0),(4 
 \rightleftarrows 5),(7\rightleftarrows 8)$ within $G_1,G_2,G_4$; leaving edges $(G_1\to 3)$,
 resp. $(G_2\to 3)$  
 resp $(G_3\to 6$, resp.   $(G_4\to 6)$, resp. $(G_5\to G_3)$ (coming from non-dominant $(2\to 3)$, resp. $(4\to 4)$,  resp. $(2\to 6)$ and $(5\to 6)$, resp. $(8\to 6)$, resp. $(6\to 5)$),
 and reentering edges $(3\to G_1)$, resp. $(3\to G_2)$, resp. $6\to G_4$ (coming from
 dominant  $(3\to 2)$, resp. $(3\to 4)$, resp. $(6\to 7)$), are dominant. 
The small weight of {\em leaving} edges is compensated  by summing  over paths of bounded
(but possibly large) length circling 
around compound vertices. This can be proved by induction on the number of steps $i$, starting from $i=1$. Assume
$G^{(i_1-1)}(\sigma_1)=\sigma_1$ and $G^{(i_1-1)}(\sigma_2)=\sigma_2$, so that vertices
are bare vertices, and there is only one barrier to cross. 
Exit probabilities $f^{\sigma_1}(\sigma_2)$ (see Lemma \ref{lem:exit-proba}) compute the sum of the weights of 
paths circulating within $G_{i_1}$ before leaving from $\sigma_1$. Then $f^{\sigma_1}(\sigma_2)
\sim w(\alpha_q)_{G_{i_1}\to \sigma_2}$, and $ w(\alpha_q)_{G_{i_1}\to \sigma_2}\sim 1$ by 
assumption (ii). Reordering the sum over paths according to increasing length, we see
that truncating to large enough rank still yields a weight $\sim 1$, for every $\alpha\sim 
\alpha_q$. The resulting path transition weights (of the same order as the edge
transition weights of  the renormalized network) may be reused to show a similar result
for $i=2$, and so on.  We call {\bf leading path set} the obtained finite path set, and write $\sum_{\gamma:\sigma\to\sigma' \ {\mathrm{leading\ path}}} 
(\cdots)$ for path sums restricted
to $\gamma$ in the leading path set.

\Medskip Thus the weight $\sum_{\gamma:\sigma\to\sigma' \ {\mathrm{leading\ path}}} $  of such a truncated sum path is $\sim 1$, times a product of 
factors corresponding to the barriers that are {\em not} crossed in the final state $\sigma'$,
on the descending sequence from $\sigma'$. Since 
  (from (\ref{eq:w(alpha)-free-autocata-regimes}))
   $ w(\alpha)_{G_{i_1}\to \sigma_2} \sim Z_{G_{i_1}}^{-1} \, \times \max_{\sigma \in V(G_{i_1})}  w_{\sigma_1\to \sigma_2}$ for a dominant leaving edge $(G_{i_1}\to\sigma_2)$
    coming from $(\sigma_1\to \sigma_2)$,  the missing $w_{\sigma_1\to\sigma_2}$-factor
     amounts to a correcting factor $Z^{-1}_{G_{i_1}}$. The product of these factors
      yields (\ref{eq:pisigma1}). 

\Medskip {\bf Leading excursions (non autocatalytic case).} A very related notion is that of {\em leading excursion}.  We still assume $\alpha_q=0$ for the moment. Namely, a {\em leading excursion from $\sigma$ to $\sigma'\not=\sigma$}  is a leading path $\gamma:\sigma=\sigma_0\to\cdots \to \sigma_{\ell}=\sigma'$ which is also an excursion from $\sigma$, i.e. such that $\sigma_i\not=\sigma$, $i\ge 1$. We write similarly $\sum_{\gamma:\sigma\to\sigma' \ {\mathrm{leading\ excursion}}} 
(\cdots)$ for path sums restricted
to $\gamma$ in the leading excursion set. The correcting factor in the weight ratio 
$\frac{\sum_{\gamma:\sigma\to\sigma' \ {\mathrm{leading\ excursion}}} w(\alpha_q)_{\gamma}}{ 
\sum_{\gamma:\sigma\to\sigma' \ {\mathrm{leading\ path}}} w(\alpha_q)_{\gamma}}$ 
is a product of $Z$-factors, equivalent to $\pi_{\sigma}^{-1}$ in the normalization
of (\ref{eq:pisigma1}). Hence, {\em when  $\sigma$ is $\pi$-dominant, so
that $\pi^*_{\sigma}\sim 1$ and all $\pi^*_{\sigma'}\preceq 1$}, we get: 
 \BEQ \sum_{\gamma:\sigma\to\sigma' \ {\mathrm{leading\ excursion}}} w(\alpha_q)_{\gamma} \sim \pi^*_{\sigma'}  \label{eq:leading-excursion-pi}
 \EEQ
  in the normalization of this paragraph. 
  
\Medskip {\bf Leading excursions (completed).} In the non-autocatalytic case,  path sums \\ $\sum_{\gamma:\sigma\to \sigma', \gamma\subset G_q^{stop}} w(\lambda^*_q)_{\gamma}$  diverge upon replacing $\alpha_q=0$ by the true Lyapunov exponent 
$\lambda_q:=\lambda^*(G^{stop}_q)$, since the sum $\Phi_{G^{stop}_q}(\lambda^*_q)_{\sigma}$ 
of the weights of excursions $\gamma:\sigma\to\sigma$ within $G^{stop}_q$ is equal to 1, see
 (\ref{eq:excursion}). Hence leading paths are actually ill-defined. The same path sums 
 diverge again and for the same reason in the autocatalytic case (leading paths can be defined again
 provided $\alpha_q\sim \lambda^*(G_q^{stop})$ is well-chosen in order to avoid divergences,
 as in \S \ref{subsubsection:exit-proba}). Fortunately, eq. (\ref{eq:leading-excursion-pi})
 survives the replacement of $\alpha_q$ by $\lambda^*_q$, and extends to the autocatalytic case.
Namely,
\begin{itemize}
	\item[(i)] {\em (subdivergences)} Sums over paths circulating in an inner SCC ${\cal G}_i\subsetneq G^{stop}_q$ are insensitive to the value of $\lambda^*_q$, for $G_i$ is not autocatalytic 
	(even if $G^{stop}_q$ is),  and $\lambda^*_q\prec k^{ext}_{{\cal G}_i} = \frac{1}{\tau_{{\cal G}_i}} \, \times Z(0)_{{\cal G}_i}$ (in the denominators of the coefficients
	of  the weight matrix	${\cal W}(\lambda^*)$)  gives a subleading correction to the $Z$-factors; 
	\item[(ii)] {\em (global divergence)} Leading excursions from any $\sigma\in V(G^{stop}_q)$ 
	are globally convergent. This means the following. Let $G_{i_1}=G^{(i_1)}(\sigma)$ be
	first non-trivial compound vertex containing $\sigma$; sums over paths circling around
	$G_{i_1}$ while avoiding $\sigma$ have bounded weight $O(1)$  -- instead of $O(Z_{G_{i_1}}^{-1})$ without the restriction --. Therefore sums over paths circling around $G_{i_1}$ and
	then leaving $G_{i_1}$ have small weight $O(Z_{G_{i_1}})$. Leading excursions may thus
	simply avoid circling around $G_{i_1}$ in the first place, i.e. avoid reentering 
	$G_{i_1}$ altogether. Proceeding by induction, one sees that leading excursions may actually avoid reentering $G_{i_{max}(\sigma)}$, which is the largest compound vertex $\subset  
	G^{stop}_q$ containing $\sigma$. In turn, this means that one may avoid circling 
	around vertices of the  last step merged diagram $G_q^{(i_{max,q})}:= (G^{stop}_q)^{(i_{max,q})}$. 'Avoiding'
	means that neglected paths have subleading weight. Then the remaining path sum is 
	insensitive to  the value of $\lambda^*_q$.
\end{itemize}

\Medskip {\em We may now relate leading excursions to the asymptotic Green kernel $G^*(\sigma,\sigma')$ (see (\ref{eq:asympt-Green})) in the case when 
   $\sigma$ is $\pi$-dominant}. We have just proved:
 $\sum_{\gamma:\sigma\to\sigma'\ {\mathrm{leading}}} w(\alpha_q)_{\gamma} \sim \pi^*_{\sigma'}$, an estimate which holds true also (non-autocatalytic case), resp. in particular (autocatalytic case)  when $\alpha_q=\lambda^*$.
  Let $\ell_{\sigma\to\sigma'} = 
\max_{\gamma:\sigma\to\sigma'\ {\mathrm{leading\ excursion}}} \ell(\gamma)$.   Then (decomposing paths
$\gamma:\sigma\to\sigma'$ into $\gamma_0\cup\gamma'$, where $\gamma_0:\sigma\to\sigma$ and 
$\gamma':\sigma\to\sigma'$ is an excursion from $\sigma$)
\BEA G^*(\sigma,\sigma') &\ge& \lim_{N\to\infty} \frac{1}{N} \sum_{\gamma_0:\sigma\to\sigma, 
\ell(\gamma)\le N-\ell_{\sigma\to\sigma'}} w(\lambda^*)_{\gamma} \, \times\, 
\sum_{\gamma':\sigma\to\sigma'\ {\mathrm{leading\ excursion}}} w(\lambda^*)_{\gamma'}
\nonumber\\
&=& G^*(\sigma,\sigma) \, \times\, \sum_{\gamma':\sigma\to\sigma'\ {\mathrm{leading\ excursion}}} w(\lambda^*)_{\gamma'} \nonumber\\
&\sim & G^*(\sigma,\sigma) \, \times\, \pi^*_{\sigma'}  \label{eq:G-decomp-ineq}
 \EEA
As we shall see presently, $G^*(\sigma,\sigma)= v^{\dagger,*}_{\sigma} \pi^*_{\sigma} \sim 1$ and $G^*(\sigma,\sigma')= v^{\dagger,*}_{\sigma} \pi^*_{\sigma'} \sim \pi^*_{\sigma'}$. Thus the inequality in (\ref{eq:G-decomp-ineq})
is actually an equivalence (i.e. both sides are of the same order): {\em leading 
excursions provide the main contribution to the Green kernel}.


\subsubsection{Full induction, second part} \label{subsubsection:full-auto}


\noindent The idea is now to incorporate the 'low-lying' 
edges, i.e. the lower-scale edges  $\sigma\overset{k}{\to}\sigma'$ replaced by degradation
rates $\sigma\overset{k}{\to} \emptyset$ through the 'stopping' procedure.  If there are none, then the hierarchical formulas are proved; this includes the case when $A$ is the generator
of a non-deficient Markov chain,  $\pi^*$ is the stationary measure of the associated skeleton
chain, and $v^{\dagger}={\bf 1}$.  

\noindent The proof relies generally speaking on an identification of the Lyapunov and adjoint 
Lyapunov eigenvector through the {\bf asymptotic Green kernel} 
\BEQ G^*(\sigma,\sigma'):= \lim_{N\to\infty} \frac{1}{N} G^{\lambda^*}_N(\sigma,\sigma'), \qquad G^{\lambda}_N(\sigma,\sigma'):= \sum_{\ell=0}^{N-1} (({\cal W}(\lambda))^{\ell})_{\sigma,\sigma'} ;  \label{eq:asympt-Green}
\EEQ
 let us spell out the argument (see heuristic section, \S \ref{subsection:heuristics}).  Let 
 $P^* := {\cal W}(\lambda^*)$. 
If $G$ is strongly connected,  there exists a unique couple $(\pi^*,v^{\dagger,*})>0$ s.t. $\pi^* P^* = \pi^*, P^* v^{\dagger,*}=v^{\dagger,*}$ and 
$\sum_{\sigma} \pi^*_{\sigma}=1$, $\langle v^{\dagger,*},\pi^*\rangle=1$. Up to a finite 
averaging procedure (see proof of Corollary \ref{cor:eta}), one may assume that all eigenvalues
other than 1 have modulus $<1$, whence $((P^*)^n)_{\sigma,\sigma'} \sim_{n\to\infty} 
v^{\dagger,*}_{\sigma} \pi^*_{\sigma'}$. Then, by Cesaro averaging, we get
\BEQ \frac{1}{N} G^*_N(\sigma,\sigma') \to_{N\to\infty} v^{\dagger,*}_{\sigma} \pi^*_{\sigma'}.
\EEQ
Now, if one can prove the existence of two vectors $v^{\dagger},\pi>0$ with $\sum_{\sigma}\pi_{\sigma}\sim 1$ such that 
$G^*(\sigma,\sigma')\sim v^{\dagger}_{\sigma} \pi_{\sigma'}$, one gets  $v^{\dagger}\sim 
v^{\dagger,*}$ (by summing over $\sigma$), and then $\pi\sim \pi^*$.

\Medskip (1) {\em Consider first the case when $\sigma,\sigma'\in V(G_q^{stop})$ are 
vertices of  the same connected component of the  'stopped' graph.}  
We claim that an  edge $e_{ext}=(\sigma,\sigma')$ not in $E(G_q^{stop})$   can be incorporated into $E(G_q^{stop})$ without changing the
renormalization features (merging steps, renormalized scales). Namely, $\sigma,\sigma'$ merge
 at some step $i\le i_{max}$, i.e. there exists $i\le i_{max}$ such that $G^{(i-1)}(\sigma)
 \not=G^{(i-1)}(\sigma'), G^{(i)}(\sigma)=G^{(i)}(\sigma')\equiv G_{\sigma,\sigma'} \subset
 G_q^{stop}$. Then (considering the $i$-step renormalized graph associated to $G_q^{stop}$)
 $w(\alpha_{G_{\sigma,\sigma'}})_{\sigma\to\sigma'}\sim 1$ receives a contribution of 
 lower order from $e_{ext}$ since $w(\alpha_G)_{e_{ext}}\preceq Z_{G_{\sigma,\sigma'}}\prec 1$.  

\Bigskip Thus, one may assume that any edge $e_{ext}=(\sigma,\sigma')\in E(G^{(0)})\setminus E(G^{stop})$ is (i) either {\em external}, i.e.  such that $\sigma,\sigma'$ do not
belong both to $V(G^{stop})$; (ii) or internal but {\em externally connecting}, i.e. $\sigma,\sigma'\in V(G^{stop})$, but  there exists no path connecting them within $G^{stop}$, i.e. 
$\sigma\in G^{stop}_q, \sigma'\in G^{stop}_{q'}$ with $q\not =q'$. By extension, an {\em externally connecting path} $\gamma:\sigma\to\sigma'$ is a path $(\sigma=\sigma_0\to\cdots
\to \sigma_{\ell}=\sigma')$  such that $\sigma_1,\ldots,\sigma_{\ell-1}\not \in V(G^{stop})$. 
By construction, $G_q^{stop}\ni \sigma$ must be autocatalytic; also, $w(\alpha_q)_{G_q^{stop}\to \sigma_1} \sim \frac{k_{G_q^{stop}\to \sigma_1}}{\alpha_q} \preceq \frac{k_{G_q^{stop}}}{\alpha_q} \prec 1$ by (\ref{eq:outgoing-autocata-small-weight}). Also, still by construction,
all {\em purely external} loops $\gamma:\sigma_0\to\cdots\to\sigma_{\ell}=\sigma_0$ (i.e.
loops 
 with $\sigma_0,\ldots,
\sigma_{\ell-1}\not\in V(G^{stop})$) have weight $w(\alpha_q)_{\gamma}\prec 1$.  This implies by an 
easy inductive argument on the length that $\sum_{\gamma:\sigma\to\sigma'\ {\mathrm{externally \ connecting}}}
w(\lambda)_{\gamma} \preceq Z_{G_q^{(i_{max,q})}}\, \times\,  w(\alpha_q)_{G_q^{stop}\to \sigma_1} \prec 1 $ for any $\lambda\succeq \lambda^*_q$ and $b$ (scale parameter)
large enough.
 
\Bigskip (2)  {\em (autocatalytic case)} Assume some $G_q^{stop}$ are autocatalytic, including
$G_1^{stop}$,
and $\alpha_1\succeq \alpha_q>0$ for all autocatalytic $G_q$, $q\not=2$, i.e. $G_1^{stop}$ is a core. Recall
 from \S \ref{section:description} that there is by assumption (banning resonance cases) only one core in that case, i.e. $\alpha_1\succ \alpha_q$
 for all autocatalytic $G_q$, while non-autocatalytic ones have $\alpha_q=0$.   

\Medskip To begin with,  $\lambda^*\sim \lambda^*_1$. Namely, the Lyapunov eigenvalue (defined through the implicit excursion equation 
(\ref{eq:excursion})) can only increase by putting back low-lying edges, hence $\lambda^*\ge 
\lambda^*_1$. On the other hand, if $\alpha\succ \lambda^*_1$, the sum 
 $\sum_{\gamma:\sigma\to \sigma} w(\alpha)_{\gamma}$ of weights of loops $\gamma$ circling around $G^{stop}_q$
becomes finite, of order $(Z_{G_q^{(i_{max,q})}})^{-1}$. Multiplying by the weight of 
an outgoing edge $w(\alpha)_{\sigma'\to\sigma_1}$, $\sigma'\in V(G_q^{stop}), \sigma_1\not 
\in V(G_q^{stop})$ yields a factor of order $w(\alpha_q)_{G_q^{stop}\to \sigma_1} \preceq 
\frac{k_{G_q^{stop}}}{\alpha_1}  \prec 1$
at most, which is an upper bound for  the multiplicative weight associated to  reentering and leaving $G_{q}^{stop}$. Thus the whole series 
$\sum_{\gamma:\sigma\to\sigma} w(\alpha)_{\gamma}$ is convergent, which implies that 
$\lambda^*<\alpha$.  

\Medskip 
 Fix a  $\pi$-dominant vertex $\sigma^*\in V(G_1^{stop})$, and let  $\sigma'\in \Sigma$. We combine the previous
 arguments to prove that the vectors $G^*(\sigma^*,\cdot) \propto \pi$ are equivalent up to 
 normalization, where $\pi$ is given by the general hierarchical formulas (\ref{eq:pisigma1}),
 (\ref{eq:pisigma2}). We decompose $\gamma:\sigma^*\to\sigma'$ as in (\ref{eq:G-decomp-ineq})
 into $\gamma_0\cup\gamma'$, where $\gamma_0:\sigma^*\to\sigma$ and $\gamma':\sigma^*\to\sigma'$
 in an excursion from $\sigma^*$. We restrict first to $\gamma'$ in the form of {\em generalized
 leading excursions}. 
 
\Medskip Assume first  $\sigma'\not\in V(G^{stop})$, a generalized leading excursion 
 $\gamma':\sigma^*\to\sigma'$ is a leading excursion within $G^{stop}_1$ from $\sigma^*$ to
 some $\sigma_1 \in V(G^{stop}_1)$, followed by a simple, purely external
  path $\gamma'':\sigma_1\to\sigma'$. The lower bound in (\ref{eq:G-decomp-ineq}) yields
 \BEQ G^*(\sigma^*,\sigma') \succeq G^*(\sigma^*,\sigma^*) \, \times\, Z_{G_1^{(i_{max,1})}}\ \pi_{\sigma'}, \label{eq:lower-bound-G-1}
 \EEQ
  with
  $\pi$ as in  the first line of  (\ref{eq:pisigma2}). Further contributions to 
  $G^*(\sigma^*,\sigma')$ , other than from subleading self-intersecting, purely external paths $\gamma''$, involve paths reentering $G^{stop} = \cup_q G^{stop}_q$ a number of times. Now, the  multiplicative weight associated to  reentering and leaving $G_{q}^{stop}$, $q\not=1$ is, as above, $\preceq \frac{k_{G_{q}^{stop}}}{\alpha_{q}} \prec 1$, since  $\lambda^*\sim \lambda^*_1\succ \lambda^*_q$ by assumption.
 Paths reentering and leaving $G^{stop}_1$, but avoiding $\sigma_1\in V(G^{stop}_1)$ (see 
 (ii) above in the discussion of leading excursions), contribute similarly a factor
     $\preceq \frac{k_{G_{1}^{stop}}}{\alpha_{1}} \prec 1$. Resumming all these contributions
   yields a subleading term  compared to the r.-h.s. of (\ref{eq:lower-bound-G-1}). 

\Medskip The case when   $\sigma'\in V(G^{stop})$ is similar, with the supplementary
multiplicative term in the form of product of $Z^{-1}$-factors in the second 
line of (\ref{eq:pisigma2}) arising as discussed in \S \ref{subsubsection:full-non-auto} just
before the introduction of leading excursions. 
 
\Medskip (3) This settles the hierarchical formulas for $\pi$ in the  autocatalytic case. 
Let us briefly discuss $v^{\dagger}$.  Let $\sigma,\sigma'\in \Sigma$. As discussed above, the sum of weights $w(\lambda^*)_{\gamma_q}$ of loops $\gamma_q$
circling around $G_q^{stop}$ is finite for $q>1$. Therefore the contribution to the 
asymptotic Green function of paths not entering $G_1^{stop}$ is identically zero. The twofold
decomposition of a path  in the path sum $\gamma:\sigma\to\sigma'$  into $\gamma_0\cup\gamma'$,
used for $\sigma\in V(G_1^{stop})$, can therefore be refined to a  fourfold decomposition
$\gamma=\gamma^{\dagger}\cup \gamma^{\dagger}_1\cup\gamma_0\cup\gamma'$, where:

\begin{itemize}
\item[(i)]  $\gamma^{\dagger}:\sigma=\sigma_0\to\cdots\to \sigma_{\ell} = \sigma_1$ is a path connecting $\sigma$ to $\sigma_1\in V(G_1^{stop})$ such that $\sigma_0,\ldots,\sigma_{\ell-1}\not\in V(G_1^{stop})$, 
i.e. $\ell$ is the entrance time into $G_1^{stop}$ ($\ell=0$ if $\sigma\in V(G_1^{stop})$). 
The leading contribution to the path sum is expressed and computed as in (2), and 
yields the hierarchical formula (\ref{eq:vdagger1}) for $v^{\dagger}_{\sigma}$. In particular,
$v^{\dagger}_{\sigma^*} \sim 1$, which is coherent with the fact that $\gamma^{\dagger}$ is 
trivial if $\sigma=\sigma^*$. 

\Medskip  
If $\sigma_1=\sigma^*$, then the decomposition of $\gamma\setminus\gamma^{\dagger}$ into
$\gamma_0\cup\gamma'$ is as before; if not,

\item[(ii)]  we sum over leading 'reverse' excursions
$\gamma^{\dagger}_1:\sigma_1\to \sigma^*$ 
connecting $\sigma_1$ to $\sigma^*$  within $V(G^{stop}_1)$ (compare to discussion of leading excursions $\sigma^*\to \sigma_1$ in \S \ref{subsubsection:full-non-auto}), which are leading paths $\gamma_1:\sigma_1\to \cdots\to \sigma_{\ell+1}=\sigma^*$ in $V(G^{stop}_1)$ such that $\sigma_1,\ldots,\sigma_{\ell}\not=\sigma^*$;  the sum over reverse leading excursions yields
a factor $\sim 1$;

\item[(iii,iv)] then the decomposition of $\gamma\setminus (\gamma^{\dagger}\cup \gamma^{\dagger}_1): \sigma^*\to \sigma'$ into  $\gamma_0\cup\gamma'$  is as in (2), 
yielding to leading order $\sim G^*(\sigma^*,\sigma^*) \, \times\, Z_{G_1^{(i_{max,1})}}\ \pi_{\sigma'}$. 

  \end{itemize}

All together, we have shown: $G^*(\sigma,\sigma')\sim v^{\dagger}_{\sigma} G^*(\sigma^*,\sigma')$, implying: $v^{\dagger,*}_{\sigma}/v^{\dagger,*}_{\sigma^*} \sim v^{\dagger}_{\sigma}/
v^{\dagger}_{\sigma^*}$.

\Bigskip (4) {\em (other cases)} Recall from \S \ref{subsection:description} that non-autocatalytic $G_q$'s are terminal, i.e. they have no outgoing edges.   If $G^{stop}$ is not autocatalytic, i.e. no  $G_q^{stop}$ is, 
then there are no low-lying edges, so that $G^{stop}=G$: the renormalization procedure has already achieved the
goal of proving hierarchical formulas. As mentioned in the beginning of this paragraph, this
includes the case  when $A$ is the generator of a non-deficient Markov chain, so that 
$\lambda^*=0$, but also some Markov chains with low degradation rates. Namely, if some
degradation reaction $\sigma\to\emptyset$ becomes dominant at some step $i\le i_{max}$, 
the cemetery state $\emptyset$ becomes a maximal dominant SCC, and our formulas simply 
give $\pi=\del_{\emptyset}$, which is a trivial statement (the cemetery state is an absorbing state). If not, renormalization yields Lyapunov data for each $G^{stop}_q=G_q$, with
$\lambda^*(G_q)\sim \lambda_{{\cal G}_q}\sim -Z(0)_{{\cal G}_q}/\tau_{{\cal G}_q} <0$,
$v^{\dagger,*}(G_q)\sim {\bf 1}_{V(G_q)}$ and $\pi^*(G_q)$ as in the hierarchical formula  (\ref{eq:pisigma1}).     
  
\Medskip  All these formulas also apply without any modification to the cut graphs $G_{cut}(i)$
of \S \ref{section:description}, and to 'restricted' cut graphs $(G_{cut}(i))_{\sigma_0}$ 
when $G_{cut}(i)$ is not strongly connected.


\subsection{Directed acyclic graph (DAG) structure}  \label{subsection:DAG}


\noindent We finally show how to describe by a finite algorithm an estimate of the maximum 
$\max_{\gamma:\sigma^*\to G^{(i)}(\sigma)} w(\alpha_{\sigma_0})_{\gamma}$  in the first 
line of (\ref{eq:pisigma1}). 

\Medskip Recall that a DAG is a directed graph 
$\T=(V,{\cal E})$
containing no cycle. Writing $x<x'$ if there exists a path from $x$ to $x'$ in $\T$ defines
a partial ordering on $\T$; minimal elements for this ordering are called roots.   An example of set endowed with a natural DAG structure is the set of 
communication classes of a Markov chain; minimal classes (in the usual definition) are also minimal
for the above ordering. 

\Medskip \begin{Definition}[depth of a path] Let $\alpha\ge 0$.  The $\alpha$-depth of a path $\gamma:\sigma_0\to\cdots\to
\sigma_{\ell}$ is the non-negative integer $D_{\gamma}(\alpha)$ defined by $D_{\gamma}(\alpha) = - \sum_{i=1}^{\ell} \lfloor \log_b w(\alpha)_{\sigma_{i-1}\to\sigma_i} \rfloor$. 
\end{Definition}
By construction, $D_{\gamma}(\alpha)$ is approximately equal to  minus the log-weight of the path $\gamma$. We skip the $\alpha$-dependence in the following.

\Medskip {\bf Leading simple excursions.} A  leading simple excursion from $\sigma$ to $\sigma'\not=\sigma$ 
is a  path $\gamma$ from $\sigma$ to $\sigma'$ minimizing 
$D_{\gamma}$ among all paths from $\sigma$ to $\sigma'$. The reason it is called: simple excursion, is that it cannot self-intersect; otherwise one could extract from $\gamma$ a cycle
$\gamma':\sigma_i\to\cdots\to \sigma_j=\sigma_i$, $j>i$ with $D_{\gamma'}=0$, i.e. such that 
edges $\sigma_i\to\sigma_{i+1}, \ldots, \sigma_{j-1}\to\sigma_j$ are
dominant, forming a forbidden dominant cycle.   

\Medskip Let $\{\gamma^*_{\sigma,1},\ldots,\gamma^*_{\sigma,p_{\sigma}}\}$
be the set of all leading simple excursions from $\sigma^*$ to $\sigma\in V$; ${\cal E}$ the
(non-necessarily disjoint) union
of the edges of these for every $\sigma\in V$; and $\T$ be the 
graph with vertex set $V$ and edge set ${\cal E}$.  

\begin{Lemma} \label{lem:DAG} 
 $\T$ is a DAG rooted in $\sigma^*$. Furthermore, if $\gamma,\gamma'$ are two paths from 
 $x$ to $y$, then $D_{\gamma} = D_{\gamma'}$. 
\end{Lemma} 

\noindent {\bf Proof.}
We note the following properties.  We say here that $\gamma$ is of 
higher weight than $\gamma'$ if $D_{\gamma}<D_{\gamma'}$. 

\Medskip (1) If $x\in \gamma^*_{\sigma,i}\cap \gamma^*_{\sigma',i'}$ with $(\sigma,i)\not=
(\sigma',i')$, then $D_{\sigma^* \overset{\gamma^*_{\sigma,i}}{\to} x} =  D_{\sigma^*\overset{\gamma^*_{\sigma',i'}}{\to} x}$. Namely, assume by absurd that 
$D_{\sigma^* \overset{\gamma^*_{\sigma,i}}{\to} x} > D_{\sigma^*\overset{\gamma^*_{\sigma',i'}}{\to} x}$. Then the mixed path $\sigma^*\overset{\gamma^*_{\sigma',i'}}{\to} x \overset{\gamma^*_{\sigma,i}}{\to} \sigma$ is of higher  weight than $\gamma^*_{\sigma,i}$, which
is contradictory with the fact that $\gamma^*_{\sigma,i}$ is leading. Thus one may define
unambiguously  $D_x$ as the common value of all $D_{\sigma^*\overset{\gamma^*_{\sigma,i}}{\to} x}$. 

\Medskip (2) Generalizing, let  $\gamma$ be a path from 
 $\sigma^*$ to $x$. It is obtained by concatenating pieces of leading simple excursions 
 $\gamma_{\sigma_1,i_1}\big|_{\sigma^*\to x_1},\ldots, \gamma_{\sigma_q,i_q}\big|_{x_{q-1}\to x_q}$, with $x_q=x$. By (1), $D_{x_1}= D_{\sigma^*\overset{\gamma^*_{\sigma_1,i_1}}{\to} x_1} =
 D_{\sigma^*\overset{\gamma^*_{\sigma_2,i_2}}{\to} x_1}$. Thus one may replace 
 $\gamma_{\sigma_1,i_1}\big|_{\sigma^*\to x_1} \uplus \gamma_{\sigma_2,i_2}\big|_{x_1\to x_2}$ by $\gamma_{\sigma_2,i_2}\big|_{\sigma^*\to x_2}$ without changing the depth. Continuing inductively
 along the path, we get $D_{\gamma} = D_{\sigma^*\overset{\gamma_{\sigma_q,i_q}}{\to} x} = D_x$. 
 
 \Medskip (3) Generalizing further, let $\gamma$ be a path from $x$ to $y$. As in (2), 
 it is obtained by concatenating pieces of leading simple excursions 
 $\gamma_{\sigma_1,i_1}\big|_{x_0\to x_1},\ldots, \gamma_{\sigma_q,i_q}\big|_{x_{q-1}\to x_q}$, with $x_0=x, x_q=y$. Concatenating it with $\gamma_0:=\gamma_{\sigma_1,i_1}\big|_{\sigma^*\to x_1}$, we 
 get: $D_{\gamma_0\uplus \gamma} = D_x$. Hence $D_{\gamma}= D_x - D_{\gamma_0} = D_x-D_y$ is 
 independent of the choice of path from $x$ to $y$.

\Medskip (4) Assume (by absurd) that $\T$ contains a cycle $x_1\to x_2\to\cdots x_p\to x_1$. 
Remove from the list $x_1,\ldots,x_p$ all "intermediate" species $x_i$ such that  
$x_{i-1}\to x_i$ and $x_i\to x_{i+1}$ are edges along the same leading simple excursion, we 
get (after renumbering) a sequence of paths $x_1\overset{\gamma_1}{\to} x_2,\cdots, 
x_{q-1}\overset{\gamma_{q-1}}{\to} x_q, x_q \overset{\gamma_q}{\to} x_1$, with 
$\gamma_j \subset \gamma^*_{\sigma_j,q_j}$. Since the depth is increasing along a path,
$D_{x_1}=D_{x_2}=\ldots=D_{x_1}$, meaning that $\gamma_q\circ\cdots\circ\gamma_1$ forms 
a forbidden dominant cycle.  \hfill\eop


\vskip 2cm



\section{Appendix 1. Proof of a priori bounds for $\lambda$}  \label{section:app1}


The setting here is that of \S \ref{subsection:cut-off-setting}.  We prove here a general double inequality for $\lambda^*$, see (\ref{eq:double-inequality-lambda*}) below, implying (\ref{eq:apriorilowerbound-lambda}) and (\ref{eq:aprioriupperbound-lambda}).  Choose
some uniform degradation rate $\alpha$.

\Medskip Recall from (\ref{eq:AAkk}) the identity  $|\tilde{A}^{int}_{\sigma,\sigma}| = |A_{\sigma,\sigma}| + \kappa_{\sigma}-k_{\sigma}^{ext}$.
Let $\alpha\ge 0$.  The transition kernels
\BEQ w(\alpha)_{\sigma\to\sigma'} = \frac{A_{\sigma',\sigma}}{|A_{\sigma,\sigma}|+\alpha}, \qquad 
\tilde{w}(\alpha)_{\sigma\to\sigma'} = \frac{A_{\sigma',\sigma}}{|\tilde{A}_{\sigma,\sigma}|+\alpha} \qquad \sigma'\not=\sigma
\EEQ
are associated to the defective adjoint generators  $A(\alpha):=A-\alpha\Id$, $\tilde{A}(\alpha):=\tilde{A}-\alpha\Id$.
 Let
$f(\alpha)_{\sigma\to \sigma'}$, resp. $\tilde{f}(\alpha)_{\sigma\to \sigma'}$ (see (\ref{eq:excursion})) be the weight of excursions from $\sigma$ to $\sigma'$ for these transition kernels. In particular, if it can be proven that $\max_{\sigma}(f(\alpha)_{\sigma\to
\sigma})>1$, then the network is autocatalytic, implying
the lower bound $\lambda^*>\alpha$.  Fix  a species index $\sigma^*\in \Sigma$.  A first-step analysis yields the  following linear system
for the vector $(\tilde{f}(\alpha)_{\sigma\to\sigma^*})_{\sigma}$, 
\BEQ  \begin{cases} \tilde{f}(\alpha)_{\sigma\to \sigma^*} = \sum_{\sigma'\not=\sigma} \tilde{w}(\alpha)_{\sigma\to \sigma'}
\tilde{f}(\alpha)_{\sigma'\to\sigma^*} + \tilde{w}(\alpha)_{\sigma\to \sigma^*} \qquad \sigma\not=\sigma^* \\
\tilde{f}(\alpha)_{\sigma^*\to \sigma^*} =\sum_{\sigma'\not=\sigma^*} \tilde{w}(\alpha)_{\sigma^*\to \sigma'}
\tilde{f}(\alpha)_{\sigma'\to\sigma^*} 
\end{cases}
\EEQ

In particular,
\BEQ \sum_{\sigma'\not=\sigma^*} \tilde{w}(0)_{\sigma^*\to
\sigma'} \tilde{f}(0)_{\sigma'\to\sigma^*} = 
\tilde{f}(0)_{\sigma^*\to\sigma^*}=1 \label{eq:proba-cons}
\EEQ
by probability conservation.

\Medskip
The solution $f=f(w)$ of the general system 
$\begin{cases} f_{\sigma\to\sigma^*} = \sum_{\sigma'\not=\sigma}
w_{\sigma\to\sigma'} f_{\sigma'\to\sigma^*} + w_{\sigma\to
\sigma^*} \qquad \sigma\not=\sigma^* \\
f_{\sigma^*\to\sigma^*} = \sum_{\sigma'\not=\sigma^*}
w_{\sigma^*\to\sigma'} f_{\sigma'\to\sigma^*} \end{cases}$
with $w\ge 0$  may be obtained by iterating the system of equations\\
$\begin{cases} f^{(n+1)}_{\sigma\to\sigma^*} = \sum_{\sigma'\not=\sigma}
w_{\sigma\to\sigma'} f^{(n)}_{\sigma'\to\sigma^*} + w_{\sigma
\to \sigma^*} \qquad \sigma\not=\sigma^* \\
f^{(n+1)}_{\sigma^*\to\sigma^*} = \sum_{\sigma'\not=\sigma^*}
w_{\sigma^*\to\sigma'} f^{(n)}_{\sigma'\to\sigma^*} \end{cases}$
with initial condition $f^{(0)}=0$ 
and taking the limit $n\to\infty$. Therefore it is monotonous
in $w$, i.e. $\Big(w'_{\sigma\to\sigma'}\ge w_{\sigma\to\sigma'},
\ \sigma\not=\sigma'\in \Sigma\Big)\Rightarrow \Big(f(w')_{\sigma\to\sigma^*}
\ge f(w)_{\sigma\to\sigma^*}, \ \sigma\in \Sigma\Big)$.
We let 
\BEQ m:=\min_{\sigma}| \tilde{A}_{\sigma,\sigma}|, \qquad 
M:=\max_{\sigma}| \tilde{A}_{\sigma,\sigma}| \EEQ
and $c := \alpha/m$, $C:=\alpha/M$ $(c>C>0)$. 
Then
\BEQ \frac{1}{1+c} \tilde{w}(0)_{\sigma\to\sigma'} \le  \tilde{w}(\alpha)_{\sigma\to\sigma'}= \frac{\tilde{A}_{\sigma',\sigma}}{|\tilde{A}_{\sigma,\sigma}|+\alpha} 
\le  \frac{1}{1+C} \tilde{w}(0)_{\sigma\to\sigma'}
\EEQ
hence (by monotonicity in $w$)
\BEQ  \frac{1}{1+c}
\tilde{f}(0)_{\sigma\to\sigma^*}\le \tilde{f}(\alpha)_{\sigma\to\sigma^*} \le  \frac{1}{1+C}
\tilde{f}(0)_{\sigma\to\sigma^*}, \qquad \sigma\in \Sigma
\label{eq:tildefalphatildef0}
\EEQ

Then, 
\BEQ w(\alpha)_{\sigma\to\sigma'} = \frac{A_{\sigma',\sigma}}{
|A_{\sigma,\sigma}|+\alpha} = \frac{|\tilde{A}_{\sigma,\sigma}|+
\alpha}{|A_{\sigma,\sigma}|+\alpha}\  \tilde{w}(\alpha)_{\sigma\to\sigma'} \EEQ
hence
\BEQ 1+d\le \frac{w(\alpha)_{\sigma\to\sigma'}}{\tilde{w}(\alpha)_{\sigma\to\sigma'}} \le \max_{\sigma} 
\frac{|\tilde{A}_{\sigma,\sigma}|+
\alpha}{|A_{\sigma,\sigma}|+\alpha} = 1+ \max_{\sigma}
\frac{\kappa_{\sigma}}{|A_{\sigma,\sigma}|+\alpha} \le 
1+D
\EEQ
where
\BEQ d:= \min_{\sigma} \Big(\frac{\kappa_{\sigma}-k_{\sigma}^{ext}}{|A_{\sigma,\sigma}|} \Big), \qquad 
 D:= \max_{\sigma} \Big(\frac{\kappa_{\sigma}-k_{\sigma}^{ext}}{|A_{\sigma,\sigma}|}\Big)
\EEQ
from which we get by monotonicity: 
\BEQ (1+d)\tilde{f}(\alpha)_{\sigma\to\sigma^*} \le  f(\alpha)_{\sigma\to\sigma^*} \le (1+D)\tilde{f}(\alpha)_{\sigma\to\sigma^*} \label{eq:f<=Dftilde}
\EEQ

\Bigskip
Now,
\BEQ f(\alpha)_{\sigma^*\to\sigma^*} = \sum_{\sigma'\not=\sigma^*} w(\alpha)_{\sigma^*\to \sigma'}
f(\alpha)_{\sigma'\to\sigma^*}  
\EEQ
Using
\BEQ w(\alpha)_{\sigma^*\to \sigma'}
= \tilde{w}(0)_{\sigma^*\to \sigma'} \ \times\ 
\frac{|\tilde{A}_{\sigma^*,\sigma^*}|}{|A_{\sigma^*,\sigma^*}| +\alpha}. \label{eq:wwAkappa}
\EEQ
and (\ref{eq:tildefalphatildef0}),
(\ref{eq:proba-cons}), we get
\BEA
f(\alpha)_{\sigma^*\to\sigma^*} &=& \frac{|\tilde{A}_{\sigma^*,\sigma^*}|}{|A_{\sigma^*,\sigma^*}|+\alpha}  \ \times\ \sum_{\sigma'\not=\sigma} \tilde{w}(0)_{\sigma^*\to\sigma'} f(\alpha)_{\sigma'\to\sigma^*} \\
&\ge & \frac{|\tilde{A}_{\sigma^*,\sigma^*}|}{|A_{\sigma^*,\sigma^*}|+\alpha}  \ \times\ \sum_{\sigma'\not=\sigma} \tilde{w}(0)_{\sigma^*\to\sigma'} \tilde{f}(\alpha)_{\sigma'\to\sigma^*} \nonumber\\
&\ge & (1+d)h(m,\alpha|\sigma^*)  \sum_{\sigma'\not=\sigma} \tilde{w}(0)_{\sigma^*\to\sigma'} \tilde{f}(0)_{\sigma'\to\sigma^*} \nonumber\\
&=& (1+d) h(m,\alpha|\sigma^*)
\EEA
where 
\BEQ h(\mu,\alpha|\sigma^*):= \frac{\mu}{\mu+\alpha} \frac{|\tilde{A}_{\sigma^*,\sigma^*}|}{|A_{\sigma^*,\sigma^*}| +\alpha}.
\EEQ 
Conversely, combining (\ref{eq:f<=Dftilde}),  (\ref{eq:tildefalphatildef0}) and  (\ref{eq:proba-cons}), we get
\BEA f(\alpha)_{\sigma^*\to\sigma^*} &=& \frac{|\tilde{A}_{\sigma^*,\sigma^*}|}{|\tilde{A}_{\sigma^*,\sigma^*}| - \kappa_{\sigma^*}+\alpha}  \ \times\ \sum_{\sigma'\not=\sigma} \tilde{w}(0)_{\sigma^*\to\sigma'} f(\alpha)_{\sigma'\to\sigma^*} \\
&\le & \frac{|\tilde{A}_{\sigma^*,\sigma^*}|}{|\tilde{A}_{\sigma^*,\sigma^*}| - \kappa_{\sigma^*}+\alpha}  \ \times\
\frac{1+D}{1+C}
\sum_{\sigma'\not=\sigma} \tilde{w}(0)_{\sigma^*\to\sigma'} \tilde{f}(0)_{\sigma'\to\sigma^*} \nonumber\\
& = & (1+D)h(M,\alpha|\sigma^*)
\EEA

\Medskip
If there exists any $\sigma^*$ for which $(1+d)h(m,\alpha|\sigma^*)>1$,
then $\lambda^*>\alpha$. Conversely, if  there exists any $\sigma^*$ for which $(1+D)h(M,\alpha|\sigma^*)<1$,
then $\lambda^*<\alpha$. Therefore  

\BEQ  \sup\{\alpha \ |\ \max_{\sigma^*}\  (1+d)h(m,\alpha|\sigma^*)>1\}\le \lambda^* \le \min\{\alpha \ |\ \min_{\sigma^*}\  (1+D)h(M,\alpha|\sigma^*)<1\}.
\EEQ

Now, $(1+d)h(m,\alpha|\sigma^*)>1$ if and only if 
$\alpha^2 + (|A_{\sigma^*,\sigma^*}| 
+m)\alpha - m( d |\tilde{A}_{\sigma^*,\sigma^*}|+\kappa_{\sigma^*})<0$, equivalently,  if
\BEQ \alpha<\alpha_{thr}(m|\sigma^*):= 
\frac{-x(m|\sigma^*)+\sqrt{x(m|\sigma^*)^2+y(m|\sigma^*)^2}}{2},
\label{eq:alphathrmsqrt}
\EEQ
where
\BEQ x(m|\sigma^*)=  |A_{\sigma^*,\sigma^*}|+m, \qquad y(m|\sigma^*)= \sqrt{4m( d |\tilde{A}_{\sigma^*,\sigma^*}|+\kappa_{\sigma^*} - k^{ext}_{\sigma^*})}
\EEQ
Conversely, $(1+D)h(M,\alpha|\sigma^*)<1$ if and only if
$\alpha^2 + (|A_{\sigma^*,\sigma^*}| 
+M)\alpha - M(D|\tilde{A}_{\sigma^*,\sigma^*}|+ \kappa_{\sigma^*})<0$, equivalently,  if
\BEQ \alpha>\alpha_{thr}(M|\sigma^*):= 
\frac{-x(M|\sigma^*)+\sqrt{x(M|\sigma^*)^2+y(M|\sigma^*)^2}}{2},
\label{eq:alphathrMsqrt}
\EEQ 
where
\BEQ x(M|\sigma^*)=  |A_{\sigma^*,\sigma^*}|+M, \qquad y(M|\sigma^*)= \sqrt{4M( D |\tilde{A}_{\sigma^*,\sigma^*}|+\kappa_{\sigma^*}- k^{ext}_{\sigma^*})}
\EEQ
Thus we get our key formula,
\BEQ \max_{\sigma^*}\alpha_{thr}(m|\sigma^*)\le \lambda^* \le  \min_{\sigma^*}\alpha_{thr}(M|\sigma^*).
\label{eq:double-inequality-lambda*}
\EEQ

Note that $|d|,|D|\preceq 1$ and $\kappa_{\sigma^*},k^{ext}_{\sigma^*}, |\tilde{A}_{\sigma^*,\sigma^*}| \preceq 
|A_{\sigma^*,\sigma^*}|$, hence $\frac{y(m|\sigma^*)}{x(m|\sigma^*)},\frac{y(M|\sigma^*)}{x(M|\sigma^*)}  \preceq 1$, implying by an elementary estimate $\alpha_{thr}(m|\sigma^*) \sim \frac{y^2(m|\sigma^*)}{x(m|\sigma^*)}, \ \alpha_{thr}(M|\sigma^*) \sim \frac{y^2(M|\sigma^*)}{x(M|\sigma^*)}$.


\section{Appendix 2. Doeblin's spectral gap lower bound}  \label{section:Doeb}


{\bf Proof of Lemma \ref{lem:Doeb1}.}
 We have the positive term decomposition: $ (\, ^t \tilde{\cal W}^{int})_{\sigma' ,\sigma} = \Big(\tilde{w}^{int}_{\sigma\to\sigma'} -  \min_{\sigma''} (\tilde{w}^{int}_{\sigma''\to\sigma'})
 \Big)  + \min_{\sigma''} (\tilde{w}^{int}_{\sigma''\to\sigma'})$.  Hence  $^t \tilde{\cal W}^{int}$ is a contraction in the $L^1$-norm when restricted to the zero-average subspace; namely, if $u= (u_{\sigma})_{\sigma \in {\cal V}^{int}}$ is a 
vector with $u^{av}= \frac{1}{|{\cal V}^{int}|}= \sum_{\sigma} u_{\sigma}=0$, then
\BEA \sum_{\sigma'} \Big| \sum_{\sigma} (\, ^t \tilde{\cal W}^{int})_{\sigma' ,\sigma}  u_{\sigma} \Big| &= & \sum_{\sigma'}   
\Big| \sum_{\sigma}   \Big(\tilde{w}^{int}_{\sigma\to\sigma'} - 
\min_{\sigma''} (\tilde{w}^{int}_{\sigma''\to\sigma'}) \Big) u_{\sigma}  \Big| \nonumber\\
& \le & \sum_{\sigma'}   
\ \sum_{\sigma}   \Big(\tilde{w}^{int}_{\sigma\to\sigma'} - 
\min_{\sigma''} (\tilde{w}^{int}_{\sigma''\to\sigma'}) \Big) |u_{\sigma}| \nonumber\\
&\le & (1-\rho) \sum_{\sigma} |u_{\sigma}|  \label{eq:contraction}
\EEA
The zero average condition for $u$ was used on the first line.  \hfill \eop


\Bigskip {\bf Proof of Lemma \ref{lem:Doeb2}.}
Let $f = (f^{\sigma})_{\sigma \in {\cal V}^{int}}$. 
Assume 
$u$ has zero average, $u^{av}=  0$.  Then
\BEA  |\langle \tilde{W}^{int} f, u\rangle| &=& 
|\langle \tilde{W}^{int} f,u\rangle | =  |\langle f, \, ^t \tilde{\cal W}^{int} u\rangle|
  \nonumber\\
&\le & (1-\rho)\,  ||f||^*_{\infty} \, ||u||_1.
\EEA
by Lemma \ref{lem:Doeb1}. 
\hfill \eop


\section{Appendix 3. Non-Hermitian perturbation theory}   \label{section:app3-pert}


Let  ${\cal H}(\eps) =  {\cal H}_{(0)} + \eps {\cal H}_{(1)} + O(\eps^2)$
 be a perturbation of 
a non-symmetric matrix ${\cal H}_{(0)}$ with a multiplicity one eigenvalue
$\lambda_{(0)}$. Extending the formulas of standard Rayleigh perturbation theory \cite{Coh},
 one can compute the  perturbed $\eps$-dependent eigenvalue $\lambda(\eps) \equiv \lambda_{(0)} +  \eps \lambda_{(1)} + O(\eps^2)$
 and associated eigenvector
$v(\eps) \equiv v_{(0)} + \eps v_{(1)} + O(\eps^2)$. Introduce
the $\lambda_{(0)}$-eigenvector  $u_{(0)}$ of the adjoint matrix 
$^t {\cal H}_{(0)}$; we choose its normalization by imposing $\langle 
u_{(0)}, v_{(0)}\rangle = 1$ and $\langle u_{(0)}, v_{(1)}\rangle =0$. If $v\in$ Ran$({\cal H}_{(0)}-\lambda_{(0)}) = $ Ker$(^t {\cal H}_{(0)}-\lambda_{(0)})^{\perp} =
$span$(u_{(0)})^{\perp}$, we let (by abuse of notation)
$u:= ({\cal H}_{(0)}-\lambda_{(0)})^{-1} v$ be the unique vector $u$ 
such that $({\cal H}_{(0)}-\lambda_{(0)})u = v$ and $\langle u_{(0)} | u\rangle
=0$. This defines a "pseudo-inverse" operator  $({\cal H}_{(0)}-\lambda_{(0)})^{-1} : $ span$(u_{(0)})^{\perp} \to  
$ span$(u_{(0)})^{\perp}$. 
Considering order by order the $\eps$-expansion of the eigenvalue equation 
\BEQ {\cal H}(\eps)v(\eps)=\lambda(\eps) v(\eps), \EEQ
 we get  $({\cal H}_{(0)}-\lambda_{(0)})v_{(1)} =   (\lambda_{(1)} - {\cal H}_{(1)}) v_{(0)}   +
 O(\eps^2)$.  
The l.-h.s. vanishes upon evaluation against $u_{(0)}$, yielding
\BEQ \lambda_{(1)} = \langle u_{(0)} , {\cal H}_{(1)}  v_{(0)} \rangle \EEQ
\BEQ v_{(1)} =- ({\cal H}_{(0)}-\lambda_{(0)})^{-1} ( {\cal H}_{(1)} -\lambda_{(1)} ) v_{(0)} \EEQ


\section{Appendix 4. Fixed point arguments}   \label{section:app-fixed-point}



\subsection{Rigorous proof for \S \ref{subsubsection:perturbation-Lyapunov}} \label{app: fixed1}


\noindent We rewrite the eigenvalue problem $(\, ^t {\cal W}(\lambda^*)-\Id)\pi^*=0$ as the solution 
$(\lambda^*,\pi^*)\equiv (\lambda,\pi)(\vec{\eps},\vec{k}^{ext}) \in \R\times \R^{{\cal V}^{int}}$ of the fixed-point
equation
\BEQ \vec{F}(\vec{\eps},\vec{k}^{ext};\lambda,\pi)=0 \EEQ
with $\vec{F} = \left(\begin{array}{c} (F_{\sigma'})_{\sigma'\in {\cal V}^{int}} \\
F_0 \end{array}\right)$, $F_{\sigma'}(\vec{\eps},\vec{k}^{ext};\lambda,\pi) = (\, ^t 
{\cal W}(\lambda)\pi)_{\sigma'} - \pi_{\sigma'}$, $F_0(\vec{\eps},\vec{k}^{ext};\lambda,\pi) = \sum_{\sigma'\in {\cal V}^{int}}\pi_{\sigma'} - 1$.   When $\vec{\eps}=\vec{k}^{ext}=0$, 
$(\lambda^*,\pi^*)=(0,\tilde{\pi}^{int})$. Thus, we are led to use the IFT (inverse function theorem).  We prove that its hypotheses are verified when the scale parameter $b$ is large enough.

\Medskip We use the standard reduction of the IFT to the local inversion theorem (LIT). Since there
are many scalings involved, we briefly review the main steps of one of the proofs of the LIT. Let 
$x=(\vec{\eps},\vec{k}^{ext})$, $y\equiv (\lambda,\pi- \tilde{\pi}^{int})$ and $\Phi(x,y):=
(x,\vec{F}(x,y + (0,\tilde{\pi}^{int})))$. The translation in the $y$-variable centers all 
neighborhoods around $0$. Provided $F$ is locally invertible in a neighborhood of $0$, 
$\Psi=\Phi^{-1}=(G_1,G_2)$, the solution $y(x)$ of the equation $\vec{F}(x,y + (0,\tilde{\pi}^{int}))=0$ is
 obtained as the second component $\Psi_2(x,0)$. Now, the solution of the equation $\Phi(z)=w$ is 
 a fixed point  of the mapping $T: \ z\mapsto D\Phi(0)^{-1}(D\Phi(0)z-\Phi(z)+w)$. We choose a homogeneous
norm,
\BEQ  ||((\vec{\eps},\vec{k}^{ext}),(\lambda,\pi))|| :=  ||\vec{\eps}|| + ||\pi|| +  
\tau_{{\cal G}} (||\vec{k}^{ext}|| + ||\lambda||),
\EEQ
 and let $|||\ \cdot \ |||$ be the 
associated operator norm .   Let $\rho_0>0$ such that 
\BEQ  |||D\Phi(z) - D\Phi(0)||| \le \half |||D\Phi(0)^{-1}|||^{-1}, \qquad ||z||\le \rho_0.
\label{eq:Phi-continuity}
\EEQ
Then $T$ is contracting in the ball $||z||\le \rho_0/2$ provided $||w||\le  \half|||D\Phi(0)^{-1}|||^{-1} \rho_0$, yielding a unique solution to the equation $\Phi(z)=w$. Coming back to the IFT and 
leaving out the translation in the $y$-variable, differential calculus yields $D_x y = -(D_y \vec{F})^{-1} D_x \vec{F}$. Furthermore, $D\Phi(x,y) = \left[\begin{array}{cc} I  & 0 \\ D_x\vec{F} & D_y \vec{F} 
\end{array}\right]$ has inverse $D\Phi(x,y)^{-1} = \left[\begin{array}{cc} I & 0 \\ - (D_y \vec{F})^{-1} D_x F & (D_y \vec{F})^{-1} \end{array}\right]$. Thus
\BEQ |||D \Phi(0,(0,\tilde{\pi}^{int}))^{-1}||| = O(1) + O(|||D_y \vec{F}(0,(0,\tilde{\pi}^{int}))|||^{-1}) \ (1+ O(||D_x \vec{F}(0,(0,\tilde{\pi}^{int}))||)).  \label{eq:DPhi-1}
\EEQ
We show below that the l.-h.s. of (\ref{eq:DPhi-1}) is $O(1)$, and that we can choose $\rho_0$ of order 1 such that (\ref{eq:rho}) is satisfied.

\Medskip Derivatives of $\vec{F}$ are estimated based on eq. (\ref{eq:w(lambda*)sim3}). Let 
\BEQ \vec{V}^{\eps}_{\sigma} := \frac{\partial \vec{F}}{\partial\eps_{\sigma}}, 
\qquad  \vec{V}^{ext}_{\sigma}:= \frac{\partial \vec{F}}{\partial k^{ext}_{\sigma}} 
\label{eq:fixed1-V} 
\EEQ
One has $(V^{\eps}_{\sigma})_{\sigma'} \sim \tilde{w}^{int}_{\sigma\to\sigma'} \pi_{\sigma}$,  $(\vec{V}^{ext}_{\sigma})_{\sigma'} \sim  -\tilde{w}^{int}_{\sigma\to\sigma'} \frac{\pi_{\sigma}}{k_{\sigma}}$, while the last components vanish, $(V^{\eps}_{\sigma})_0 = 
(V^{est}_{\sigma})_0=0$. It holds: $|||D_x \vec{F}|||  = O(\max_{\sigma} ||\vec{V}^{\eps}_{\sigma}||) 
+ O(\max_{\sigma}  (\frac{1}{\tau_{{\cal G}^{int}}} ||\vec{V}^{ext}_{\sigma}||)) =  O(1) $.

\BEQ M := \frac{D\vec{F}}{D(\lambda^*,\pi^*)} = \left[\begin{array}{cc} -\vec{m}_{\lambda} & 
M_{\pi} \\ 0 & ^t {\bf 1} \end{array}\right]  \label{eq:fixed1-M} \EEQ
with
\BEA && (m_{\lambda})_{\sigma'} = -\sum_{\sigma} \pi_{\sigma} \frac{\partial}{\partial\lambda}
(w(\lambda)_{\sigma\to\sigma'}) \sim  \sum_{\sigma} \pi_{\sigma} \frac{\tilde{w}^{int}_{\sigma\to\sigma'}}{k_{\sigma}}, \nonumber\\
&& \qquad 
(M_{\pi})_{\sigma',\sigma} = w(\lambda^*)_{\sigma\to\sigma'} - \del_{\sigma,\sigma'} = 
(\, ^t \tilde{\cal W}^{int} - \Id)_{\sigma',\sigma} + {\cal H}_{\sigma',\sigma}, 
\label{eq:fixed1-Mpi}
\EEA
where 
\BEQ {\cal H}_{\sigma',\sigma} \sim \tilde{w}^{int}_{\sigma\to\sigma'} (\eps_{\sigma} - \frac{k^{ext}_{\sigma}}{k_{\sigma}} - \frac{\lambda}{k_{\sigma}}) \label{eq:7H}
\EEQ
is small, $|||{\cal H}||| = O(b^{-1})$.   
Then, the IFT states that derivatives  $\vec{U}^{\eps}_{\sigma} := 
\left[\begin{array}{c} \partial_{\eps_{\sigma}} \lambda^* \\ \partial_{\eps_{\sigma}} \pi^* 
\end{array}\right], 
\ \vec{U}^{ext}_{\sigma} := 
\left[\begin{array}{c} \partial_{k^{ext}_{\sigma}} \lambda^* \\ \partial_{k^{ext}_{\sigma}} \pi^* 
\end{array}\right] $ are solutions of the linear systems 
\BEQ M\vec{U}^{\eps}_{\sigma} = - 
\vec{V}^{\eps}_{\sigma}, \ {\mathrm{resp.}}  \ M \vec{U}^{ext}_{\sigma} = - \vec{V}^{ext}_{\sigma}.
\label{eq:fixed1-MUV}
\EEQ 

\Medskip Let $(\vec{U},\vec{V})= (\vec{U}^{\eps}_{\sigma},\vec{V}^{\eps}_{\sigma})$ or $(\vec{U}^{ext}_{\sigma},\vec{V}^{ext}_{\sigma})$ indifferently. Write $\vec{U} = \left[\begin{array}{c} 
u_{\lambda} \\ u \end{array}\right], \vec{V} = \left[\begin{array}{c} v \\ 0 
\end{array}\right]$, with $u= (u_{\sigma'})_{\sigma'\in {\cal V}^{int}}, v= (v_{\sigma'})_{\sigma'\in {\cal V}^{int}}$. Then
\BEQ M\vec{U} = -\vec{V} \Longleftrightarrow \begin{cases} 
(\, ^t \tilde{\cal W}^{int} - \Id)u = -(v+ {\cal H}u) + u_{\lambda} \vec{m}_{\lambda} \\
\langle {\bf 1}, u \rangle =0 \end{cases} 
\EEQ
Let 
\BEQ \tau_{\lambda} := \sum_{\sigma'} (m_{\lambda})_{\sigma'} \sim \sum_{\sigma} \frac{\pi_{\sigma}}{k_{\sigma}} 
\sim \tau_{{\cal G}}, \EEQ
see (\ref{eq:tauGint}). This system of equations can be solved if and only if 
$\langle {\bf 1}, -(v+{\cal H}u) + u_{\lambda} \vec{m}_{\lambda}\rangle=0$, i.e. 
\BEQ \tau_{\lambda} u_{\lambda} = \langle {\bf 1}, v+ {\cal H}u \rangle. \label{eq:taulambdaulambda}
\EEQ
 In other words 
(using the pseudo-inverse $(\, ^t \tilde{\cal W}^{int} - \Id)^{-1}$, see main text) 
\BEQ u= -(\, ^t \tilde{\cal W}^{int} - \Id)^{-1}
(\Id - P)(v+{\cal H}u),  \label{eq:uv0}
\EEQ
 where $P$  is a projection operator  onto the normalized vector $\frac{\vec{m}_{\lambda}}{\tau}$, $P(w):=\langle {\bf 1},w\rangle\ \frac{\vec{m}_{\lambda}}{\tau}$.  Thus 
 $\frac{||u||}{||v||} \le ||| (\, ^t \tilde{\cal W}^{int} - \Id + (\Id-P){\cal H} )^{-1} ||| \le 
 |||(\, ^t \tilde{\cal W}^{int} - \Id)^{-1}  |||) \ \times\  ||| \Big(I + (\Id-P){\cal H} (\, ^t \tilde{\cal W}^{int} - \Id)^{-1} \Big)^{-1} ||| = O(1)$ since $|||(\, ^t \tilde{\cal W}^{int} - \Id)^{-1}  ||| = O(1)$ and $|||{\cal H}||| = O(b^{-1})=o(1)$. This and (\ref{eq:taulambdaulambda})
 imply: $|||(D_y\vec{F})^{-1}||| = |||M^{-1}|||=O(1)$.   Then continuity of $\vec{V}^{\eps}$, $\vec{V}^{ext}, \vec{m}_{\lambda}, M_{\pi}$ 
  as defined above easily imply the continuity property (\ref{eq:Phi-continuity}) for a ball
   radius $\rho_0$ of order 1.

\Medskip We finally check compatibility of the general differential formula $D_x y = -(D_y \vec{F})^{-1} D_x \vec{F}$ with the perturbative results of \S \ref{subsubsection:perturbation-Lyapunov}.  Consider first  $(\vec{U},\vec{V})= (\vec{U}^{\eps}_{\sigma},\vec{V}^{\eps}_{\sigma})$. 
To leading order, see (\ref{eq:taulambdaulambda}), $\tau_{\lambda} u_{\lambda} \sim \langle
{\bf 1},v\rangle \sim \sum_{\sigma'} \tilde{w}^{int}_{\sigma\to\sigma'} \pi_{\sigma}  
\sim \pi_{\sigma} $, so that 
\BEQ \partial_{\eps_{\sigma}} \lambda^* = u_{\lambda^*} \sim \frac{\pi^*_{\sigma}}{\tau_{\lambda^*}} \sim \frac{\tilde{\pi}^{int}_{\sigma}}{\tau_{\cal G}}, \EEQ
in agreement with (\ref{eq:lambda1-gen}). Then, 
\BEA \Big( (\, ^t \tilde{\cal W}^{int}-\Id) \partial_{\eps_{\sigma}}\pi^*\Big)_{\sigma'} &=&
\Big( (\Id - P)(v+ {\cal H}u) \Big)_{\sigma'} \sim \Big( (\Id-P)v\Big)_{\sigma'}  \nonumber\\
&\sim& \tilde{w}^{int}_{\sigma\to\sigma'} \pi^*_{\sigma} - \langle {\bf 1},
\tilde{w}^{int}_{\sigma\to\sigma'} \pi^*_{\sigma} \rangle_{\sigma'} \ \frac{(m_{\lambda^*})_{\sigma'}}{\tau_{\lambda^*}}  \nonumber\\
&\sim &  \tilde{w}^{int}_{\sigma\to\sigma'} \pi^*_{\sigma} - 
\partial_{\eps_{\sigma}}\lambda^*\ \times\  \tilde{w}^{int}_{\sigma\to\sigma'} \frac{\pi^*_{\sigma}}{k_{\sigma}}.
\EEA
This agrees to leading order with (\ref{eq:pi1}), (\ref{eq:H1}).   

\Medskip Consider next $(\vec{U},\vec{V})= (\vec{U}^{ext}_{\sigma},\vec{V}^{ext}_{\sigma})$. 
Computations are very similar.
First, $\tau_{\lambda}u_{\lambda} \sim \langle \bf{1},v\rangle \sim -\frac{\pi_{\sigma}}{k_{\sigma}}$, implying $\partial_{k^{ext}_{\sigma}}\lambda^* \sim -\frac{1}{\tau_{{\cal G}}} 
\frac{\tilde{\pi}^{int}_{\sigma}}{k_{\sigma}}$, in agreement with  (\ref{eq:lambda1-gen}). Then,
$\Big( (\, ^t \tilde{\cal W}^{int}-\Id) \partial_{k^{ext}_{\sigma}}\pi^*\Big)_{\sigma'} 
\sim \Big( (\Id-P)v\Big)_{\sigma'}  \sim -\tilde{w}^{int}_{\sigma\to\sigma'} \frac{\pi^*_{\sigma}}{k_{\sigma}} + \langle {\bf 1},
\tilde{w}^{int}_{\sigma\to\sigma'} \frac{\pi^*_{\sigma}}{k_{\sigma}} \rangle_{\sigma'} \ \frac{(m_{\lambda^*})_{\sigma'}}{\tau_{\lambda^*}} \sim  - \tilde{w}^{int}_{\sigma\to\sigma'} \frac{\pi^*_{\sigma}}{k_{\sigma}} - 
\partial_{k^{ext}_{\sigma}}\lambda^*\ \times\  \tilde{w}^{int}_{\sigma\to\sigma'} \frac{\pi^*_{\sigma}}{k_{\sigma}}$, 
which agrees once again with (\ref{eq:pi1}), (\ref{eq:H1}).


\subsection{Rigorous proof for \S \ref{subsubsection:one-renormalization-step}} \label{app:fixed2}


The procedure is the same as in \S \ref{app: fixed1}. The main difference is that we replace
the variable $\pi$ by the renormalized measure $\underline{\rho}:= z(0)^{-1} {\cal Z}(\underline{\lambda})\underline{\pi}$.  We rewrite the eigenvalue problem $(\, ^t \underline{\cal W}(\underline{\lambda})-\Id)\underline{\pi}=0$ as the solution 
$(\underline{\lambda}^*,\underline{\rho}^*)\equiv (\underline{\lambda},\underline{\rho})(\vec{\eps},\vec{k}^{ext}) \in \R\times \R^{{\cal V}^{int}}$ of the fixed-point
equation
\BEQ \vec{F}(\vec{\eps},\vec{k}^{ext};\underline{\lambda},\underline{\rho})=0 \EEQ
with $\vec{F} = \left(\begin{array}{c} (F_{\sigma'})_{\sigma'\in {\cal V}^{int}} \\
F_0 \end{array}\right)$, $F_{\sigma'}(\vec{\eps},\vec{k}^{ext};\underline{\lambda},\underline{\rho}) = ( (\, ^t 
{\cal W}(\underline{\lambda})-\Id) ({\cal Z}(\underline{\lambda}))^{-1}\underline{\rho})_{\sigma'} $, $F_0(\vec{\eps},\vec{k}^{ext};\underline{\lambda},\underline{\rho}) = z(0)\sum_{\sigma'\in {\cal V}^{int}} (({\cal Z}(\underline{\lambda}))^{-1}\underline{\rho})_{\sigma'} - 1$.   The homogeneous norm is here 
\BEQ ||((\vec{\eps},\vec{k}^{ext}),(\underline{\lambda},\underline{\pi}))|| :=  ||\vec{\eps}|| + ||\underline{\pi}|| +  
\tau^{ren} (||\vec{k}^{ext}|| + ||\underline{\lambda}||).
\EEQ
  When $\vec{\eps}=\vec{k}^{ext}=0$, 
$(\underline{\lambda},\underline{\rho})=(0,z(0)^{-1}{\cal Z}(0)\tilde{\pi}_{(0)}) =
 (0, \left(\begin{array}{c} \underline{\pi}^{\perp} \\ \underline{\pi}^{ren}\end{array}\right))$ is 
of order $O(1)$. 

\Medskip Let $\vec{V}=(\vec{V}^{\eps},\vec{V}^{ext})$ as in (\ref{eq:fixed1-V}). By definition, letting
$x = (\vec{\eps},\vec{k}^{ext})$, we get
\BEQ \vec{V} =  \Big(\partial_x \, \left(\begin{array}{cc} 0 & (-P^{\perp} \tilde{\cal W}^{int} 
(\partial_x \Del) P^{av} \ \ 0) \\ 0 & 0 \end{array}\right) + 
\left(\begin{array}{ccc} -P^{\perp} \, ^t \tilde{\cal W}^{int} (\partial_x D) P^{\perp} & 0 & 0\\
P^{av} (\partial_x D)P^{\perp} & 0 & 0  \\ (\partial_x \underline{\cal W}^{int\to ext}(\underline{\lambda}) P^{\perp} & 0 & 0 \end{array}\right) \Big) \underline{\rho} \EEQ
which is $O(1)$ as in \S \ref{app: fixed1} since $\frac{\partial}{\partial\eps_{\sigma}} d_{\sigma}(\alpha) = -1,  \tau^{ren}\frac{\partial}{\partial k^{ext}_{\sigma}} d_{\sigma}(\alpha) \sim \frac{k^{ren}_{min}}{k_{\sigma}}$ and $\frac{\partial}{\partial\eps_{\sigma}}  \underline{\cal W}^{int\to ext}(\underline{\lambda})_{\sigma\to\sigma'} \sim \tilde{\underline{w}}^{int}_{\sigma\to\sigma'} \underline{\rho}_{\sigma},  \tau^{ren}\frac{\partial}{\partial k^{ext}_{\sigma}} \underline{\cal W}^{int\to ext}(\underline{\lambda})_{\sigma\to\sigma'}\sim -\tilde{\underline{w}}^{int}_{\sigma\to\sigma'}  \underline{\rho}_{\sigma} \frac{k^{ren}_{min}}{k_{\sigma}} $ are all $O(1)$. 

Then $M:=\frac{D\vec{F}}{D(\underline{\lambda}^*,\underline{\rho}^*)}$ has the same general structure as in (\ref{eq:fixed1-M}), with 
\BEQ \vec{m}_{\underline{\lambda}} = (^t {\cal H}_{(1)})'(\lambda) \sim \left(\begin{array}{cc} 
0 & (-P^{\perp}\, ^t \tilde{\cal W}^{int} \Del'(0) P^{av} \ \   (\, ^t {\cal W}^{ext,\perp})'(\underline{\lambda})) \\ 0 & (\, ^t \underline{\cal W}^{ren})'(\underline{\lambda}) \end{array}\right) + (\, ^t  \underline{\cal W}_{(1)})'(\underline{\lambda})
\EEQ
of order $(\tau^{ren})^{-1}$ at most (by the same arguments as in (\ref{eq:Z'HHvpi}) and below), and 
\BEQ M_{\pi} = (\, ^t 
{\cal W}(\underline{\lambda})-\Id) ({\cal Z}(\underline{\lambda}))^{-1} = \, ^t {\cal H}_{(0)} + 
\, ^t {\cal H}_{(1)}(\underline{\lambda})
\EEQ
splits as in the r.-h.s. of (\ref{eq:fixed1-Mpi}) into the sum of two terms, with 
$|||(^t {\cal H}_{(0)})^{-1}|||=O(1)$, see below (\ref{eq:fixed1-tH0inv}), and 
$|||^t {\cal H}_{(1)}(\underline{\lambda})|||=o(1)$. Finally, 
 continuity of $\vec{V}^{\eps}$, $\vec{V}^{ext}, \vec{m}_{\lambda}, M_{\pi}$ 
  as defined above easily imply the continuity property (\ref{eq:Phi-continuity}) for  a ball radius of order 1.


\newpage\begin{center}  {\bf\large Index of main definitions and notations} 
\end{center}

\Bigskip

\begin{tabular}{cc}
barrier, \S \ref{subsection:merging} & circulating path, \S \ref{subsection:merging} \\
collapse, \S \ref{subsection:description} & cut-graph $G_{cut}(i)$, \S \ref{subsection:description} \\
descending sequence/path, \S \ref{subsection:merging} &  directed acyclic graph (DAG),  
\S\ref{subsection:DAG}  \\
dominant edge,  \S \ref{section:notations}  &  dominant SCC, \S \ref{section:notations}  \\
excursion, \S \ref{subsection:resolvent}  &   Green  kernel, \S \ref{subsection:heuristics}, 
\S \ref{subsubsection:full-auto} \\
hierarchical formulas for $\pi$, $v^{\dagger}$,  \S \ref{subsection:description} & 
leading path/excursion, \S \ref{subsubsection:full-non-auto} \\
leaving edge, \S  \ref{subsection:merging} & 
merging, \S \ref{subsection:merging}  \\
reentering edge, \S \ref{subsection:merging} &  
renormalization,  \S \ref{subsection:description},  \S \ref{subsubsection:one-renormalization-step} \\
resonance, \S \ref{subsection:aim}, \S \ref{subsection:multi-scale-method}, \S \ref{subsection:cycle-elementary}, \S \ref{subsection:description} &
scale parameter $b$, \S \ref{subsection:aim}, \S \ref{section:notations} 
\\ 
\\
$A$, adjoint generator, \S \ref{section:notations} & $\tilde{A}$, probability preserving version,
\S \ref{section:notations} \\
$\underline{A}^{ren}(\alpha)$, renormalized generator \S \ref{subsubsection:an-introductory-computation} & $\underline{\tilde{A}}^{ren}$,   probability preserving version  \S \ref{subsubsection:an-introductory-computation}  \\
$b$, scale parameter, \S \ref{subsection:aim}, \S \ref{section:notations} &
$f_{av}$, adjoint average, \S 
  \ref{subsubsection:perturbation-Lyapunov} \\
$G^{\lambda},G^*$, Green kernel, \S \ref{subsubsection:full-auto}  &
$k_{\sigma}$, total outgoing rate, \S \ref{section:notations} \\
 $k_{\sigma\to\sigma'}$, 
transition (kinetic) rate, \S \ref{section:notations} &
$k^{ext}_{\cal G}$, external rate, \S \ref{subsection:description} \\ 
$n_i$, scales, \S \ref{section:notations}  &               
$n(j)$, step $j$  cut-off scale, \S \ref{subsection:description} \\
  $P_{av},P^{av}, P_{\perp},P^{\perp}$, 
projection operators, \S \ref{subsubsection:perturbation-Lyapunov} &
  $\Sigma$, species set, \S \ref{section:notations} \\
   $u^{av}$, average, \S 
  \ref{subsubsection:perturbation-Lyapunov} &
 $v^*$, Lyapunov eigenvector,  \S \ref{section:notations}  \\  
 $v^{\dagger,*}$, adjoint,  \S \ref{section:notations} & 
 $v^{\dagger}$, hierarchical formula for $v^{\dagger,*}$, \S \ref{section:description} \\
$w_{\sigma\to\sigma'}, \tilde{w}_{\sigma\to\sigma'}$, transition weights, \S \ref{section:notations} & $w(\alpha)_{\sigma\to\sigma'}, \tilde{w}(\alpha)_{\sigma\to\sigma'}$,  \S \ref{section:notations} \\
${\cal W}, \tilde{\cal W}$, transition matrix, \S \ref{section:notations} & ${\cal W}(\alpha), \tilde{\cal W}(\alpha)$,  \S \ref{section:notations} \\
$Z(\eps,\alpha)_{\cal G}$, renormalization factor, \S \ref{subsection:description} &
 $Z_{G}$, renormalization weight, \S \ref{subsection:description}
\\  \\
$\alpha$, resolvent parameter, \S \ref{subsection:resolvent} &  $\alpha_{\sigma_0}$, threshold
rate, \S \ref{subsection:description}   \\
$\beta_{\sigma}$, degradation rate,  \S \ref{section:notations} & 
$\gamma$, path, \S \ref{section:notations}  \\ 
$\eps_{\sigma}$, deficiency weight, \S \ref{section:notations} &  $\bar{\eps}_{\cal G}$,
average deficiency weight, \S \ref{subsection:description} \\
$\kappa_{\cal G}$, 
$\kappa_{\sigma}$, deficiency rate, \S  \ref{section:notations} & 
$\lambda^*$, Lyapunov exponent, \S \ref{section:notations}  \\  
$\lambda_{\cal G}$, Lyapunov exponent of $\cal G$, \S \ref{subsection:description}   &  
$\tilde{\pi}$,  discrete-time stationary measure, \S \ref{section:notations} \\
 $\pi^*$, 
Lyapunov weights, \S\ref{section:notations} & 
$\pi$, hierarchical formulas for $\pi^{*}$, \S \ref{section:description}  \\
$\pi$-dominant, \S \ref{subsubsection:example2}, \S \ref{subsubsection:full-non-auto} &
$\rho$, reaction index \\ $\sigma$, species index  &
$1/\tau_{\cal G}$, characteristic rate,  \S \ref{subsection:description}  \\ 
$\Phi(\lambda^*)_{\sigma}$, excursion weight, \S \ref{subsection:resolvent} \\ \\
$\sim$, same order, \S \ref{section:notations} & $\prec,\succ$,  different orders, \S \ref{section:notations}  \\
\end{tabular}


\end{document}